\documentclass[12pt,twoside]{amsart}
\usepackage[margin=3cm]{geometry}
\usepackage[colorlinks=false]{hyperref}
\usepackage[english]{babel}
\usepackage{graphicx,titling}
\usepackage{float}
\usepackage{amsmath,amsfonts,amssymb,amsthm}
\usepackage{lipsum}
\usepackage[T1]{fontenc}
\usepackage{fourier}
\usepackage{color}
\usepackage[latin1]{inputenc}
\usepackage{esint}
\usepackage{caption}
\usepackage{ccicons}

\makeatletter
\def\blfootnote{\xdef\@thefnmark{}\@footnotetext}
\makeatother

\newcommand\ccnote{
    \blfootnote{\copyright\,\, Assaf Naor and Kevin Ren}
    \blfootnote{\ccLogo\, \ccAttribution\,\, Licensed under a \href{https://creativecommons.org/licenses/by/4.0/}{Creative Commons Attribution License (CC-BY)}.}
}

\usepackage[export]{adjustbox}
\numberwithin{equation}{section}
\usepackage{setspace}
\renewcommand{\le}{\leqslant}

\renewcommand{\ge}{\geqslant}

\renewcommand{\mathbb}{\varmathbb}
\usepackage{fancyhdr}
\newtheorem{theorem}{Theorem}[section]
\newtheorem{lemma}[theorem]{Lemma}

\newtheorem{definition}[theorem]{Definition}
\newtheorem{remark}[theorem]{Remark}
\usepackage[shortlabels]{enumitem}
\usepackage[mathscr]{euscript}
\usepackage{dutchcal}
\usepackage{upgreek}
\usepackage{comment}

\hypersetup{hyperfootnotes=false}

\allowdisplaybreaks

\makeatletter
\def\resetMathstrut@{%
  \setbox\z@\hbox{%
    \mathchardef\@tempa\mathcode`\(\relax
    \def\@tempb##1"##2##3{\the\textfont"##3\char"}%
    \expandafter\@tempb\meaning\@tempa \relax
  }%
  \ht\Mathstrutbox@1.2\ht\z@ \dp\Mathstrutbox@1.2\dp\z@
}
\makeatother

\newtheorem{defn}[theorem]{Definition}

\newtheorem{observation}[theorem]{Observation}

\newtheorem{conjecture}[theorem]{Conjecture}
\newtheorem{problem}[theorem]{Problem}
\newtheorem{notation}[theorem]{Notation}

\newtheorem{question}[theorem]{Question}

\newcommand{\net}{\mathrm{net}}
\newcommand{\C}{\mathbb{C}}

\newcommand{\cG}{\mathcal{G}}
\newcommand{\cS}{\mathscr{S}}
\newcommand{\cT}{\mathcal{T}}

\newcommand{\R}{\mathbb{R}}
\newcommand{\Z}{\mathbb{Z}}
\newcommand{\N}{\mathbb{N}}

\newcommand{\cd}{\mathcal{d}}
\DeclareMathOperator{\Lip}{Lip}
\newcommand{\norm}[1]{\| #1 \|}

\renewcommand{\le}{\leqslant}
\renewcommand{\ge}{\geqslant}

\renewcommand{\setminus}{\smallsetminus}
\renewcommand{\subset}{\subseteq}
\newcommand{\cc}{\mathsf{c}}
\newcommand{\ee}{\mathsf{e}}
\newcommand{\ii}{\mathsf{i}}
\newcommand{\dd}{\mathsf{d}}

\newcommand{\eqdef}{\stackrel{\mathrm{def}}{=}}
\newcommand{\cZ}{\mathcal{Z}}

\newcommand{\ck}{\mathcal{k}}
\newcommand{\fp}{\mathfrak{p}}

\newcommand{\SEP}{\mathsf{SEP}}
\newcommand{\PAD}{\mathsf{PAD}}

\DeclareMathOperator{\sgn}{sign}

\newcommand{\U}{\mathbb{U}}
\renewcommand{\i}{\mathsf{i}}
\newcommand{\ul}{\uplambda}

\newcommand{\cU}{\mathscr{U}}
\newcommand{\NN}{\mathcal{N}}
\newcommand{\MM}{\mathcal{M}}
\newcommand{\XX}{\mathscr{X}}
\newcommand{\YY}{\mathscr{Y}}

\newcommand{\sub}{\mathscr{C}}
\newcommand{\cD}{\mathscr{D}}
\newcommand{\bX}{\mathbf{X}}

\newcommand{\bY}{\mathbf{Y}}
\newcommand{\bZ}{\mathbf{Z}}
\DeclareMathOperator{\prob}{\mathbb{P}}
\newcommand{\e}{\varepsilon}
\newcommand{\f}{\varphi}
\renewcommand{\d}{\delta}
\newcommand{\ud}[0]{\,\mathrm{d}}
\newcommand{\sfG}{\mathsf{G}}

\newcommand{\n}{\{1,\ldots,n\}}

\newcommand{\Part}{\mathscr{P}}

\newcommand{\1}{\mathbf{1}}
\newcommand{\diam}{\mathrm{diam}}
\newcommand{\rad}{\mathrm{rad}}

\renewcommand{\nu}{\upnu}

\newcommand{\TSEP}{\widehat{\mathsf{SEP}}}

\address{Assaf Naor, Department of Mathematics, Princeton University, Princeton, NJ 08544-1000}
\email{naor@math.princeton.edu}
\medskip
\address{Kevin Ren, Department of Mathematics, Princeton University, Princeton, NJ 08544-1000}
\email{kr5621@princeton.edu}

\begin{document}

\thispagestyle{empty}

\begin{minipage}{0.28\textwidth}
\begin{figure}[H]
%\centering
\includegraphics[width=2.5cm,height=2.5cm,left]{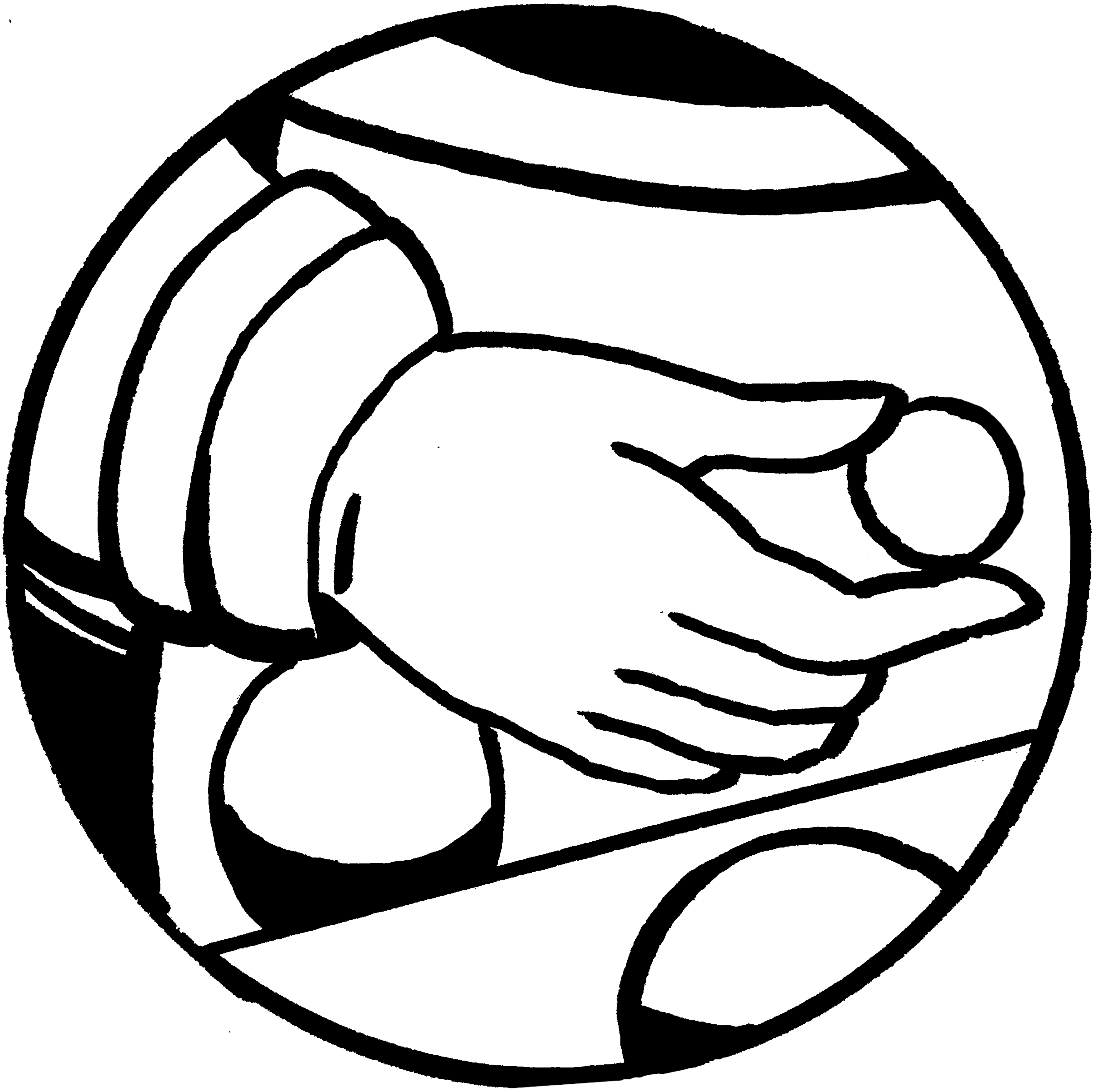}
\end{figure}
\end{minipage}%
\begin{minipage}{0.7\textwidth} 
\begin{flushright}
%% The following metadata, in particular
%% the Paper No. and the DOI will be inserted by the journal
Ars Inveniendi Analytica (2026), Paper No. 2, 60 pp.
\\
DOI 10.15781/d72xfr59
\\
ISSN: 2769-8505
\end{flushright}
\end{minipage}

\ccnote

\vspace{1cm}

%%      -------------------------------------------------------------------------------
%%      -------------------------- TITLE ----------------------------
%%      -------------------------------------------------------------------------------
%% Authors, please put here the full title of the article

\begin{center}
\begin{huge}
\textit{Euclidean embedding, randomized clustering, and Lipschitz extension for finite and doubling subsets of $L_p$ when $p>2$}

\end{huge}
\end{center}

\vspace{1cm}

%%      -------------------------------------------------------------------------------
%%      -------------------------- AUTHORS AND AFFILIATIONS ----------------------------
%%      -------------------------------------------------------------------------------
%% Authors, please put here your full names and affiliations

\begin{minipage}[t]{.45\textwidth}
\begin{center}
{\large{\bf{Assaf Naor}}} \\
\vskip0.15cm
\footnotesize{Princeton University}
\end{center}
\end{minipage}
\hfill
\begin{minipage}[t]{.45\textwidth}
\begin{center}
{\large{\bf{Kevin Ren}}} \\
\vskip0.15cm
\footnotesize{Princeton University}
\end{center}
\end{minipage}

\blfootnote{A.~N.~was supported  by NSF grant DMS-2453936, BSF grant 2018223, and a Simons Investigator award. K.~R.~was supported by an NSF GRFP fellowship. An extended abstract announcing part of this work and titled ``Optimal randomized clustering for subsets of $L_p$ when $p>2$'' appeared in the proceedings of the 37th ACM-SIAM Symposium on Discrete Algorithms (SODA26).}

\vspace{1cm}

%%% Please replace "James Mustard" below 
%%% with the name of the managing editor for your submission.
%%% If you are unsure about their identity
%%% please ask an editor-in-chief about.

\begin{center}
\noindent \em{Communicated by Larry Guth}
\end{center}
\vspace{1cm}

%%      -------------------------------------------------------------------------------
%%      -------------------------- BEGIN ABSTRACT ----------------------------
%%      -------------------------------------------------------------------------------
%% Authors, please put here the ABSTRACT and KEYBOARDS

\noindent \textbf{Abstract.} \textit{Fix $p>2$. We prove that the Euclidean distortion of every $n$-point subset of $L_p$  is $p^3(\log n)^{\frac12+o(1)}$, thus, in particular, demonstrating that all $n$-point subsets of $L_p$  exhibit  an asymptotic improvement over the $O(\log n)$ Euclidean distortion guarantee that Bourgain's embedding theorem provides for arbitrary $n$-point metric spaces.  We also prove that the separation modulus of every $n$-point subset of $ L_p$ is  $O(p^2\sqrt{\log n})$, which is sharp up to the dependence on $p$. We deduce from (a refinement of) this asymptotic evaluation of the finitary separation modulus of $ L_p$ that for any $n$-point subset $\sub$ of $ L_p$, any Banach space $\bZ$, and any $1$-Lipschitz function $f:\sub\to \bZ$, there exists a $O(p^2\sqrt{\log n})$-Lipschitz  function $F:L_p\to \bZ$ that extends $f$. We obtain analogous separation and extension statements  for doubling subsets of $L_p$.}
\vskip0.3cm

\vspace{0.5cm}

%%      -------------------------------------------------------------------------------
%%      -------------------------- BEGIN ARTICLE ----------------------------
%%      -------------------------------------------------------------------------------
%% Authors, copy the body of your paper here

\tableofcontents

\section{Introduction}

Prior to passing to detailed descriptions of concepts, results, methods  and history, we open with the following quick list that describes in broad strokes the main three outcomes of the present article which are resolutions  of  longstanding open problems in metric embedding theory:
\begin{enumerate}
\item For every  $2<p<\infty$, every $n$-point subset of $ L_p$ is shown to embed into $L_2$ with bi-Lipschitz distortion that grows to $\infty$ as $n\to \infty$ asymptotically slower $\log n$. Whether or not this holds has been a well-known open problem ever since the 1985 work~\cite{bourgain1985lipschitz} proved that any $n$-point metric space whatsoever embeds into a Hilbert space with distortion $O(\log n)$. In other words, it was unknown (and resolved herein) if---in terms of their bi-Lipschitz embeddability into $ L_2$--- anything can be said  about finite subsets of $ L_p$  beyond the fact that they are metric spaces.
\item For every $2<p<\infty$,  any $n$-point subset $\sub$ of $ L_p$, any Banach space $\bZ$, and any $1$-Lipschitz function $f:\sub\to \bZ$, there is an $L$-Lipschitz  function $F: L_p\to \bZ$ whose restriction to $\sub$ coincides with $f$, where  $L$  grows to $\infty$ as $n\to \infty$ asymptotically slower than $(\log n)/\log\log n$. Whether or not this holds has been  open ever since the 2004 work~\cite{LN04-comptes} proved that this extension statement  holds with $L=O((\log n)/\log\log n)$ and with $ L_p$  replaced by any metric space whatsoever.
\item  For every $2<p<\infty$, answering a question that was posed in the 2017 work~\cite{naor2017probabilistic}, it is proved that the largest possible separation modulus of an $n$-point subset of $ L_p$ is bounded from above and from below by positive multiples (which may depend only on $p$) of $\sqrt{\log n}$. This yields asymptotically   optimal randomized clustering of finite subsets of $ L_p$, where the quality of the clustering is measured by the ratio between the  probability that it separates points and their distance, which is  an influential  method that was introduced~\cite{bartal1996probabilistic} in the mid-1990s in the computer  science literature and has since substantially impacted both algorithm design and pure mathematics. 
\end{enumerate} 

In addition to the above listed answers to known questions, we obtain  improved Lipschitz extension and randomized clustering  results for doubling subsets of $L_p$ when $2\le p<\infty$,  which are new even when $p=2$ but (to the best of our knowledge) they have not been previously posed as open problems. Specifically, setting $C_p=p^2\sqrt{\log p}$, we prove that if $\cD\subset  L_p$ is $\ul$-doubling for some $\ul\ge 2$, then for any Banach space $\bZ$ and any $1$-Lipschitz function $f:\cD\to \bZ$, there exists a $O(C_p\sqrt{\log \ul})$-Lipschitz  function $F: L_p\to \bZ$ that extends  $f$, and furthermore the separation modulus of $\cD$ is $O(C_p\sqrt{\log \ul})$.  

All of the aforementioned results are proved through an induction on scales and localization  argument that we develop herein, in combination with a 
 novel property of the  Mazur map that we introduce.

We will next turn to a more technical (but entirely self-contained) description of the above statements:

\subsection{Bi-Lipschitz embeddings}\label{sec:bilip intro} The Euclidean distortion of a finite metric space $(\MM,d_\MM)$, which is commonly denoted $\cc_2(\MM)$ following~\cite{linial1995geometry}, is the smallest $D\ge 0$ such that there exists $f:\MM\to  L_2$ satisfying 
$$\forall x,y\in \MM,\qquad d_\MM(x,y)\le \|f(x)-f(y)\|_{\! L_2}\le Dd_\MM(x,y).$$ 
The Euclidean distortion growth $\big\{\cc_2^n(\MM)\big\}_{n=1}^\infty$ of an infinite metric space $(\MM,d_\MM)$ is:
\begin{equation}\label{eq:def finitary c2}
\forall n\in \N,\qquad \cc_2^n(\MM)\eqdef \sup_{\substack{\sub\subset \MM\\ |\sub|\le n}} \cc_2(\sub). 
\end{equation}

One says that $(\MM,d_\MM)$ has nontrivial Euclidean distortion growth if 
\begin{equation}\label{eq:def on nontrivial distortion growth}
\lim_{n\to \infty} \frac{\cc_2^n(\MM)}{\log n}=0.
\end{equation}
The above use of the word ``nontrivial'' arises from the Bourgain embedding theorem~\cite{bourgain1985lipschitz}, which asserts that $\cc_2^n(\MM)=O(\log n)$ for {\em every} metric space $(\MM,d_\MM)$.  A well-known  open question  (which has been prominently on researchers' minds ever since~\cite{bourgain1985lipschitz} appeared but to the best of our knowledge was not stated in print)   is whether $ L_p$ has nontrivial Euclidean distortion growth for every $2<p<\infty$.\footnote{For concreteness,  $L_p$ will always be the space of (equivalence classes of) $p$-integrable functions on the interval $[0,1]$, equipped with Lebesgue measure, but all of our results  hold for any $L_p(\mu)$ space; this follows formally as any separable  $L_p(\mu)$ space embeds isometrically into $L_p$ (see e.g.~\cite[Chapter~III.A]{Woj91}), though the proofs extend effortlessly to arbitrary $L_p(\mu)$ spaces.  Other such standard Banach space-theoretic notations and conventions that will be used herein are according to~\cite{LT77,LT79}. } In other words, in terms of their bi-Lipschitz  embeddability into a Hilbert space, can anything be said about finite subsets of $ L_p$  beyond merely that they are metric spaces?  Here we prove that the answer is affirmative:

\begin{theorem}\label{thm:all >2} If $2<p<\infty$, then $ L_p$ has nontrivial Euclidean distortion growth. More precisely,\footnote{We will use throughout the ensuing text the following (standard) conventions for asymptotic notation, in addition to  the usual $O(\cdot),o(\cdot),\Omega(\cdot), \Theta(\cdot)$ notation. Given $a,b>0$, by writing
$a\lesssim b$ or $b\gtrsim a$ we mean that $a\le \kappa b$ for some
universal constant $\kappa>0$, and $a\asymp b$
stands for $(a\lesssim b) \wedge  (b\lesssim a)$.  When we will need to allow for dependence on parameters, we will indicate it by subscripts. For example, in the presence of auxiliary objects $q,U,\phi$, the notation $a\lesssim_{q,U,\phi} b$ means that $a\le \kappa(q,U,\phi)b$, where $\kappa(q,U,\phi)>0$ may depend only on $q,U,\phi$, and similarly for the notations $a\gtrsim_{q,U,\phi} b$ and $a\asymp_{q,U,\phi} b$. Also, in what follows when expressions like, say, $\sqrt{\log \log\log n}$ appear for some integer $n$, it will be assumed tacitly that $n$ is greater than a sufficiently large universal constant, so that they make sense (thus, $n>15$ in the above example). In all such occurrences, the corresponding statement will be self-evident for smaller values of $n\in \N$ (by suitably adjusting an implicit universal constant factor). }
\begin{comment}
$$
\cc_2^n ( L_p)\lesssim \left\{\begin{array}{ll}\sqrt{\log n}\log\log n &\mathrm{if}\quad  2<p\le 3,\\
    (\log n)^{\frac{p}{2}-1}\log\log n & \mathrm{if}\quad 3<p< 4. \end{array}\right.
$$
\end{comment}
\begin{equation}\label{eq:p cubed version}
\forall n\in \{3,4,\ldots\},\qquad \cc_2^n( L_p)\lesssim p^3\! {\textstyle\sqrt{\log n}}\log\log n. 
\end{equation}
\end{theorem}

If $1\le p\le 2$, then $\cc_2^n( L_p)\lesssim \sqrt{\log n}$ for every integer $n\ge 2$ by~\cite{ANR24}. In fact, together with the matching lower bound from~\cite{Enf69}, it follows from~\cite{ANR24} that $\cc_2^n( L_1)\asymp \sqrt{\log n}$; at present, $ L_1$ is the only classical  Banach space for which the Euclidean distortion growth has been computed up to universal constant factors, other than the (trivial) case  $\cc_2^n( L_2)=1$ and the (substantial) case $\cc_2^n( L_\infty)\asymp \log n$ (the latter statement   is a combination of~\cite{bourgain1985lipschitz} and~\cite{linial1995geometry,AR98} for the upper and lower bounds on $\cc_2^n( L_\infty)$, respectively). 

If $p>2$, then prior to Theorem~\ref{thm:all >2} it was known that $ L_p$ has nontrivial Euclidean distortion growth only when $2< p<3\sqrt{e}\in [4.94,4.95]$. Specifically, it was announced in Section~5.3 of the 2014 arXiv posting~\cite{BG14} that $ L_p$ has nontrivial Euclidean distortion growth when $2<p<4$. We warn that~\cite{BG14} was replaced by a version in which all discussion of distortion growth was removed (and, to the best of our knowledge what was removed from~\cite{BG14} did not appear elsewhere). We do not know why the authors of~\cite{BG14}  chose to remove that material, as it seems to us that the sketch of proof in~\cite{BG14} is sound, albeit with missing steps.  In our initial arXiv posting~\cite{NR25-v1}, which is superseded by the present article, we succeeded to prove that $ L_p$ has nontrivial Euclidean distortion growth only when $2<p<4$, while being unaware of~\cite{BG14}  (the bound that we obtained in~\cite{NR25-v1} is stronger than what is implied by~\cite{BG14}).  The subsequent work~\cite{KPS25} used a clever iterative application of our main embedding theorem in~\cite{NR25-v1} to extend to $4\le p<3\sqrt{e}$ the range for which $L_p$ has nontrivial Euclidean distortion growth, but that method of~\cite{KPS25} breaks down when $p\ge 3\sqrt{e}$. From the quantitative perspective, the bound~\eqref{eq:p cubed version} of Theorem~\ref{thm:all >2} is the best-known such bound for all fixed $p>2$.

\begin{conjecture} {The bound~\eqref{eq:p cubed version} gives $\cc_2^n( L_p)=o(\log n)$ when $p=o\big(\sqrt[6]{\log n}/\sqrt[3]{\log\log n}\big)$. If $p\gtrsim \log n$, then $\ell_{\!\!\infty}^n$ embeds with $O(1)$ distortion into $ L_p$, so  $\cc_2^n( L_p)\asymp \log n$ as all $n$-point metric spaces are isometric to a subset of  $\ell_{\!\!\infty}^n$. We conjecture that $\cc_2^n( L_p) =o(\log n)$ in the remaining range $\sqrt[6]{\log n}/\sqrt[3]{\log\log n}\lesssim p=o(\log n)$.  }
\end{conjecture}

\subsubsection{Toward a sharp Euclidean distortion growth rate for \texorpdfstring{$L_p$}{Lp}}\label{sec:lewis} Now that it is established that $ L_p$ has nontrivial Euclidean distortion growth for every $1\le p<\infty$, one can turn attention to determining the growth rate of $\cc_2^n( L_p)$ as $n\to \infty$, which is currently known only for $p\in \{1,2\}$. This seems to be a highly nontrivial matter. The Lewis theorem~\cite{Lew78} asserts that $\cc_2(\bX)\le k^{|1/2-1/p|}$ for every $k$-dimensional subspace $\bX$ of $ L_p$ (this upper bound holds as equality when $\bX=\ell_{\!\! p }^k$; see e.g.~\cite{Tom89,JS01}). In accordance  with the longstanding Ribe research program~\cite{Bourgain-superreflexivity,Kal08,Naor-Ribe,Ball-Ribe,Ostrovskii-book,Ost16,God17,Nao18}), this naturally leads to the following  open question:

\begin{question}[Nonlinear Lewis problem]\label{Q:lewis} Is it true that for every $1\le p\le \infty$ we have 
\begin{equation}\label{eq:lewis conjectrure} 
\forall n\in \{2,3,\ldots\},\qquad \cc_2^n( L_p)\asymp_p (\log n)^{\left|\frac12-\frac{1}{p}\right|}.
\end{equation}
In fact, it is not even known whether~\eqref{eq:lewis conjectrure} holds for {\em any} fixed $p\in [1,\infty]\setminus \{1,2,\infty\}$.  
\end{question}

We expect  that answering Question~\ref{Q:lewis} would be very challenging (based in part on how the only known nontrivial cases $p\in \{1,\infty\}$ of Question~\ref{Q:lewis} were resolved), especially if the answer is positive, which would likely be a major achievement that would entail introducing a significant new idea/embedding method. 

The challenge within Question~\ref{Q:lewis}  is to prove the upper bound on $\cc_2^n( L_p)$ in~\eqref{eq:lewis conjectrure},  as the lower bound 
\begin{equation}\label{eq:lower bound on p growth rate}
\forall n\in \{2,3,\ldots,\},\qquad \cc_2^n( L_p)\gtrsim (\log n)^{\left|\frac12-\frac{1}{p}\right|}.
\end{equation}
follows from known results (observe that thanks to~\eqref{eq:lower bound on p growth rate}  our new upper bound~\eqref{eq:p cubed version} is optimal up to lower order factors whenever $p=(\log n)^{o(1)}\to \infty$). To justify~\eqref{eq:lower bound on p growth rate}, if $1\le p<2$, then for each $k\in \N$ consider  the $n$-point hypercube $\{0,1\}^k$ as a subset of $\ell_{\!\! p }^k$, thus $n=|\{0,1\}^n|=2^k$.  By~\cite{Enf69} its Euclidean distortion is 
$$
\cc_2\big(\{0,1\}^k,\|\cdot\|_{\!\ell_{\!\!p }^k}\big) =k^{\frac{1}{p}-\frac12}\asymp(\log n)^{\frac{1}{p}-\frac12},
$$ 
which shows that~\eqref{eq:lower bound on p growth rate} holds when $n$ is a power of $2$; it is straightforward to deduce from this that~\eqref{eq:lower bound on p growth rate} holds for every integer $n\ge 2$ (e.g., by augmenting the hypercube with auxiliary very distant points  that form an equilateral metric). If $p\ge 2$, then by~\cite{NR03} there is an $n$-vertex connected planar graph $\sfG_n$, equipped with its shortest-path metric, whose Euclidean distortion satisfies $\cc_2(\sfG_n)\gtrsim \sqrt{\log n}$. As $p\ge 2$,  the proof of the main result of~\cite{Rao99} (see also the exposition in~\cite{GKL03, KLMN05}) shows that there is an embedding $f$ of $\sfG_n$ (indeed, of any connected planar graph) into $ L_p$ whose bi-Lipschitz distortion is at most a universal constant multiple of $(\log n)^{1/p}$; the image $f(\sfG_n)\subset  L_p$ exhibits the validity of~\eqref{eq:lower bound on p growth rate} in the remaining range $p\ge 2$.

The lower bound~\eqref{eq:lower bound on p growth rate} is only part of the picture, as there {\em must} be some dependence on $p$ in~\eqref{eq:lewis conjectrure} that diverges  as $p\to \infty$, which demonstrates that there is a qualitative difference between the nonlinear setting of Question~\ref{Q:lewis} and its aforementioned linear counterpart~\cite{Lew78}.  Indeed, consider any $n$-point metric space $\MM_n$ with $\cc_2(\MM_n)\gtrsim \log n$, which exists thanks to~\cite{linial1995geometry,AR98}. Because $\MM_n$ embeds isometrically into $\ell_{\!\!\infty}^n$, while $\ell_{\!\!\infty}^n$ embeds with $O(1)$ distortion into $L_{\Omega(\log n)}$, it follows that for $p\asymp \log n$ we necessarily  have
\begin{equation}\label{eq:sqrt p lower}
\cc_2^n( L_p)\gtrsim \textstyle{\sqrt{\min\{p,\log n\}}}(\log n)^{\frac12-\frac{1}{p}}. 
\end{equation}
This suggests (but does not prove!) that perhaps~\eqref{eq:sqrt p lower} holds for all $p>2$. If so, then it would follow that  a power-type dependence on $p$ in our new bound~\eqref{eq:p cubed version} is unavoidable if $\alpha(\log\log n)/\log\log\log n\le p\lesssim \log n$ for some fixed $\alpha>2$, e.g., combining~\eqref{eq:p cubed version} and~\eqref{eq:sqrt p lower} gives $\sqrt{(\log n)\log \log n}\lesssim \cc_2^n(L_{\log\log n})\lesssim \sqrt{\log n} (\log\log n)^3$.

\subsection{Lipschitz extension}\label{sec:lip ext} Given a (source) metric space $(\MM,d_\MM)$, a subset $\sub$ of $\MM$ such that $|\sub|\ge 2$, and a (target) metric space $(\cT,d_\cT)$, one denotes (following~\cite{MAT90}) by $\ee(\MM,\sub;\cT)\in [1,\infty]$ the infimum over those $K>0$ such that for every function $f:\sub\to \cT$ there exists a function $F:\MM\to \cT$ whose restriction to $\sub$ coincides with $f$,  and the Lipschitz constant of $F$ satisfies:
\begin{equation}\label{eq:extended F is Lip}
\|F\|_{\Lip(\MM;\cT)}\le K\|f\|_{\Lip(\sub;\cT)}, 
\end{equation}
where~\eqref{eq:extended F is Lip} uses the following notation for Lipschitz constants, which will occur throughout what follows:
$$
\|f\|_{\Lip(\sub;\cT)}\eqdef \sup_{\substack{x,y\in \sub\\x\neq y}}\frac{d_\cT\left(f(x),f(y)\right)}{d_\MM(x,y)}.
$$
The supremum of $\ee(\MM,\sub;\cT)$ over all the  subsets  $\sub$ of  $\MM$ containing at least two points is denoted $\ee(\MM;\cT)$. 

For each integer  $n\ge 2$ one defines $\ee_n(\MM)$ to be the supremum of  $\ee(\MM,\sub;\bZ)$ over all $\sub\subset \MM$ satisfying $2\le |\sub|\le n$, and all Banach spaces $(\bZ,\|\cdot\|_\bZ)$.
%One defines $\ee(\MM;\cZ)$ to be the supremum of $\ee(\MM,\sub;\cZ)$ over all  $\sub\subset \MM$ with $|\sub|\ge 2$. For each $n\in \{2,3,\ldots\}$, one %analogously lets $\ee^n(\MM;\cZ)$ be the supremum of  $\ee(\MM,\sub;\cZ)$ over all $\sub\subset \MM$ satisfying $2\le |\sub|\le n$. 
By~\cite[Theorem~1.10]{LN05}, every metric space $(\MM,d_\MM)$ satisfies:
\begin{equation}\label{eq:ae bound quote}
\forall n\in \{3,4,\ldots\},\qquad \ee_n(\MM)\lesssim \frac{\log n}{\log \log n}. 
\end{equation}
\begin{remark}\label{rem:superscript version of finatary extension} {\em Given a metric space $(\MM,d_\MM)$, its {\em Lipschitz extension modulus} $\ee(\MM)$ is defines to be  the supremum of $\ee(\MM,\sub;\bZ)$ over all  $\sub\subset \MM$ with $|\sub|\ge 2$ and all Banach spaces $(\bZ,\|\cdot\|_\bZ)$.  When the (typical) situation $\ee(\MM)=\infty$ occurs, one can consider analogously to~\eqref{eq:def finitary c2} the finitary invariant   $\ee^n(\MM)$ which is defined for each   $n\in \{2,3,\ldots\}$ to be the supremum of $\ee(\sub)$ over all subsets $\sub\subset \MM$ with $2\le |\sub|\le n$.\footnote{More generally, we will maintain the following notational convention. Given an invariant  $\ii(\MM)\in \R\cup\{\infty\}$ of metric spaces $(\MM,d_\MM)$, for each $n\in \{2,3,\ldots\}$ the superscript notation $\ii^n(\MM)$ is reserved  for the supremum of $\ii(\sub)$ over $\sub\subset \MM$ with $2\le |\sub|\le n$. In contrast, the subscript notation $\ii_n(\MM)$ is used less consistently in the literature to denote more subtle finitary invariants that still consider arbitrary subsets $\sub$ of $\MM$ with $2\le |\sub|\le n$, but do not depend only on the intrinsic geometry of $(\sub,d_\MM)$.} It is important to stress  the difference between  $\ee_n(\MM)$ and $\ee^n(\MM)$, namely, $\ee_n(\MM)$ measures the extent to which for any $\sub\subset \MM$ with $2\le |\sub|\le n$, any  Lipschitz function from $\sub$ to a Banach space can be extended to a Lipschitz function on  $\MM$, while $\ee^n(\MM)$ disregards how $\sub$ is situated in the super-space $\MM$, asking for such extensions to $\sub$ of Banach space-valued functions from arbitrary subsets of $\sub$. While  $\ee_n(\MM)$ was studied in the literature for a long time (starting with~\cite{JLS86},  inspired by~\cite{Gru60,Lin64}),  to the best of our knowledge    $\ee^n(\MM)$ was not considered before. Nevertheless,  $\ee^n(\MM)$ has a key role in our investigations herein; see Section~\ref{sec:sep ext to embed in intro}. }
\end{remark}

By~\cite[Theorem~1.12]{LN05} for every $1<p\le 2$  we have: 
\begin{equation}\label{eq:quote LN p<2}
\forall n\in \{2,3,\ldots,\},\qquad \ee_n( L_p)\lesssim_p (\log n)^{\frac{1}{p}}. 
\end{equation}
The following Lipschitz extension theorem treats the analog of~\eqref{eq:quote LN p<2} in the previously unknown range $p>2$: 
\begin{theorem}\label{thm:new ext intro} If $2<p<\infty$ and $n\in \{2,3,\ldots,\}$, then $$\ee_n( L_p)\lesssim p^2{\textstyle \sqrt{\log n}}.$$
\end{theorem}
Theorem~\ref{thm:new ext intro} answers the natural question that was left open by~\cite{LN05} (and stated elsewhere, e.g.~\cite{NR17}) whether for every (indeed, any) fixed $2<p<\infty$ we have 
\begin{equation}\label{eq:extension growth}
 \lim_{n\to \infty} \ee_n( L_p) \frac{\log \log n}{\log n}=0. 
\end{equation}
In other words, prior to Theorem~\ref{thm:new ext intro} it was not even known if one could provide an upper bound on $\ee_n( L_p)$ that is asymptotically better than merely using the fact that $ L_p$ is a metric space through~\eqref{eq:ae bound quote}. 

\begin{remark}\label{rem:extension growth}{\em We hesitate to call~\eqref{eq:extension growth} ``nontrivial extension growth'' for $ L_p$ partly due to  the distinction that is noted in Remark~\ref{rem:superscript version of finatary extension}, and mainly because in contrast to the analogous setting in~\eqref{eq:def on nontrivial distortion growth}, it is unknown if~\eqref{eq:ae bound quote} is asymptotically sharp, nor is it conjectured that this is the case (though, it could very well be so). Determining the largest possible growth rate as $n\to \infty$ of $\ee_n(\MM)$ over all metric spaces $\MM$  is a major open question for which the currently best-known lower bound~\cite{NR17} is that $\ee_n(\MM)$ must sometimes be at least a positive universal constant multiple of  $\sqrt{\log n}$. In other words, the currently best-known bounds on $\ee_n( L_\infty)$ are $\sqrt{\log n}\lesssim \ee_n(L_{
\infty})\lesssim (\log n)/\log\log n$. For finite $p$, at present there is insufficient information for there to be a widely accepted conjecture what the sharp asymptotic dependence as $n\to \infty$ should be in Theorem~\ref{thm:new ext intro}, as well as in the estimate~\eqref{eq:quote LN p<2} of~\cite{LN05}. The best-available bounds in this context are presented in Section~\ref{sec:finitary impossibility} below, from which it follows that if $p$ is allowed to depend on $n$ so that it tends to $\infty$ at a  $(\log n)^{o(1)}$ rate, then Theorem~\ref{thm:new ext intro} is sharp up to lower order factors; e.g. we now know that:  
$${\textstyle \sqrt{\log n}}\lesssim \ee_n(L_{\log\log n})\lesssim {\textstyle \sqrt{\log n}}(\log\log n)^2.$$ Conceivably $\ee_n( L_p)\asymp_p \sqrt{\log n}$ or even $\ee_n( L_p)\asymp \sqrt{\log n}$ for all $p\ge 2$, but proving this would be a spectacular  achievement; in particular, the former statement would settle the Hilbertian case $p=2$ and the latter statement would also settle the aforementioned case $p=\infty$ of extension from finite subsets of general metric spaces, while both of these remain tantalizingly unknown despite major efforts over many years.  }
\end{remark}

The implicit constant in~\eqref{eq:quote LN p<2} that the proof in~\cite{LN05} provides tends to $\infty$  as $p\to 1^+$; this seems inherent to the currently available approach due to its reliance on~\cite{MP84}.\footnote{The best available bound as $p\to 1^+$ on the implicit constant in~\eqref{eq:quote LN p<2} can be deduced from the proof  in~\cite[Section~4.2]{naor2024extension}, which implies that one can take it to be at most a universal constant multiple of $1/(p-1)$.} Thus the case $p=1$ remains the last holdout for the question of improving~\eqref{eq:ae bound quote} when $\MM$ is $ L_p$ for some fixed $1\le p <\infty$.\footnote{The case $p=\infty$ is also open, as we recalled in Remark~\ref{rem:extension growth}, since it corresponds to arbitrary metric spaces $(\MM,d_\MM)$.} 

\begin{problem}\label{prob:l1 ext} Determine whether for $n\in \{3,4,\ldots\}$ we have  $\ee_n( L_1)=o\left(\frac{\log n}{\log \log n}\right)$. 
\end{problem}

It is worthwhile to recall in the context of Problem~\ref{prob:l1 ext} that by~\cite{MN06}  $\ee_n(\MM;\bZ)\lesssim_\bZ{\textstyle \sqrt{\log n}}$ for every metric space $(\MM,d_\MM)$ and every Banach space  $(\bZ, \|\cdot\|_\bZ)$ that has an equivalent norm whose modulus of uniform convexity has power type $2$ (see e.g.~\cite{BCL94} for background on this class of spaces); the Hilbertian special case $\bZ= L_2$ of this statement is the famous Johnson--Lindenstrauss extension theorem~\cite{JL82}.

A positive answer to Problem~\ref{prob:l1 ext} would have significant algorithmic consequences by improving over the best-known upper bound on the existence of vertex cut sparsifiers of weighted graphs~\cite{Moi09}, which has stood for a long time despite substantial efforts to improve it in the computer science literature (see e.g.~\cite{CLLM10,EGKRTT14,MM16}). The connection between Lipschitz extension and graph sparsification that yields this potential application of Problem~\ref{prob:l1 ext}  was discovered in~\cite{MM16}; we omit the details and the definition of the relevant  sparsification notion as they are covered thoroughly in~\cite{MM16}, as well as in e.g.~\cite[Section~1.3.3]{ANR24}.

\begin{remark} {\em Even though it is of secondary importance  given the current state of  knowledge, we note that if one allows $p\ge 2$ to depend on $n$, then Theorem~\ref{thm:new ext intro} improves over the general estimate~\eqref{eq:ae bound quote} if and only if $p=o(\sqrt[4]{\log n}/\sqrt{\log \log n})$. Understanding this for $p\gtrsim \sqrt[4]{\log n}/\sqrt{\log \log n}$ is open. For that matter, by examining its behavior as $p\to 1^+$, one sees that the bound~\eqref{eq:quote LN p<2} of~\cite{LN05}  improves over~\eqref{eq:ae bound quote} if and only if  
$$(p-1)\log\log n=2(\log\log \log n)-\log\log\log\log n+O(1);$$ understanding this for other $1<p<2$ is open.}
\end{remark}

\subsubsection{Known Lipschitz extension impossibility results}\label{sec:finitary impossibility}\hfill\break
Unlike the discussion in Section~\ref{sec:lewis}, there is no widely accepted conjectural growth rate (as $n\to \infty$) for $\ee_n( L_p)$ when  $1\le p\le \infty$ is fixed; this growth rate is not even known in the Hilbertian setting $p=2$. A na\"ive appeal to~\cite{Lew78} from the perspective of the Ribe program leads to the prediction  $\ee_n( L_p)\asymp_p (\log n)^{|1/2-1/p|}$. This was indeed a possibility for quite some time (see e.g.~the questions posed in~\cite{JLS86,Bal92} for $p=2$),  until it was shown in~\cite{Nao01} to fail even when $p=2$. 

The best-known lower bound in the range $2\le p\le 4$ is $\ee_n( L_p)\gtrsim \sqrt[4]{(\log n)/\log\log n}$. In fact, for every infinite dimensional Banach space $\bX$ we have $\ee_n(\bX)\gtrsim \sqrt[4]{(\log n)/\log\log n}$  because by~\cite{MN13-ext} (see also~\cite{Nao01}),
\begin{equation}\label{eq:l2 lower ext}
\forall m\in \N,\qquad \ee_{m^{O(m)}}\!(\ell_{\!\!2}^m)\gtrsim \sqrt[4]{m},
\end{equation}
while for all $m\in \N$ by Dvoretzky's theorem~\cite{Dvo60}  $\ell_{\!\!2 }^m$ is $O(1)$-isomorphic to a subspace of $\bX$. 

For the remaining values of $p$, the best-known lower bounds are $\ee_n(L_{\infty})\gtrsim  \sqrt{\log n}$ by~\cite{NR17}, and $\ee_n( L_p)\gtrsim ((\log n)/\log\log n)^{1/2-1/p}$ for $p\in [1,2)\cup (2,\infty)$; if $1\le p <2$, then this coincides with~\cite[Theorem~3]{JL82}, and if $p>2$, then it also follows from the reasoning in~\cite{JL82}, though it is not stated there explicitly. Indeed, by~\cite{Sob41} for every $1\le p\le \infty$ and every  $m\in \N$ there is a linear subspace $\bY=\bY_p$ of $\ell_{\!\! p }^m$ such that $\|\mathsf{Proj}\|_{\ell_{\!\! p }^m\to \ell_{\!\! p }^m}\gtrsim m^{|1/2-1/p|}$ for every projection $\mathsf{Proj}$ from $\ell_{\!\! p }^m$ onto $\bY$.  By~\cite{Lin64}, this implies: 
\begin{equation}\label{eq:no extension to Y}
\ee(\ell_{\!p}^m,\bY;\bY)\gtrsim m^{\left|\frac12-\frac{1}{p}\right|}.
\end{equation}
The discretization method of~\cite{JL82} deduces from~\eqref{eq:no extension to Y} that there is a subset $\NN=\NN_p$ of $\ell_{\!\! p }^m$ with $|\NN|=m^{O(m)}$ (specifically, $\NN$ is a $(1/m^{O(1)}$)-net of the unit sphere of $\bY$) such that $\ee(\ell_{\!\! p }^m,\NN;\bY)\gtrsim m^{|1/2-1/p|}$.

\subsubsection{Lipschitz extension from doubling subsets}\label{sec:doubling}\hfill\break
Given $\ul\in \{2,3,\ldots\}$, a  metric space $(\MM,d_\MM)$ is $\ul$-doubling if for every $r\ge 0$ and every $x\in \MM$ there exist  $y_1,\ldots,y_\ul\in \MM$ with $B_\MM(x,2r)\subset B_\MM(y_1,r)\cup\ldots\cup B_\MM(y_\ul,r)$.\footnote{The notion of $\ul$-doubling  when $\ul\ge 2$ is not necessarily an integer coincides with $\lfloor \ul\rfloor$-doubling.} Here, as well as throughout the ensuing discussion,  
$B_\MM(u,\rho)=\{v\in \MM:\ d_\MM(u,v)\le \rho\}$ denotes  the {\em closed} $d_\MM$-ball centered at a point $u\in \MM$ of radius $\rho\ge 0$.

By~\cite[Theorem~1.6]{LN05}, for every metric space $(\MM,d_\MM)$, for every $\ul$-doubling subset $\cD$ of $\MM$, and for every Banach space $(\bZ,\|\cdot\|_\bZ)$ we have 
\begin{equation}\label{eq:LN doubling case}
\ee(\MM,\cD;\bZ)\lesssim \log \ul.
\end{equation}
It is an important  open problem to determine if in the above stated generality the  right hand side of~\eqref{eq:LN doubling case} can be reduced to $o(\log \ul)$ as $\ul\to \infty$. It is also natural to investigate if such an  improvement could be achieved for specific metric spaces $\MM$, though to the best of our knowledge this has not been previously posed as an open question. In particular, it was not  known if~\eqref{eq:LN doubling case} could be improved when $\MM$ is a Hilbert space; this is answered as the special case $p=2$  of the following theorem:

\begin{theorem}\label{thm:doubling ext} If $2\le p<\infty$, then the following estimate holds for every $\ul\in \{2,3,\ldots\}$, for every subset $\cD$ of $ L_p$ that is $\ul$-doubling, and for every Banach space $(\bZ,\|\cdot\|_\bZ)$:
\begin{equation}\label{eq:doubling ext}
\ee( L_p,\mathscr{D};\bZ)\lesssim \Big(p^2{\textstyle\sqrt{\log p}}\Big){\textstyle\sqrt{\log \ul}}.
\end{equation}
\end{theorem}

We expect that the right hand side of~\eqref{eq:LN doubling case}  could be improved to $o(\log \ul)$ also in the remaining range $1\le p<2$, i.e.,  for every $\ul\in \{2,3,\ldots\}$, for every subset $\cD$ of $ L_p$ that is $\ul$-doubling, and for every Banach space $(\bZ,\|\cdot\|_\bZ)$, we conjecture that $\ee( L_p,\cD;\bZ)= o(\log \ul)$. The case $p=1$ here is most tenuous in terms of available methods and evidence. When $1<p<2$, an attempt to combine our approach herein with the proof of~\eqref{eq:quote LN p<2} in~\cite{LN05} leads to probabilistic issues that we currently do not know how to address but they could be quite tractable. In terms of the  best-known lower bounds in the context of Theorem~\ref{thm:doubling ext}, as every finite dimensional normed space $\bY$ is $\ul$-doubling for $\ul=e^{O(\dim(\bY))}$, it follows from~\eqref{eq:l2 lower ext} and~\eqref{eq:no extension to Y} that for every $\ul\in\{2,3,\ldots\}$ and every $1\le p\le \infty$  there exists a $\ul$-doubling subset $\cD=\cD_p$ of $ L_p$ and a Banach space $\bZ=\bZ_p$ for which $\ee( L_p,\mathscr{D};\bZ)\gtrsim \sqrt[4]{\log \ul}$ if $2\le p\le 4$ and $\ee( L_p,\mathscr{D};\bZ)\gtrsim (\log \ul)^{|1/2-1/p|}$ if $p\in [1,2)\cup (4,\infty]$.

\subsection{Randomized clustering}\label{sec:sep intro} Given a metric space $(\MM,d_\MM)$ and a partition $\Part$ of $\MM$, denote for each $x\in \MM$ the element of $\Part$ containing $x$ by $\Part(x)$. For $\Delta>0$, a partition $\Part$ of $\MM$ is said to be $\Delta$-bounded if $\diam_\MM(\Part(x))\le \Delta$ for every $x\in \MM$. Here, as well as  throughout the ensuing discussion, the $d_\MM$-diameter of  $\emptyset\neq \sub\subset\MM$ will be  denoted   $\diam_\MM(\sub)=\sup_{x,y\in \sub} d_\MM(x,y)\in [0,\infty]$.

Following~\cite{naor2024extension}, one says that $\Part=\{\Gamma^1,\Gamma^2,\ldots\}$ is  a  {\em random partition} of a metric space $(\MM,d_\MM)$ if there exists a probability space $(\Omega,\prob)$  such that $\Gamma^1:\Omega\to 2^\MM,\Gamma^2:\Omega\to 2^\MM,\ldots$ is a sequence of set-valued mappings that are strongly measurable,\footnote{This notion of measurability of set-valued functions is called here ``strongly measurable'' even though parts of the literature calls it more simply  ``measurable'' to distinguish it from the notion of a ``weakly measurable'' set-valued function which is also commonly used in the literature but is  {\em not} what we need herein. See~\cite{Him75} for a treatment of these classical concepts.} namely, for every $i\in \N$ and every {\em closed} subset $E$ of $\MM$ the set $\{\omega\in \Omega:\ E \cap \Gamma^i(\omega)\neq \emptyset\}$ is $\prob$-measurable, and   if we write $\Part^\omega=\{\Gamma^k(\omega):\ k\in \N\}$ for each $\omega\in \Omega$, then the mapping $\omega\mapsto \Part^\omega$ takes values in partitions of $\MM$. Note that we are formally discussing here   random {\em ordered} partitions of $\MM$ into countably many clusters, but this nuance will not have a role in what follows (it is important only for some of the results from the literature that we  will need to quote). Measurability is not relevant when $\MM$ is finite, which is the setting that we will initially discuss  below, but we will quickly need to  also  treat random partitions of infinite spaces, at which point we will verify measurability as required.

Given $\Delta>0$, one says  that $\Part$ is a $\Delta$-bounded random partition of $(\MM,d_\MM)$ if $\Part^\omega$ is a $\Delta$-bounded partition of $(\MM,d_\MM)$ for every $\omega\in \Omega$. Given $\sigma\ge 0$,  a random $\Delta$-bounded partition $\Part$  of a metric space  $(\MM,d_\MM)$ is said to be $\sigma$-separating if the following requirement holds:
\begin{equation}\label{eqn:sep}
          \forall x,y\in \MM,\qquad   \prob \big[\Part(x) \neq \Part(y)\big] \le  \frac{\sigma}{\Delta}d_\MM(x, y).
        \end{equation}

The separation modulus of $(\MM,d_\MM)$, denoted $\SEP(\MM)$, is the infimum over $\sigma\ge 0$ such that for every $\Delta>0$ there exists a random $\Delta$-bounded   $\sigma$-separating partition $\Part_{\!\!\Delta}$  of $\MM$; if no such random partition exists, then set $\SEP(\MM)=\infty$. This important concept has been introduced by~\cite{bartal1996probabilistic}, see~\cite{naor2024extension} for the history.  

For an infinite metric space $(\MM,d_\MM)$, define its separation growth $\big\{\SEP^n(\MM)\big\}_{n=1}^\infty$ by
$$
\forall n\in \N,\qquad \SEP^n(\MM)\eqdef \sup_{\substack{\sub\subset \MM\\ |\sub|\le n}} \SEP(\sub). 
$$
We say that $(\MM,d_\MM)$ has nontrivial separation growth if 
\begin{equation}\label{eq:nontrivial sep}
\lim_{n\to \infty} \frac{\SEP^n(\MM)}{\log n}=0.
\end{equation} 
The term ``nontrivial'' is used here since it was proved in~\cite{bartal1996probabilistic} that {\em every} metric space $(\MM,d_\MM)$ satisfies $\SEP^n(\MM)=O(\log n)$. 

By~\cite{ccg98,LN03,naor2024extension}, if $1\le p\le 2$, then $\SEP^n( L_p)=o(\log n)$ if $$\lim_{n\to \infty} (p-1)\frac{\log \log n}{\log\log\log n}=\infty,$$ while if $\liminf_{n\to \infty} (p-1)\log\log n<\infty$, then  $\SEP^n( L_p)\gtrsim \log n$. See~\cite[Section~1.7.6]{naor2024extension} for more on this (including  sharp bounds for fixed $1<p\le 2$), where it is conjectured that $\SEP^n( L_p)=o(\log n)$  if and only if $\lim_{n\to \infty} (p-1)\log\log n =\infty$. 

It was asked in~\cite[Question~1]{naor2017probabilistic}  (reiterated in~\cite[Question~83]{naor2024extension}) if $ L_p$ has nontrivial separation growth when $2<p<\infty$  (see~\cite[Section~1.7.6]{naor2024extension} for the relation to metric dimension reduction). This question was answered affirmatively  in~\cite{krauthgamer2025lipschitz} by showing that $\SEP^n( L_p)\lesssim (\log n)^{1-1/p}$ for every $(p,n)\in (2,\infty)\times \N$. It was also asked in~\cite{naor2017probabilistic,naor2024extension} if, in fact,  for every $(p,n)\in (2,\infty)\times \N$ we have $\SEP^n( L_p)\lesssim_p \sqrt{\log n}$, which would be asymptotically optimal as $n\to \infty$ for every fixed $2<p<\infty$ (see below).  Here we prove that this sharp evaluation of the largest possible  separation modulus of an $n$-point subset of $ L_p$ indeed holds:
\begin{theorem}\label{cor:lp} For every $2<p<\infty$ and every  $n\in \{2,3,
\ldots\}$ we have  $\SEP^n( L_p)\asymp_p\sqrt{\log n}$. More precisely,
\begin{equation}\label{eq:our sep bounds}
\forall (p,n)\in (2,\infty)\times \N,\qquad {\textstyle \sqrt{\log n}} \lesssim \SEP^n( L_p) \lesssim p^2 {\textstyle \sqrt{\log n}}.
\end{equation}
\end{theorem}

The new content of~\eqref{eq:our sep bounds} is its upper bound on $\SEP^n( L_p)$; the lower bound on $\SEP^n( L_p)$ in~\eqref{eq:our sep bounds} holds since $\SEP^n( L_2)\asymp \sqrt{\log n}$, by~\cite{ccg98}, while   $\SEP^n(\bX)\ge \SEP^n( L_2)$ for every infinite dimensional Banach space $\bX$,  by Dvoretzky's theorem~\cite{Dvo60}. So, we now know that $\SEP^n( L_p)=o(\log n)$  as $n\to \infty$  if $2<p=o\left(\sqrt[4]{\log n}\right)$, and it fails if $p\gtrsim \log n$ because in that case $\ell_{\!\! \infty}^n$ embeds with distortion $O(1)$ into $ L_p$, any $n$-point metric space $\MM$ embeds into $\ell_{\!\! \infty}^n$, and we already recalled that by~\cite{bartal1996probabilistic} there exists such a space for which $\SEP(\MM)\gtrsim \log n$. It remains open  (likely requiring a substantially new idea) to understand what happens in the remaining range $\sqrt[4]{\log n}\lesssim p=o(\log n)$. 

\begin{remark}{\em  The first arXiv posting~\cite{NR25-v1} of our work (which the present article supersedes) suppressed the dependence on $p$ in Theorem~\ref{cor:lp} because the main matter at hand is determining the asymptotic growth rate of $\SEP^n( L_p)$ as $n\to \infty$; this is what was asked in~\cite{naor2017probabilistic} and what Theorem~\ref{cor:lp}  answers. Nevertheless, understanding the dependence on $p$ is of value in its own right (partially because in applications sometimes $p$ itself is allowed to depend on $n$), and an inspection of the proof in~\cite{NR25-v1} reveals that it yields the estimate $\SEP^n( L_p)\le e^{O(p)}\sqrt{\log n}$.   This was achieved in~\cite{NR25-v1} by combining the approach of~\cite{BG14,krauthgamer2025lipschitz} with a bootstrapping argument, which, as explained in~\cite[Remark~8]{NR25-v1} can also be realized as an iterative procedure. The question of improving the dependence on $p$ was broached in the subsequent work~\cite{KPS25}, which enhanced the recursion in a novel and interesting way to get the bound $\SEP^n( L_p)\lesssim p^4\sqrt{\log n}$. The proof herein of Theorem~\ref{cor:lp}  incorporates a more geometric approach by examining radially bounded random partitions and relying on a new property of the Mazur map that could be of use elsewhere; an overview of the ideas and steps of that proof appears in Section~\ref{intro sketch} below.  }
\end{remark}

We will also prove the following theorem for doubling subsets of $L_p$:

\begin{theorem}\label{thm:doubling sep} If  $p,\ul\ge 2$, then   every $\ul$-doubling subset $\cD$ of $L_p$ satisfies:
$$\SEP(\mathscr{D})\lesssim \Big(p^2{\textstyle \sqrt{\log p}}\Big){\textstyle \sqrt{\log \ul}.}$$
\end{theorem}

\subsubsection{From separation of neighborhoods  to Lipschitz extension} The similarity between the conclusions of Theorem~\ref{cor:lp}  and Theorem~\ref{thm:new ext intro}, as well as Theorem~\ref{thm:doubling sep} and Theorem~\ref{thm:doubling ext}, are not coincidental. The link is provided by Theorem~\ref{thm:ext from sep of neighborhood} below, which we will deduce quickly in Section~\ref{sec:exte neighborhood version proof}  from~\cite{LN05} (with input from~\cite{naor2024extension}).   

In the formulation of Theorem~\ref{thm:ext from sep of neighborhood}, as well as throughout the ensuing discussions, we will use the following (natural but nonstandard) notation. Given a metric space $(\MM,d_\MM)$ and $\sub\subset \MM$, for every $r\ge 0$ we will denote the $r$-neighborhood $\sub$ in $\MM$ by: 
\begin{equation}\label{eq:def of r neighborhood notation}
B_\MM(\sub,r)\eqdef \bigcup_{x\in \sub} B_\MM(x,r).
\end{equation}
Given $\Delta>0$ we will denote by $\SEP_\Delta(\MM)$ the infimum over $\sigma\ge 0$ such that there is a random $\Delta$-bounded   $\sigma$-separating partition of $\MM$ (again, if no such random partition exists, then write $\SEP_\Delta(\MM)=\infty$). Thus, $$\SEP(\MM)=\sup_{\Delta>0}\SEP_\Delta(\MM).$$

\begin{theorem}\label{thm:ext from sep of neighborhood} Suppose that $(\MM,d_\MM)$ is a metric space and that $\sub\neq \emptyset$ is a locally compact subset of $\MM$. Then, the following Lipschitz extension estimate holds for every Banach space $(\bZ,\|\cdot\|_\bZ)$  and every $L>0$: 
\begin{equation}\label{eq:lip extension from sep of neighborhood}
\ee(\MM,\sub;\bZ)\lesssim L+\sup_{\Delta>0} \SEP_\Delta\Big(B_\MM\big(\sub,\frac{1}{L}\Delta\big)\Big).
\end{equation}
\end{theorem}

Our main contribution to randomized clustering of (neighborhoods of) subsets of $L_p$ is:

\begin{theorem}\label{thm:neighborhoods in Lp} There exists a universal constant $\gamma>0$ with the following property. Suppose that $2\le p<\infty$. For every $n\in \{2,3,\ldots\}$, if $\sub$ is an $n$-point subset of $L_p$, then for every $\Delta>0$ we have:
\begin{equation}\label{eq:finite neighborhood is theorem}
 \SEP_\Delta\Big(B_{L_p}\big(\sub,\frac{\gamma}{p}\Delta\big)\Big)\lesssim p^2{\textstyle\sqrt{\log n}}. 
\end{equation}
Furthermore, for every $\ul \ge 2$, if $\cD$ is a subset of $L_p$ that is $\ul$-doubling, then for every $\Delta>0$ we have:
\begin{equation}\label{eq:doubling neighborhood is theorem}
\SEP_\Delta\Big(B_{L_p}\big(\cD,\frac{\gamma}{p}\Delta\big)\Big)\lesssim \Big(p^2{\textstyle\sqrt{\log p}}\Big){\textstyle\sqrt{\log \ul}}. 
\end{equation}
\end{theorem}

Because any  neighborhood of a set contains the set itself, Theorem~\ref{thm:neighborhoods in Lp} strengthens (the upper bound on $\SEP^n(L_p)$ in) Theorem~\ref{cor:lp} and Theorem~\ref{thm:doubling sep}.  Furthermore, Theorem~\ref{thm:new ext intro} and Theorem~\ref{thm:doubling ext} follow from Theorem~\ref{thm:neighborhoods in Lp} by invoking Theorem~\ref{thm:ext from sep of neighborhood}. Indeed, any Lipschitz function from subset of a metric space  that take values in a complete metric space automatically extends to the closure of its domain, so in Theorem~\ref{thm:neighborhoods in Lp} we may assume that $\cD$ is a closed  subset of $L_p$.  As $\cD$ is also assumed in Theorem~\ref{thm:neighborhoods in Lp} to be doubling, it is locally compact (see e.g.~\cite[Lemma~4.1.14]{HK15}, or notice that any ball in $\cD$ is compact because it is totally bounded and complete). The assumptions of Theorem~\ref{thm:ext from sep of neighborhood} therefore hold with $L\asymp p$ by  Theorem~\ref{thm:neighborhoods in Lp}.    Note here  that even when $\sub\subset L_p$ is finite, the above justification of Theorem~\ref{thm:new ext intro} involves random partitions of infinite sets (neighborhoods of $\sub$ in $L_p$), for which  measurability considerations are pertinent.

\subsubsection{From separation and extension to Euclidean embedding}\label{sec:sep ext to embed in intro}   The link between Theorem~\ref{thm:new ext intro} and Theorem~\ref{thm:all >2}  is furnished by Theorem~\ref{thm:sep plux ext gives c} below, which is an embedding statement of independent interest; one can view it  as a variant of the ``measured descent'' embedding method~\cite{KLMN05} that incorporates separating partitions rather than the padded  partitions that occur in~\cite{KLMN05}. The proof of Theorem~\ref{thm:sep plux ext gives c}, which appears in Section~\ref{sec:proof of multiscale}, combines multiple ingredients, many of which (but not all) refine the reasoning in~\cite{ALN08}.

For Theorem~\ref{thm:sep plux ext gives c}, recall that $\ee^k(\MM)=\sup\{\ee(\sub):\ \sub\subset \MM\ \wedge\  |\sub|\le k\}$  was defined in Remark~\ref{rem:superscript version of finatary extension}.
\begin{theorem}\label{thm:sep plux ext gives c}  For every $n\in \{3,4,\ldots\}$, every $n$-point metric space $(\MM,d_\MM)$  satisfies
\begin{equation}\label{eq:decednt for separataion}
\cc_2(\MM)\lesssim {\textstyle\sqrt{(\log n)\log\log n}}\bigg(\sum_{k=2}^n \frac{\SEP^k(\MM)^2\ee^k(\MM; L_2)^4}{k(\log k)^2}\bigg)^{\frac12}.
\end{equation}
\end{theorem}
 For every $k\in \{2,3,\ldots\}$ we have   $\SEP^k(L_p)\lesssim p^2\sqrt{\log k}$ by    Theorem~\ref{cor:lp}, and  by~\cite{NPSS06} we have $\ee^k(L_p)\le \ee(L_p)\lesssim \sqrt{p}$ . A substitution of these two estimates into Theorem~\ref{thm:sep plux ext gives c} implies Theorem~\ref{thm:all >2} as follows: 
\begin{align*}
\cc_2^n(L_p)\lesssim {\textstyle\sqrt{(\log n)\log\log n}}&\bigg(\sum_{k=2}^n \frac{p^4(\log k)\cdot p^2}{k(\log k)^2}\bigg)^{\frac12}\\&\asymp p^3{\textstyle\sqrt{(\log n)\log\log n}}\bigg(\int_1^{n}\frac{\ud s}{s\log s}\bigg)^{\frac12}= p^3{\textstyle\sqrt{\log n}}\log\log n. 
\end{align*}

\subsubsection{Localization and induction on scales, and a radial property of the Mazur map}\label{intro sketch} The purpose of this section is to provide an overview of the key ingredients of our proof of Theorem~\ref{thm:neighborhoods in Lp}.  

Suppose that $(\MM,d_\MM)$ is a metric space and $\sub\subset \MM$. The $d_\MM$-(circum)radius of $\sub$ is defined as follows:
\begin{equation}\label{eq:radius notation}
\rad_\MM(\sub)\eqdef \inf\big\{r\in [0,\infty]:\ \exists x\in \MM\ \mathrm{such\ that\ }   B_\MM(x,r)\supseteq \sub\big\},
\end{equation}
Given $\Delta>0$, we will say that a random partition $\Part$ of $\sub$ is radially $\Delta$-bounded with respect to $\MM$ if: 
\begin{equation}
\forall x\in \sub,\qquad \rad_\MM\big(\Part(x)\big)\le \Delta.
\end{equation}
A random partition $(\omega\in \Omega)\mapsto \Part^\omega$ of $\sub$, defined on some probability space $(\Omega,\prob)$, is radially $\Delta$-bounded with respect to $\MM$ if $\rad_\MM(\Part^\omega(x))\le \Delta$ for every $\omega\in \Omega$ and every $x\in \sub$. Given a random partition $\Part$ of $\sub$ that is radially $\Delta$-bounded with respect to $\MM$, say that it is $\sigma$-separating for some $\sigma\ge 0$ if~\eqref{eqn:sep} holds.

Denote by $\TSEP_\Delta(\sub;\MM)$ the infimum over $\sigma\ge 0$ such that there exists a $\sigma$-separating  random partition of $\sub$ that is radially $\Delta$-bounded with respect to $\MM$ (as always, when no such partition exists we define this parameter to be $\infty$). 
Recalling that $\diam_\MM(\cdot)$ denotes the $d_\MM$-diameter,  the following  bounds hold: 
\begin{equation}\label{eq:rad2diam}
\forall \emptyset\neq \sub\subset \MM,\qquad \rad_\MM(\sub)\le\diam_\MM(\sub)\le 2\rad_\MM(\sub).
\end{equation} 
Consequently, the following simple general relations between the separation and radial separation moduli are satisfied for every metric space $(\MM,d_\MM)$, every $\emptyset\neq \sub\subset \MM$, and every $\Delta>0$:  
\begin{equation}\label{eq:TSEP relations}
\TSEP_\Delta (\sub; \MM) \le \SEP_\Delta (\sub) \le 2\TSEP_{\frac{\Delta}{2}} (\sub; \MM),
\end{equation}

Even though the above   radial variants  of the notions of random $\Delta$-bounded and $\sigma$-separating partitions may seem to be  nuanced minor tweaks of the standard (by now classical) definitions, we will see that they influence our results in  a way that is more dramatic   than what one might initially expect. 

The following lemma is a localization and induction on scales principle for random radial separation: 

\begin{lemma}\label{lem:induction and localization} Suppose that $(\MM,d_\MM)$ is a separable metric space and that $\emptyset\neq \sub\subset \MM$ satisfies 
 \begin{equation}\label{decay assumption on sep}
    \lim_{\Delta \to \infty} \frac{1}{\Delta}\SEP_{\Delta} (\sub) =0.
    \end{equation}
    Then, the following estimate holds for every $\Delta>0$ and $K>1$:  
      \begin{equation}\label{eqn:sum geometric series}
        \TSEP_{\Delta} (\sub; \MM) \le \sum_{s=0}^\infty \frac{1}{K^s}\lim_{\e\to 0^+} \sup_{z \in \MM} \TSEP_{ K^s \Delta} \big(\sub \cap B_\MM (z, K^{s+1}\Delta+\e); \MM\big).
    \end{equation} 
\end{lemma}
The limits in~\eqref{eqn:sum geometric series} exist as the summands are nondecreasing with $\e$. We used the term ``localization'' to describe Lemma~\ref{lem:induction and localization} as the $s$-summand in~\eqref{eqn:sum geometric series} depends only on the intersection  $\sub \cap B_\MM (z, K^{s+1}\Delta+\e)$, which is a  ``local snapshot''  of $\sub$. We used the term ``induction on scales'' to describe Lemma~\ref{lem:induction and localization} as the left hand side of~\eqref{eqn:sum geometric series} treats partitions at scale $\Delta$ while the right hand side of~\eqref{eqn:sum geometric series} considers partitions at an increasing sequence of  larger scales. Our proof of Lemma~\ref{lem:induction and localization} is  an iterative use of the following observation. One can partition $\sub$ into clusters of radius at most $\Delta$ by first partitioning it into clusters of radius at most $K\Delta$; each of those clusters is thus contained in a ball of radius $K\Delta+\e$, so we can proceed to  refine the aforementioned initial  partition by partitioning each of the enclosing balls  into pieces of radius  at most $\Delta$. 

A standard tool for analysing the geometry  $L_p$ is the classical Mazur map~\cite{mazur1929remarque} to $L_2$ (and shifts thereof), whose restriction to balls has well-understood (and widely-used)  quite favorable uniform continuity properties, i.e.,  it is well-behaved on the local snapshots of $\sub$ that appear in the right hand side of~\eqref{eqn:sum geometric series}. 

For $1\le p,q<\infty$, the Mazur map $M_{p\to q}:L_p\to L_q$ is defined~\cite{mazur1929remarque} by: 
\begin{equation}\label{eq:def mazur map}
\forall \phi\in L_p,\ \forall t\in [0,1], \qquad M_{p\to q}(\phi)(t)\eqdef |\phi(t)|^{\frac{p}{q}}\sgn\big(\phi(t)\big).
\end{equation}
If $p>2$, then $M_{p\to 2}$ is Lipschitz by~\cite{mazur1929remarque}, so we  may consider the following normalization of the Mazur map that makes it be a $1$-Lipschitz function from $L_p$ into the Hilbert space $L_2$:  
\begin{equation}\label{eq:normalized mazur}
\widetilde{M}_{p\to 2}\eqdef \frac{1}{\|M_{p\to 2}\|_{\Lip( L_p; L_2)}}M_{p\to 2}. 
\end{equation}
The normalization factor in~\eqref{eq:normalized mazur} satisfies $\|M_{p\to 2}\|_{\Lip( L_p; L_2)}< p$, as computed in~\cite[equation~(5.32)]{naor2014comparison}.

We will prove herein that there exists a universal constant $\beta>0$ such that the following inclusion holds: 
\begin{equation}\label{eq:radial mazur normalized}
\forall \phi\in B_{L_p},\qquad \big(\widetilde{M}_{p\to 2}\big)^{-1} \bigg(B_{L_2} \Big(\widetilde{M}_{p\to 2}(\phi),\frac{\beta}{p}\Big)\bigg)\subset B_{L_p}\Big(\frac{1}{p}\phi,1-\frac{1}{4p}\Big),
\end{equation}
where $B_{L_p}$ is the ball in $L_p$ of radius $1$ centered at the $0$-function; more generally, the unit ball centered at the origin of a Banach space $(\bX,\|\cdot\|_\bX)$ will always be denoted below by $B_\bX=\{x\in \bX:\ \|x\|_\bX\le 1\}=B_\bX(0,1)$.  

The upshot of~\eqref{eq:radial mazur normalized} is that the restriction of $\widetilde{M}_{p\to 2}$ to  $B_{L_p}$ gives a $1$-Lipschitz function $f$ into a Hilbert space with the property that for any point $x\in B_{L_p}$, the pullback under $f$ of the Hilbertian ball centered at $f(x)$ of radius $\beta/p$ is contained in an $L_p$-ball whose radius is a definite number (independent of $x$) smaller than $1$, albeit ever so slightly as $p\to \infty$. 

Returning to the  setting of Theorem~\ref{thm:neighborhoods in Lp}, the above discussion sets the stage for the following procedure. First  apply Lemma~\ref{lem:induction and localization} with $K=1/(1-1/(4p))$; the choice of this value of $K$ will allow us to take advantage of the radius decrease in~\eqref{eq:radial mazur normalized} from $1$ (the radius of the domain of $f$) to $1-1/(4p)$. 

Thus, suppressing (only for the purpose of the present proof sketch) the small additive $\e>0$ correction of the radius that appears in~\eqref{eqn:sum geometric series} (the complete reasoning in Section~\ref{sec:radial preserves sep}  takes this correction into account), we may now fix $z\in \MM$, $\Delta>0$ and an integer $s\ge 0$ and focus our attention on obtaining a separating partition  that is $K^s\Delta$-radially bounded with respect to $L_p$  of the local snapshots $\sub \cap B_{L_p} (z, K^{s+1}\Delta)$ and $\cD \cap B_{L_p} (z, K^{s+1}\Delta)$. By translating by $z$ and rescaling by $K^{s+1}\Delta$, we are thus interested in   obtaining a separating partition  that is $(1/K)$-radially bounded with respect to $L_p$  of  $\sub \cap B_{L_p}$ and $\cD \cap B_{L_p}$.

 As these local snapshots are contained in the unit ball of $L_p$, we may apply the Mazur map to them, for which~\eqref{eq:radial mazur normalized} holds. The resulting images are now subsets of $L_2$, on which we may use Euclidean geometric considerations  to randomly partition them with good separation  into clusters of $L_2$-diameter   at most $\beta/p$ (using the partitioning results of~\cite{ccg98,LN05,naor2024extension} as well as the Kirszbraun extension theorem~\cite{Kir34} and the Johnson--Lindenstrauss dimension reduction lemma~\cite{JL82}). The pullback under the normalized Mazur map of the resulting partition will then have good separation, and by~\eqref{eq:radial mazur normalized} it will have $L_p$-radius a most $1-1/(4p)=1/K$, as desired. Note that the center $(1/p)\phi$ of the ball in the right hand side of~\eqref{eq:radial mazur normalized} need not belong to the corresponding local snapshot of $\sub$ and $\cD$, but this is permitted by the definition of radially bounded random partitions with respect to the super-space $L_p$.  By substituting the resulting bounds on the radial separation moduli of the snapshots and summing over $s$,  Theorem~\ref{thm:neighborhoods in Lp} follows.  

\begin{remark} 
 {\em As we already mentioned, good bounds on the moduli of uniform continuity of the Mazur map are well known~\cite[equation~(5.32)]{naor2014comparison} and widely used. If one incorporates into the above reasoning those bounds in place of the new radial property~\eqref{eq:radial mazur normalized} of the Mazur map that we obtain herein, one arrives at exponential dependence on $p$ in the right had side of~\eqref{eq:finite neighborhood is theorem} and of~\eqref{eq:doubling neighborhood is theorem}. In Section~\ref{sec:vanilla localized} we will prove that such an exponential loss is inherent to the aforementioned alternative route, and it will be incurred by any  uniform homeomorphism from the unit ball of $L_p(\mu)$ into Hilbert space,  not only by the Mazur map. What~\eqref{eq:radial mazur normalized} achieves is power-type dependence on $p$ if one is willing to settle for a small gain in the radius of the pre-image. This forces us to work with $K$ close to $1$ in Lemma~\ref{lem:induction and localization}, so the geometric convergence in~\eqref{eqn:sum geometric series} is slow, but as $K-1$ is of order $1/p$, this leads to a multiplicative  loss which is also just $O(p)$. 
 }
\end{remark}

\subsection{Beyond \texorpdfstring{$L_p$}{Lp}} Here we will discuss generalizations of our results to Banach spaces that need not be $ L_p$ (or any  $L_p(\mu)$ space, for which all the statements are identical). Those who are interested only the aforementioned statements for $L_p$ can harmlessly skip this material, which assumes some (entirely standard, per e.g.~\cite{BL00}) background from Banach space theory that is not pertinent to the ensuing proofs.

An inspection of the reasoning herein reveals that if  $(\bX,\|\cdot\|_\bX)$ is an infinite dimensional Banach space, then the conclusions of Theorem~\ref{thm:all >2}, Theorem~\ref{thm:new ext intro},  Theorem~\ref{thm:doubling ext}, Theorem~\ref{thm:doubling sep}, Theorem~\ref{cor:lp}, and Theorem~\ref{thm:neighborhoods in Lp}  hold with $L_p$ replaced by $\bX$ provided that there is an injective  Lipschitz function from the unit ball of $\bX$ into a Hilbert space whose inverse is uniformly continuous, with the only difference being that in this (much) more general setting the dependence on $p$ should be replaced by a dependence on the Lipschitz constant and the modulus of uniform continuity of the inverse  in the aforementioned  assumption on $\bX$. There is  substantial literature~\cite{mazur1929remarque,CL93,OS94,Cha95,Dah95,BL00,Ray02,AP10,Ric15,Cor24,CFGT25} that obtains uniform homeomorphisms from the unit ball of certain Banach spaces into a Hilbert space, which we will next describe.

 The important work~\cite{OS94} implies that if $(\bX,\|\cdot\|_\bX)$ is a Banach space with an unconditional basis whose modulus of uniform smoothness has power-type $2$ and whose modulus of uniform convexity has power type $p$, then  there is an injective  Lipschitz function from the unit ball of $\bX$ into a Hilbert space whose inverse is $(2/p)$-H\"older,\footnote{\cite{OS94} obtains such a function into $ L_1$ rather than into $ L_2$, but it is straightforward to adapt the reasoning so as to obtain a function  into $ L_2$. Furthermore, these Lipschitz and H\"older assertions follow from an inspection of the proofs in~\cite{OS94},  but they are not stated there explicitly. It is perhaps simplest to verify them by examining the exposition of~\cite{OS94} in the monograph~\cite{BL00}.} whence the theorems that we listed above hold for such $\bX$. The work~\cite{Cha95} generalizes the aforementioned result of~\cite{OS94} to Banach lattices, though in a manner that is lossy with respect to the  H\"older estimates that are obtained; we expect, however, that with not much more care the reasoning in~\cite{Cha95} could be adapted to show that if  $(\bX,\|\cdot\|_\bX)$ is a Banach lattice whose modulus of uniform smoothness has power-type $2$ and whose modulus of uniform convexity has power type $p$, there is an injective  Lipschitz function from the unit ball of $\bX$ into a Hilbert space whose inverse is $(2/p)$-H\"older. 
 
 By~\cite{Ric15}, if $(\bX,\|\cdot\|_\bX)$ is the Schatten--von Neumann trace class $\mathsf{S}_p$, or more generally if it  is a noncommutative $L_p$ space over any von Neumann algebra, then  there is   an injective  Lipschitz function from the unit ball of $\bX$ into a Hilbert space whose inverse is $(2/p)$-H\"older, thus showing that the theorems that we listed above  can be generalized to this noncommutative setting. We conjecture that the same holds for any unitary ideal $\mathsf{S}_{\mathbf{E}}$ over any Banach space $(\mathbf{E},\|\cdot\|_{\mathbf{E}})$ with a $1$-symmetric basis, provided that $\mathbf{E}$ has modulus of uniform smoothness of power-type $2$ and modulus of uniform convexity has power type $p$; we leave this as an open question that seems quite accessible (perhaps even straightforward) with current technology (the work~\cite{HSZ22} could be especially helpful here).  
 
There are valuable works that construct uniform homeomorphisms into Hilbert space of unit balls in interpolation spaces (see~\cite{Dah95,Cor24,CFGT25} and~\cite[Section~3 of Chapter~9]{BL00}), but they do not yield the above Lipschitz estimate and it would be interesting to understand (perhaps requiring a substantial new idea) whether such a Lipschitz estimate  could be derived in this context as well (under suitable assumptions).

\begin{problem} {\em The extent to which the inclusion~\eqref{eq:radial mazur normalized} generalizes beyond $L_p$ is uncharted (and currently unstudied) territory. In particular, it would be interesting to determine if~\eqref{eq:radial mazur normalized} holds for the noncommutative Mazur map on the  Schatten--von Neumann trace class $\mathsf{S}_p$ when $p>2$. And, it would be worthwhile to understand how~\eqref{eq:radial mazur normalized} should be changed when $L_p$ is replaced by a Banach space $(\bX,\|\cdot\|_\bX)$ with an unconditional basis whose modulus of uniform smoothness has power-type $2$, the pertinent question being how the radii of the balls that occur in~\eqref{eq:radial mazur normalized} should depend on geometric characteristics of $\bX$.  }
\end{problem}

\subsection{Roadmap} The rest of the ensuing text is organized  as follows. We will start by proving Theorem~\ref{thm:sep plux ext gives c}  in Section~\ref{sec:proof of multiscale}. As we explained in Section~\ref{sec:sep ext to embed in intro}, this will complete the reduction of Theorem~\ref{thm:all >2} to Theorem~\ref{cor:lp}. Thus,  what will remain to be done after Section~\ref{sec:proof of multiscale}  is to prove Theorem~\ref{thm:neighborhoods in Lp}, which we have already seen implies the rest  of the new results that are obtained herein; the corresponding extension statements rely on Theorem~\ref{thm:ext from sep of neighborhood}, which is not stated in the literature but readily follows from the link between randomized separation and Lipschitz extension that was discovered in~\cite{LN05}. The justification of Theorem~\ref{thm:ext from sep of neighborhood} using~\cite{LN05} appears in Section~\ref{sec:exte neighborhood version proof}. In Section~\ref{sec:vanilla localized} we will set the stage for proving the inclusion~\eqref{eq:radial mazur normalized}, which we will do in Section~\ref{sec:all p>2}. While the concepts that are discussed in  Section~\ref{sec:vanilla localized} are not needed for the proof Theorem~\ref{thm:neighborhoods in Lp}, those motivate their radial counterparts that are used for this purpose and are discussed in Section~\ref{sec:all p>2}. Also, Section~\ref{sec:vanilla localized} proves an impossibility result that explains the need for the variants that Section~\ref{sec:all p>2} provides. Lemma~\ref{lem:induction and localization} is proved in Section~\ref{sec:induction on scales}, and the proof of Theorem~\ref{thm:neighborhoods in Lp} is completed in Section~\ref{sec:radial preserves sep}.

\section{Measured descent for separated partitions (in the presence of Lipschitz extension) }\label{sec:proof of multiscale}

The purpose of this section is to prove Theorem~\ref{thm:sep plux ext gives c}. As we explained in Section~\ref{sec:sep ext to embed in intro}, this will reduce Theorem~\ref{thm:all >2} to  Theorem~\ref{thm:neighborhoods in Lp}, which will be proven later,  in Section~\ref{sec:radial preserves sep}  below.

\subsection{Weakly bi-Lipschitz embeddings}\label{sec:single scale}
Section~\ref{sec:bilip intro} recalled (for brevity of the Introduction) only the notion of Euclidean distortion, but it is  standard (and fruitful) to study the analogous concept for embeddings into an arbitrary metric space.  A metric space $(\MM,d_\MM)$ has a  bi-Lipschitz embedding of distortion $D\ge 1$ into a metric space $(\cZ,d_\cZ)$ if there is a nonconstant Lipschitz function $f:\MM\to \cZ$ satisfying
\begin{equation}\label{eq:def bilip general}
\forall x,y\in \MM,\qquad d_\cZ\big(f(x),f(y)\big)\ge \frac{\|f\|_{\Lip(\MM;\cZ)}}{D}d_\MM(x,y),
\end{equation}
The $\cZ$-distortion $\cc_\cZ(\MM)$ of $\MM$ is  the infimum over those $D\in [1,\infty]$ for which an embedding $f$ as above exists. We also denote  $\cc_\cZ^n(\MM)=\sup\{\cc_\cZ(\sub):\ |\sub|\subset \MM\ \wedge\ 2\le |\sub|\le n\}$ for every integer $n\ge 2$. 

%where in~\eqref{eq:def bilip general}, as well as throughout what follows, $\|f\|_{\Lip(\MM;\cZ)}$ denotes the Lipschitz constant of $f$, namely, it is the %supremum of $d_\cZ(f(x),f(y))/d_\MM(x,y)$ over all distinct $x,y\in \MM$. 

There is a  ``one scale at a time'' variant of the above classical setup, which is an important and commonly used concept in the study of metric embeddings; the relevant terminology, which we next recall, was introduced in~\cite[p.~192]{NPSS06} (the reason for the nomenclature is further articulated in~\cite[Section 7.2]{naor2014comparison}).
\begin{definition}\label{def:npss} A metric space $(\MM,d_\MM)$ admits a weakly bi-Lipschitz embedding of distortion $D\ge 1$ into a metric space $(\cZ,d_\cZ)$ if for any $\Delta>0$ there is a nonconstant Lipschitz function $f=f_\Delta:\MM\to \cZ$ satisfying
\begin{equation*}
\forall x,y\in \MM,\qquad d_\MM(x,y)\ge \Delta\implies d_\cZ\big(f(x),f(y)\big)\ge \frac{\|f\|_{\Lip(\MM;\cZ)}}{D}\Delta. 
\end{equation*}
\end{definition}

Analogously  to the above notation, in the setting of Definition~\ref{def:npss}  let $\dd_\cZ(\MM)$ be the infimum over those $D\in [1,\infty]$ for which $(\MM,d_\MM)$ has a weakly bi-Lipschitz embedding of distortion $D$ into $(\cZ,d_\cZ)$, and write 
$$
\forall n\in \{2,3,\ldots\},\qquad \dd_\cZ^n(\MM)\eqdef \sup_{\substack{\sub \subset \MM\\ 2\le |\sub|\le n}}\dd_\cZ(\sub).
$$
We will also use the shorter notation  $\dd_{ L_p}(\MM)=\dd_p(\MM)$ and $\dd_{ L_p}^n(\MM)=\dd_p^n(\MM)$ for  $p\ge 1$ and $n\in \{2,3,\ldots\}$. 

Suitably defined ``spatially localized'' variants of Definition~\ref{def:npss} will be discussed in Section~\ref{sec:vanilla localized} below, as they are important for proving the results that obtained herein. In this section, though, it suffices to consider only the above ``vanilla'' notion of ``one scale at a time'' embedding, as well as the following ad hoc notation for an obvious (also standard)  ``two-sided scale-localized'' version thereof. 

Given metric spaces $(\MM,d_\MM)$ and $(\cZ,d_\cZ)$, let $\widehat{\dd}_\cZ(\MM)$ denote the infimum over those $D\in [1,\infty]$ with the property that  for any $\Delta>0$ there exists a nonconstant Lipschitz function $f=f_\Delta:\MM\to \cZ$ satisfying
\begin{equation}\label{eq:2Delta}
\forall x,y\in \MM,\qquad \Delta\le d_\MM(x,y)\le 2\Delta\implies d_\cZ\big(f(x),f(y)\big)\ge \frac{\|f\|_{\Lip(\MM;\cZ)}}{D}\Delta. 
\end{equation}
We also use the shorter notation  $\widehat{\dd}_{ L_p}(\MM)=\widehat{\dd}_p(\MM)$ and $\widehat{\dd}_{ L_p}^n(\MM)=\widehat{\dd}_p^n(\MM)$ for  $p\ge 1$ and $n\in \{2,3,\ldots\}$, where  
$$
\forall n\in \{2,3,\ldots\},\qquad \widehat{\dd}_\cZ^n(\MM)\eqdef \sup_{\substack{\sub \subset \MM\\ 2\le |\sub|\le n}}\widehat{\dd}_\cZ(\sub).
$$

\begin{remark}\label{rem:2Delta} {\em The factor $2$  in~\eqref{eq:2Delta} was fixed for concreteness, but it is an arbitrary choice of minor significance to the context of what we study herein. Specifically, suppose that $1<\alpha<\beta<\infty$ and for any $\Delta>0$ there exists a $1$-Lipschitz function $f_\Delta:\MM\to  L_2$ satisfying
\begin{equation}\label{eq:alpha delta version}
\forall x,y\in \MM,\qquad \Delta\le d_\MM(x,y)\le \alpha\Delta\implies \|f_\Delta(x)-f_\Delta(y)\|_{\! L_2}\ge \frac{\Delta}{D}. 
\end{equation}
Set $k= \lceil (\log\beta)/\log \alpha\rceil-1$. If $x,y\in \MM$ satisfy $\Delta\le d_\MM(x,y)\le \beta\Delta$, then  $$\exists i\in \{0,\ldots,k\},\qquad \alpha^i\Delta\le d_\MM(x,y)\le \alpha^{i+1}\Delta,$$ whence $\|f_{\alpha^i\Delta}(x)-f_{\alpha^i\Delta}(y)\|_{\! L_2}\ge \alpha^i\Delta/D\ge \Delta/D$. So, if we define $g_\Delta:\MM\to  L_2\otimes  L_2\cong  L_2$ by
$$
g_\Delta(x)\eqdef \frac{1}{\sqrt{k+1}}\sum_{i=0}^k f_{\alpha^i\Delta}(x)\otimes v_i,
$$
where $v_0,\ldots, v_k$ is an arbitrary orthonormal system in $ L_2$, then $g_\Delta$ is $1$-Lipschitz and
\begin{equation*}
\forall x,y\in \MM,\qquad \Delta\le d_\MM(x,y)\le \beta\Delta\implies \|g_\Delta(x)-g_\Delta(y)\|_{\! L_2}\ge \frac{\Delta}{D\sqrt{k+1}}. 
\end{equation*}
Therefore, the existence for every $\Delta>0$ of a $1$-Lipschitz function satisfying~\eqref{eq:alpha delta version} implies the existence for every $\Delta>0$ of such a function having the analogous property with $\alpha$ replaced by $\beta$ and $D$ replaced by $$D\sqrt{\left\lceil \frac{\log\beta}{\log \alpha}\right\rceil}\asymp_{\alpha,\beta} D.$$}
\end{remark}

\subsubsection{From separation and extension to a weakly bi-Lipschitz embedding}

Lemma~\ref{lem:sep ext product}  below bounds  $\dd_2(\cdot)$ by a product of separation moduli and (squares of) Lipschitz extension moduli, as those that appear in the statement of  Theorem~\ref{thm:sep plux ext gives c}. The  statement of Lemma~\ref{lem:sep ext product} uses the following common  notation and terminology, as well as the nonstandard notation that we will introduce  in~\eqref{eq:def DM}. 

Given a metric space $(\MM,d_\MM)$ and $r>0$, a subset  $\NN$ of $\MM$ is called an $r$-net of $\MM$ if  $\NN$ is $r$-separated, i.e., $d_\MM(a,b)\ge r$ for every distinct $x,y\in \NN$, and $\NN$ is also $r$-dense in $\MM$, i.e., for every $x\in \MM$ there is $a\in \NN$ such that $d_\MM(x,a)\le r$. When we say that a subset $\NN$ of $\MM$ is a net of $\MM$ we mean that there is some $r>0$ such that $\NN$ is an $r$-net of $\MM$.  For each target metric space $(\cT,d_\cT)$,  the Lipschitz extension modulus from nets  $\ee_{\net}(\MM;\cT)$ is  defined to be  the supremum of $\ee(\MM,\NN;\cZ)$ over all possible nets $\NN$ of $\MM$, where we recall that the (standard) notation for Lipschitz extension moduli was already introduced in Section~\ref{sec:lip ext}.   Clearly $\ee_{\net}(\MM;\cT)\le \ee(\MM;\cT)$, but it is worthwhile to single out extension from nets because it is an especially important and useful instance of the Lipschitz extension problem, and there are situations in which bounds on $\ee_{\net}(\MM;\cT)$ are known that are stronger than the available upper bounds on $\ee(\MM;\cT)$; examples of works that either treat or rely on  extension from nets include~\cite{Bou87,Beg99,MN13-ext,Nao15-nonextend,NR17,Nao21-almost,ANN25,BN25,NS25}. 

For Theorem~\ref{thm:sep plux ext gives c}, it is convenient to introduce the following notation for every metric space $(\MM,d_\MM)$: 
\begin{equation}\label{eq:def DM}
\Pi(\MM)=\Pi(\MM,d_\MM)\eqdef \max \Big\{\SEP(\NN)\ee(\MM,\NN; L_2)^2:\  
 \NN\ \mathrm{is\ a\ net\ of\ }\MM\Big\}.
\end{equation}
Thus,
\[
\Pi(\MM)\le \SEP(\MM)\ee_{\net}(\MM;L_2)^2\le \SEP(\MM)\ee(\MM;L_2)^2.
\]
The following lemma  relates the ``separation-extension product'' $\Pi(\MM)$ to the weakly bi-Lipschitz Euclidean distortion $\dd_2(\MM)$:

\begin{lemma}\label{lem:sep ext product} Every finite metric space $(\MM,d_\MM)$ satisfies $\dd_2(\MM)\le 8\Pi(\MM)+1$.
\end{lemma}

\begin{proof} Our goal is to prove that for any $\Delta>0$ there is a $1$-Lipschitz function $f=f_{\Delta}:\MM\to  L_2$ such that
\begin{equation}\label{eq:desired Phi}
\forall x,y\in \MM,\qquad d_\MM(x,y)\ge \Delta\implies \|f(x)-f(y)\|_{\!  L_2}\ge \frac{\Delta}{8\Pi(\MM)+1}.
\end{equation}

Set $m=|\MM|$ and fix $\Delta>0$.  Let  $0<\e<1/2$ and $0<\d<1-2\e$ be auxiliary parameters whose  values will be specified later  to optimize the ensuing reasoning. Fix from now an arbitrary $(\e\Delta)$-net $\NN$ of $\MM$. 

Let $\{v_\mathscr{S}\}_{\cS\subset\NN}\subset L_2$ be an orthonormal system of $2^{|\NN|}$ vectors in $ L_2$, indexed by the subsets of $\NN$. Denoting by $\Omega$  the collection of all the $(\d\Delta)$-bounded partitions of $\NN$, let $\mu$ be a probability measure on $\Omega$ such that: 
\begin{equation}\label{eq:our mu separates}
\forall a,b\in \NN,\qquad \mu \big[\Part\in \Omega:\ \Part(a)\neq \Part(b)\big] \le \frac{\SEP(\NN)}{\d\Delta}d_\MM(a,b).
\end{equation}
We can now define a function $\uppsi:\NN\to L_2(\mu; L_2)$ as follows: 
\begin{equation}\label{eq:here is our phi normalized}
\forall a\in \NN, \forall \Part\in \Omega,\qquad \uppsi(a)(\Part)\eqdef \frac{\Delta\sqrt{\e\d}}{\sqrt{2\Pi(\MM)}}v_{\Part(a)}.
\end{equation}

As $\|v_{\mathscr{S}}-v_{\mathscr{T}}\|_{\! L_2}^2=2$ whenever  $\mathscr{S},\mathscr{T}\subset \NN$ are distinct, every $a,b\in \NN$ satisfy 
\begin{align}\label{eq:lip on net}
\begin{split}
\|\uppsi(a)-\uppsi(b)\|_{\! L_2(\mu; L_2)}&\stackrel{\eqref{eq:here is our phi normalized}}{=}\frac{\Delta\sqrt{\e\d}}{\sqrt{2\Pi(\MM)}} \sqrt{2\mu \big[\Part\in \Omega:\ \Part(a)\neq \Part(b)\big] }\\&\stackrel{\eqref{eq:our mu separates}}{\le} \frac{\sqrt{\e\Delta\SEP(\NN) d_\MM(a,b) }}{\sqrt{\Pi(\MM)}} \stackrel{\eqref{eq:def DM}}{\le}
\frac{\sqrt{\e\Delta d_\MM(a,b)}}{\ee(\MM,\NN;\ell_{\!2})} \le \frac{d_\MM(a,b)}{\ee(\MM,\NN;\ell_{\!2})},
\end{split}
\end{align}
where the last step of \eqref{eq:lip on net} holds as $d_\MM(a,b)\ge \e\Delta$ for distinct $a,b\in \NN$. By \eqref{eq:lip on net}, the Lipschitz constant of $\uppsi$ is at most $1/\ee(\MM,\NN;\ell_{\!2})$, so there is a $1$-Lipschitz function $$\Psi:\MM\to L_2(\mu; L_2)$$ that extends $\uppsi$.

Consider $x,y\in \MM$ with $d_\MM(x,y)\ge \Delta$. Take $a,b\in \NN$ such that $d_\MM(x,a),d_\MM(y,b)\le \e\Delta$. Then, 
\begin{equation}\label{eq:in different atoms}
d_\MM(a,b)\ge d_\MM(x,y)-d_\MM(x,a)-d_\MM(y,b)\ge (1-2\e)\Delta>\d\Delta,
\end{equation}
where the last step of \eqref{eq:in different atoms} is where the assumption $\d<1-2\e$ is used.  Since $\diam_\MM(\Part(a))\le \d\Delta$ for every $\Part\in \Omega$, necessarily $b\notin \Part(a)$. In other words, $\Part(a)\neq \Part(b)$ for every $\Part\in \Omega$. Equivalently, every $\Part\in \Omega$ satisfies $\|v_{\Part(a)}-v_{\Part(b)}\|_{\! L_2}=\sqrt{2}$. Recalling the definition \eqref{eq:here is our phi normalized} of $\uppsi$, because $a,b\in \NN$ and  $\Psi$ extends $\uppsi$ we consequently have $$\|\Psi(a)-\Psi(b)\|_{\! L_2(\mu; L_2)}=\|\uppsi(a)-\uppsi(b)\|_{\! L_2(\mu; L_2)}=\Delta\sqrt{\frac{\e\d}{\Pi(\MM)}}.$$ Therefore,
\begin{align}\label{eq:pass tpo xy}
\begin{split}
\|\Psi(x)-\Psi(y)\|_{\! L_2(\mu; L_2)}\ge {}&
\|\Psi(a)-\Psi(b)\|_{\! L_2(\mu; L_2)}\\
&\quad -\|\Psi(x)-\Psi(a)\|_{\! L_2(\mu; L_2)}
-\|\Psi(y)-\Psi(b)\|_{\! L_2(\mu; L_2)}\\
\ge {}& \frac{\Delta\sqrt{\e\d}}{\sqrt{\Pi(\MM)}}-d_\MM(x,a)-d_\MM(y,b)\\
\ge {}& \bigg(\frac{\sqrt{\e\d}}{\sqrt{\Pi(\MM)}}-2\e\bigg)\Delta,
\end{split}
\end{align}
where the penultimate step of~\eqref{eq:pass tpo xy} uses the fact that $\Psi$ is $1$-Lipschitz. 

The right hand side of~\eqref{eq:pass tpo xy} is maximal for $\e=\e_{\mathrm{opt}}\eqdef \d/(16\Pi(\MM))$, which is a valid choice for the above reasoning provided that the requirement $\d<1-2\e_{\mathrm{opt}}$ is satisfied, i.e., if $\d<8\Pi(\MM)/(8\Pi(\MM)+1)$.  Assuming from now on that this restriction on $\d$ holds, substitute $\e_{\mathrm{opt}}$ into~\eqref{eq:pass tpo xy}. Since  $\Psi$ takes values in an $m$-dimensional subspace of a Hilbert space, by composing $\Psi$ with an isometry between that subspace and $\ell_{\!\!2 }^{m}$ and translating the resulting function so that $0\in  \ell_{\!\! 2}^{m}$ is in its image, we conclude that for  every $0<\d<8\Pi(\MM)/(8\Pi(\MM)+1)$ there is a  $1$-Lipschitz function $\Phi_\d:\MM\to B_{\ell_2^{m}}\big(0,\diam(\MM)\big)$ that satisfies:
\begin{equation*}
\forall x,y\in \MM,\qquad d_\MM(x,y)\ge \Delta\implies \|\Phi_\d(x)-\Phi_\d(y)\|_{\! \ell_{\!\!2}^m}\ge \frac{\d\Delta}{8D(\Pi(\MM)}.
\end{equation*}
Now~\eqref{eq:desired Phi} follows by taking a   $\d\to \big(8\Pi(\MM)/(8\Pi(\MM)+1)\big)^-$ limit of $\Phi_\d$ along a convergent subsequence. 
\end{proof}

\subsection{Strongly bi-Lipschitz from weakly bi-Lipschitz} Here we will prove the following embedding result: 

\begin{theorem}\label{thm:ALN Dk version}  For every $n\in \{3,4,\ldots\}$, every $n$-point metric space $(\MM,d_\MM)$ satisfies
\begin{equation}\label{eq:combine scales}
\cc_2(\MM)\lesssim \max\Bigg\{{\textstyle\sqrt{\log n}}\bigg(\sum_{k=2}^{n-1} \frac{\widehat{\dd}_2^k(\MM)^2}{k(\log k)^2}\bigg)^{\frac12},\dd_2(\MM) \Bigg\}{\textstyle\sqrt{\log\log n}}.
\end{equation}
\end{theorem}

Thanks to Lemma~\ref{lem:sep ext product}, Theorem~\ref{thm:ALN Dk version}  implies Theorem~\ref{thm:sep plux ext gives c}. In fact, this  yields the following estimate:
\begin{align}\label{eq:refine intro product embedding thereom}
\begin{split}
\cc_2(\MM)&\lesssim \max\Bigg\{{\textstyle\sqrt{\log n}}\bigg(\sum_{k=2}^{n-1} \frac{\Pi^k(\MM)^2}{k(\log k)^2}\bigg)^{\frac12},\Pi(\MM) \Bigg\}{\textstyle\sqrt{\log\log n}}\\&\stackrel{\eqref{eq:def DM}}{\le} \max\Bigg\{{\textstyle\sqrt{\log n}}\bigg(\sum_{k=2}^{n-1} \frac{\SEP^k(\MM)^2\ee_{\net}^k(\MM; L_2)^4}{k(\log k)^2}\bigg)^{\frac12},\SEP(\MM)^2\ee_{\net}(\MM; L_2)^4 \Bigg\}{\textstyle\sqrt{\log\log n}},
\end{split}
\end{align}
which is a strengthening of~\eqref{eq:decednt for separataion}, where~\eqref{eq:refine intro product embedding thereom} uses the following notation:
$$
\forall k\in \N,\qquad \Pi^k(\MM)\eqdef \sup_{\substack{\emptyset\neq \sub\subset \MM\\ |\sub|\le k}}\Pi(\sub)\qquad \mathrm{and}\qquad \ee_{\net}^k(\MM; L_2) \eqdef \sup_{\substack{\emptyset\neq \sub\subset \MM\\ |\sub|\le k}}\ee_{\net}(\sub; L_2).
$$
(This adheres to the notational convention  that we have been maintaining throughout by which a superscript indicates  the hereditary version of an invariant of metric spaces which is obtained by considering the  subset of the metric space of a given size for which the invariant in question  is maximal.)

For the sole purpose of proving Theorem~\ref{thm:all >2} (as a consequence of our Theorem~\ref{cor:lp} and the Lipschitz extension theorem of~\cite{NPSS06}), one could use~\cite{ALN08} as a ``black box'' rather than using Theorem~\ref{thm:ALN Dk version}. Specifically,    Theorem~\ref{thm:all >2} follows by fixing $p>2$, using Theorem~\ref{cor:lp} and the estimate $\ee( L_p; L_2)\lesssim \sqrt{p}$ of~\cite{NPSS06} to deduce that $\Pi^k( L_p)\lesssim  p^3\sqrt{\log k}$ for every integer $k\ge 2$, whence $\dd_2^k( L_p)\lesssim p^3\sqrt{\log k}$ by Lemma~\ref{lem:sep ext product}, and now a substitution of this conclusion into~\cite[Theorem~4.1]{ALN08} yields Theorem~\ref{thm:all >2}. 

Thus, those who wish to verify Theorem~\ref{thm:all >2} while appealing to the published literature can stop reading Section~\ref{sec:proof of multiscale} here, and continue reading the present article from Section~\ref{sec:vanilla localized} onwards; that material does not rely on the contents of the rest of Section~\ref{sec:proof of multiscale}.

Notwithstanding the above discussion, it is worthwhile to revisit the ideas of~\cite{ALN08}  so as to derive Theorem~\ref{thm:ALN Dk version} because it is an independently interesting geometric statement that did not appear elsewhere. In comparison to the reasoning in~\cite{ALN08},  the ensuing proof of Theorem~\ref{thm:ALN Dk version} is  more modular and its steps are more general and flexible, so they could be of use for further investigations  elsewhere. 

Observe  that Theorem~\ref{thm:ALN Dk version} implies the following embedding result, which holds for every $\frac12\le \theta\le 1$ and every  $C\ge 1$: 
\begin{equation}\label{eq:theta aln version}
\begin{aligned}
 \max_{k\in \{2,\ldots,n\}} \frac{\widehat{\dd}^k_2(\MM)}{ (\log k)^\theta}\le C
\implies \cc_2(\MM)\lesssim C(\log n)^\theta\min \bigg\{\frac{\sqrt{\log\log n}}{\sqrt{2\theta-1}},\log\log n\bigg\}.
\end{aligned}
\end{equation}
The special case $\theta=1/2$ of~\eqref{eq:theta aln version}  coincides with the case $\e=1/2$ of~\cite[Theorem~4.1]{ALN08}. However, if $1/2<\theta\le 1$ is independent of $n$ or more generally if $\lim_{n\to \infty}(\theta-1/2)\log \log n=\infty$, then~\eqref{eq:theta aln version}   improves asymptotically over~\cite[Theorem~4.1]{ALN08}, though this could also be deduced by a careful inspection of the proof in~\cite{ALN08}.  

 In fact, the ensuing justification of Theorem~\ref{thm:ALN Dk version}  is an involved but quite direct refinement and generalization of the reasoning in~\cite{ALN08}; it is also not entirely self-contained because  it uses~\cite[Theorem~4.5]{ALN08} as a ``black box,'' since we do not have anything novel to add to ~\cite[Theorem~4.5]{ALN08} (there are enhancements of~\cite[Theorem~4.5]{ALN08} that appeared in~\cite{ALN07,ANR24}, but those are not needed for the present  purposes). 

\subsubsection{Weakly bi-Lipschitz embeddings of controlled local growth centers} A key part of the proof of Theorem~\ref{thm:ALN Dk version} studies the geometry of the set of points in a metric space which are centers of balls whose size increases in a controlled manner for radii that belong to a prescribed range. The pertinent object is:

\begin{notation}[controlled local growth centers]\label{notation:G set}
Suppose that $(\MM,d_\MM)$ is a locally finite metric space (thus, all balls in $\MM$ are finite). For $K\ge 1$ and $R\ge r\ge 0$,  let $\cG_{\le K}(r,R)$ denote the set of centers in $\MM$ at which the growth rate of balls from radius $r$ to radius $R$ is at most $K$: 
\begin{equation}\label{eq:local groth set def}
\cG_{\le K}(r,R)=\cG_{\le K}^{d_\MM}(r,R)\eqdef \bigg\{x\in \MM:\ \frac{|B_\MM(x,R)|}{|B_\MM(x,r)|}\le K\bigg\}. 
\end{equation}
\end{notation}
Note in passing that the sets of Notation~\ref{notation:G set} obey the following immediate inclusions: 
\begin{equation}\label{eq:growth sets inclusions}
\forall 1\le K\le K',\ \forall 0\le r\le r'\le R'\le R,\qquad  \cG_{\le K}(r,R)\subset \cG_{\le K'}(r',R').
\end{equation}
The following theorem derives  favorable Euclidean embedding properties of the sets of controlled local growth centers from Notation~\ref{notation:G set}, which refine and generalize results that were obtained in~\cite{ALN08}:

\begin{theorem}\label{thm:single scale on controlled growth} Fix an integer $n\ge 3$ and suppose that $(\MM,d_\MM)$ is an $n$-point metric space. Then, 
\begin{itemize}
\item There is a $1$-Lipschitz function $\psi:\MM\to  L_2$ satisfying the following for every $0<r<R$ and $K\ge 1$: 
\begin{equation}\label{eq:goak GK bourgain}
\forall (x,y)\in \cG_{\le K}(r,R)\times \MM,\qquad d_\MM(x,y)> \frac12 r+\frac32 R\implies \|\psi(x)-\psi(y)\|_{\! L_2}\gtrsim \frac{R-r}{\sqrt{K\log n}}.
\end{equation}
\item There is a universal constant $C\ge 1$ with the following property. Fix $K,D,\beta\ge 1$. Suppose  that for any $\sub\subset \MM$ with $|\sub|\le K$ and any $\Delta>0$ there is a $1$-Lipschitz function $f=f_{\sub,\Delta}:\MM\to  L_2$ satisfying:
\begin{equation}\label{eq:assume single scale in theorem}
\forall x,y\in \sub,\qquad \Delta\le d_\MM(x,y)\le 3\beta\Delta\implies \|f(x)-f(y)\|_{\! L_2}\ge\frac{\Delta}{D}.
\end{equation}
Then, for any $\Delta>0$ and any  $R>r>0$ that satisfy the restrictions
\begin{equation}\label{eq:Rr restrictions}
\Delta\ge 9Dr\qquad\mathrm{and}\qquad   R-r\ge C\beta  \big( \log |\cG_{\le K} (r,R)|\big)\Delta,
\end{equation}
 there is a $1$-Lipschitz function $\phi=\phi_{\Delta,R,r}:\MM\to  L_2$ for which the following property holds:
\begin{equation}\label{eq:x in local growth y general}
\forall (x,y)\in \cG_{\le K}(r,R)\times \MM,\qquad \Delta \le d_\MM(x,y)\le \beta \Delta\implies \|\phi(x)-\phi(y)\|_{\! L_2}\gtrsim \frac{\Delta}{D}.
\end{equation}
\end{itemize}
\end{theorem}

Prior to proving Theorem~\ref{thm:single scale on controlled growth} we will explain how it implies Theorem~\ref{thm:ALN Dk version}:

\begin{proof}[Proof of Theorem~\ref{thm:ALN Dk version} assuming Theorem~\ref{thm:single scale on controlled growth}] Let $C\ge 1$ be the universal constant from (the second part of) Theorem~\ref{thm:single scale on controlled growth}. For the ensuing reasoning, it will be convenient to set the following notation: 
\begin{equation}\label{eq:rR notation ik}
\forall i\in \Z,\  \forall k\in \n,\qquad r_i^k\eqdef \frac{2^i}{27\widehat{\dd}_2^k(\MM)}\qquad \mathrm{and}\qquad R_i^n\eqdef 3C(\log n)2^i. 
\end{equation}

Fix $k\in \n$. By the definition of $\widehat{\dd}_2^k(\MM)$ and  Remark~\ref{rem:2Delta}  for $\beta=6$ and $\alpha=2$, for every subset $\sub\subset \MM$ with $2\le |\sub|\le k$ and every $\Delta>0$ there exists a $1$-Lipschitz function $f=f_{\sub,\Delta}:\MM\to  L_2$ such that 
$$
\forall x,y\in \sub,\qquad \Delta\le d_\MM(x,y)\le 6\Delta\implies \|f(x)-f(y)\|_{\! L_2}\ge\frac{\Delta}{3\widehat{\dd}_2^k(\MM)}.
$$
The assumption of the second part of Theorem~\ref{thm:single scale on controlled growth} therefore holds with the following parameters:
\begin{equation}\label{eq:beta K D coices}
\beta=2\qquad \mathrm{and}\qquad K=k\qquad \mathrm{and}\qquad  D=3\widehat{\dd}_2^k(\MM).
\end{equation}
Consequently, for every $i\in \Z$ there exists a $1$-Lipschitz function $\phi_i^k:\MM\to  L_2$ such that 
\begin{equation}\label{apply second part of theorem}
\forall (x,y)\in \cG_{\le k}\big(r_i^k,R_i^n\big)\times \MM,\qquad 2^i \le d_\MM(x,y)\le 2^{i+1}\implies \|\phi_{i}^k(x)-\phi_i^k(y)\|_{\! L_2}\gtrsim r_i^k,
\end{equation}
This is indeed a valid application of Theorem~\ref{thm:single scale on controlled growth} for the parameter choices
\begin{equation}\label{eq:Delta,rR choices}
\Delta=27\widehat{\dd}_2^k(\MM)2^i\qquad\mathrm{and}\qquad  r=r_i^k\qquad   \mathrm{and}\qquad  R=R^n_i, 
\end{equation}
because $\Delta=9Dr$ by~\eqref{eq:beta K D coices} and~\eqref{eq:Delta,rR choices}, and furthermore since $C\ge 1$ and $|\cG_{\le K} (r,R)|\le |\MM|=n$ we have 
$$
R-r\stackrel{\eqref{eq:rR notation ik}}{\ge} 2C(\log n)2^i\ge 2C\big(\log |\cG_{\le K} (r,R)|\big)2^i\stackrel{\eqref{eq:beta K D coices}\wedge \eqref{eq:Delta,rR choices}}{=}\beta C \big(\log |\cG_{\le K} (r,R)|\big)\Delta. 
$$

Apply~\cite[Theorem~4.5]{ALN08} to $\{\phi_i^k\}_{k\in \Z}$ with the parameters $A= 3C(\log n)/2\ge 1$ and $B=27\widehat{\dd}_2^k(\MM)\lesssim \log n$, while using $B\lesssim  \log n$, which holds by~\cite{bourgain1985lipschitz}. This yields  $\f^k:\MM\to  L_2$ that has the following two properties: 
\begin{equation}\label{eq:lip requirement for phik}
\|\f^k\|_{\Lip(\MM; L_2)}\lesssim {\textstyle \sqrt{(\log n)\log\log n}},
\end{equation}
and for every $x,y\in \MM$ and every $i\in \Z$ we have:
\begin{equation}\label{eq:our phik lower}
\begin{aligned}
\|\f^k(x)-\f^k(y)\|_{\! L_2}^2
\gtrsim 
\bigg\lfloor\log_2 \frac{|B_\MM(x,R_i^n)|}{|B_\MM(x,r_i^k)|}\bigg\rfloor
\min \bigg\{\frac{2^{2i}}{\widehat{\dd}_2^k(\MM)^2}, \|\phi_i^k(x)-\phi_i^k(y)\|_{\! L_2}^2\bigg\}.
\end{aligned}
\end{equation}

For every distinct $x,y\in \MM$ define $i(x,y)\in \Z$ by the requirement $2^{i(x,y)}\le d_\MM(x,y)<2^{i(x,y)+1}$, i.e., 
\begin{equation}\label{eq:def ixy}
i(x,y)\eqdef \big\lfloor \log_2d_\MM(x,y)\big\rfloor.   
\end{equation}
We will also set $i(x,x)=-\infty$ for every $x\in \MM$. As $|\MM|=n$ the definition~\eqref{eq:local groth set def} of $\cG_{\le n}(r,R)$  implies that  $\cG_{\le n}(r,R)=\MM$ for every $R\ge r\ge 0$. Hence, the following index $\ck(x,y)\in \n$ is well-defined:  
\begin{equation}\label{eq:def ck}
\ck(x,y)\eqdef \min \Big\{k\in \{1,\ldots,n\}:\ x\in \bigcap_{\ell=k}^n \cG_{\le \ell}\big(r^\ell_{i(x,y)},R_{i(x,y)}^n\big)\Big\}.
\end{equation}
For every $x\in \MM$, the  above convention $i(x,x)=-\infty$ is consistent with setting $\ck(x,x)=1$. We claim that:
\begin{equation}\label{eq:our phik lower'}
\forall x,y\in \MM,\ \forall k\in \{\ck(x,y),\ck(x,y)+1,\ldots,n\},\qquad \|\f^k(x)-\f^k(y)\|_{\! L_2}\gtrsim \frac{\sqrt{\lfloor\log \ck(x,y)\rfloor}}{\widehat{\dd}_2^k(\MM)}d_\MM(x,y).
\end{equation}

To justify~\eqref{eq:our phik lower'}, fix $x,y\in \MM$ and  $k\in \{\ck(x,y),\ck(x,y)+1,\ldots,n\}$. If $\ck(x,y)\in \{1,2\}$, then the right hand side of the inequality in~\eqref{eq:our phik lower'}  vanishes, so~\eqref{eq:our phik lower'} holds vacuously. We can therefore assume that $\ck(x,y)\in \{3,4,\ldots,n\}$. As $k\ge \ck(x,y)$,  the definition~\eqref{eq:def ck} of $\ck(x,y)$ and the definition~\eqref{eq:def ixy} of $i(x,y)$ show that 
\begin{equation}\label{eq:x in rikght gfrowth set}
x\in \cG_{\le k}\big(r^k_{i(x,y)},R_{i(x,y)}^n\big)\qquad \mathrm{and}\qquad 2^{i(x,y)}\le d_\MM(x,y)< 2^{i(x,y)+1}.
\end{equation}
Consequently,
\begin{equation}\label{eq:phiik lower at index}
\|\phi_{i(x,y)}^k(x)-\phi_{i(x,y)}^k(y)\|_{\! L_2} \stackrel{ \eqref{eq:x in rikght gfrowth set}\wedge  \eqref{apply second part of theorem}}{\gtrsim}  r_{i(x,y)}^k\stackrel{\eqref{apply second part of theorem}}{\asymp} \frac{d_\MM(x,y)}{\widehat{\dd}^k_2(\MM)}.
\end{equation}
At the same time, by the minimality of $\ck(x,y)$ per~\eqref{eq:def ck}, we have 
$$
x\notin \cG_{\le \ck(x,y)-1}\big(r^{\ck(x,y)-1}_{i(x,y)},R_{i(x,y)}^n\big),
$$
which by~\eqref{eq:local groth set def} is equivalent to the following lower bound on the growth rate of balls centered at $x$:
\begin{equation}\label{eq:ck minimum}
\frac{\big|B_\MM\big(x,R_{i(x,y)}^n\big)\big|}{\big|B_\MM\big(x,r_{i(x,y)}^{\ck(x,y)-1}\big)\big|}>\ck(x,y)-1. 
\end{equation}
Since $\widehat{\dd}_2^k(\MM)\ge \widehat{\dd}_2^{\ck(x,y)-1}(\MM)$ we have $r^{k}_{i(x,y)}\le r^{\ck(x,y)-1}_{i(x,y)}$ by~\eqref{eq:rR notation ik}, and therefore
\begin{equation}\label{eq:ratio of balls decrease}
\big|B_\MM\big(x,r_{i(x,y)}^{k}\big)\big|\le \big|B_\MM\big(x,r_{i(x,y)}^{\ck(x,y)-1}\big)\big|.
\end{equation}
Consequently,
\begin{equation}\label{eq:bound ratio of balls by critical indiex}
\begin{aligned}
\bigg\lfloor\log_2 \frac{\big|B_\MM\big(x,R_{i(x,y)}^n\big)\big|}{\big|B_\MM\big(x,r_{i(x,y)}^k\big)\big|}\bigg\rfloor
&\stackrel{\eqref{eq:ratio of balls decrease}}{\ge}
\bigg\lfloor\log_2 \frac{\big|B_\MM\big(x,R_{i(x,y)}^n\big)\big|}{\big|B_\MM\big(x,r_{i(x,y)}^{\ck(x,y)-1}\big)\big|}\bigg\rfloor\\
&\stackrel{\eqref{eq:ck minimum}}{\ge}
\big\lfloor \log_2(\ck(x,y)-1)\big\rfloor
\asymp \log \ck(x,y).
\end{aligned}
\end{equation}
where the last step of~\eqref{eq:bound ratio of balls by critical indiex} is valid since we are now working under the assumption $\ck(x,y)\ge 3$. The desired assertion~\eqref{eq:our phik lower'} follows by substituting~\eqref{eq:phiik lower at index} and~\eqref{eq:bound ratio of balls by critical indiex} into~\eqref{eq:our phik lower} while using $2^{i(x,y)}\asymp d_\MM(x,y)$. 

Next, observe that the function $\psi:\MM\to  L_2$  from the first part of Theorem~\ref{thm:single scale on controlled growth} satisfies:
\begin{equation}\label{eq:psi properties1}
\begin{aligned}
\|\psi\|_{\Lip(\MM; L_2)}\le 1,
\end{aligned}
\end{equation}
and:
\begin{equation}\label{eq:psi properties2}
\begin{aligned}
\forall x,y\in \MM,\qquad 
\|\psi(x)-\psi(y)\|_{\! L_2}\gtrsim \frac{1}{\sqrt{\ck(x,y)\log n}}d_\MM(x,y).
\end{aligned}
\end{equation}
Indeed, \eqref{eq:psi properties1} is the upper bound on the Lipschitz constant of $\psi$ from statement of Theorem~\ref{thm:single scale on controlled growth}. For~\eqref{eq:psi properties2},  fix distinct $x,y\in \MM$ and record  the following (very) crude bounds:
\begin{equation}\label{eq:crude for psi}
\begin{aligned}
R^n_{i(x,y)} \stackrel{\eqref{eq:rR notation ik}}{>} 2^{i(x,y)-1}\qquad\mathrm{and}\qquad 
d_\MM(x,y)\stackrel{\eqref{eq:def ixy}}{\ge} 2^{i(x,y)}
\stackrel{\eqref{eq:rR notation ik}}{>} \frac12r^{\ck(x,y)}_{i(x,y)}+\frac32 2^{i(x,y)-1}.
\end{aligned}
\end{equation}
Thus,
\begin{equation}\label{eq: x also in smaller R}
x\in \cG_{\le \ck(x,y)}\big(r^{\ck(x,y)}_{i(x,y)},R_{i(x,y)}^n\big)\subset \cG_{\le \ck(x,y)}\big(r^{\ck(x,y)}_{i(x,y)},2^{i(x,y)-1}\big),
\end{equation}
where the first step of~\eqref{eq: x also in smaller R} follows from the case $k=\ck(x,y)$ of~\eqref{eq:x in rikght gfrowth set}  and the second step of~\eqref{eq: x also in smaller R} follows from the first part of~\eqref{eq:crude for psi} and the (trivial) monotonicity property~\eqref{eq:growth sets inclusions} of the sets of controlled local growth centers.  We can thereofre complete the justification of~\eqref{eq:psi properties2} as follows:
\begin{equation}\label{eq:use GK bourgain}
\|\psi(x)-\psi(y)\|_{\! L_2}\gtrsim  \frac{2^{i(x,y)-1}-r^{\ck(x,y)}_{i(x,y)}}{\sqrt{\ck(x,y)\log n}}\asymp \frac{2^{i(x,y)-1}}{\sqrt{\ck(x,y)\log n}}\asymp \frac{d_\MM(x,y)}{\sqrt{\ck(x,y)\log n}},
\end{equation}
where the first step of~\eqref{eq:use GK bourgain} uses the first part~\eqref{eq:goak GK bourgain} of Theorem~\ref{thm:single scale on controlled growth} with $r=r^{\ck(x,y)}_{i(x,y)}$, $R=2^{i(x,y)-1}$,  $K=\ck(x,y)$, which is valid   thanks to~\eqref{eq: x also in smaller R} and the second part of~\eqref{eq:crude for psi}.

If $v_1,\ldots,v_n$ is an orthonormal system in $ L_2$, then define $f:\MM\to L_2\otimes  L_2\cong L_2$ by setting for every $x\in \MM$:
\begin{equation}\label{eq:def ourt final aln f}
 f(x)\eqdef {\textstyle\sqrt{(\log n)\log\log n}}\cdot \psi(x)\otimes v_1+\sum_{k=2}^{n-1} \frac{\widehat{\dd}_2^k(\MM)}{\sqrt{k}\log k}\f^k(x)\otimes v_k + \frac{\widehat{\dd}_2(\MM)}{\sqrt{\log n}}\f^n(x)\otimes v_n.
\end{equation}
The orthonormality of $v_1,\ldots,v_n$  and the the bounds on the Lipschitz constants in~\eqref{eq:lip requirement for phik} and~\eqref{eq:psi properties1} give 
\begin{equation}\label{eq:or multiscale f lip constant}
\|f\|_{\Lip(\MM; L_2\otimes  L_2)}\lesssim \max\Bigg\{{\textstyle\sqrt{\log n}}\bigg(\sum_{k=2}^{n-1} \frac{\widehat{\dd}_2^k(\MM)^2}{k(\log k)^2}\bigg)^{\frac12},\dd_2(\MM) \Bigg\}{\textstyle\sqrt{\log\log n}}. 
\end{equation}
Furthermore, every distinct $x,y\in \MM$ satisfy the following lower bound:
\begin{align}\label{eq:lower multiscale ck}
\begin{split}
&\!\!\!\!\!\!\!\!\!\!\!\!\!\!\frac{\|f(x)-f(y)\|_{ L_2\otimes  L_2}}{d_\MM(x,y)}\\&\stackrel{\eqref{eq:def ourt final aln f}\wedge \eqref{eq:our phik lower'}\wedge \eqref{eq:psi properties2}}{\gtrsim} \bigg(\frac{\log \log n}{\ck(x,y)}+\lfloor \log \ck(x,y)\rfloor\sum_{k=\ck(x,y)}^{n-1} \frac{1}{k(\log k)^2}+\frac{\lfloor \log \ck(x,y)\rfloor}{\log n}\bigg)^{\frac12}\\
&\qquad \quad \ge \bigg(\frac{\log \log n}{\ck(x,y)}+\lfloor \log \ck(x,y)\rfloor\int_{\max\{\ck(x,y),2\}}^n\frac{\ud s}{s(\log s)^2}+\frac{\lfloor \log \ck(x,y)\rfloor}{\log n}\bigg)^\frac12\\
&\qquad \ \ \ \ = \bigg(\frac{\log \log n}{\ck(x,y)}+\lfloor \log \ck(x,y)\rfloor\bigg(\frac{1}{\log \max\{\ck(x,y),2\} }-\frac{1}{\log n}\bigg)+\frac{\lfloor \log \ck(x,y)\rfloor}{\log n}\bigg)^\frac12\gtrsim 1,
\end{split}
\end{align}
where the last step of~\eqref{eq:lower multiscale ck} is  verified by considering separately the cases $1\le \ck(x,y)\le \max\{\log\log n,3\}$, $\max\{\log\log n,3\}\le \ck(x,y)\le \sqrt{n}$ and $\sqrt{n}\le \ck(x,y)\le n$. The desired distortion bound~\eqref{eq:combine scales} of Theorem~\ref{thm:ALN Dk version} follows by combining~\eqref{eq:or multiscale f lip constant} and~\eqref{eq:lower multiscale ck} for the embedding $f$ of $\MM$ into Hilbert space. 
\end{proof}

Passing to the proof of Theorem~\ref{thm:single scale on controlled growth}, we will first justify its first part by following the reasoning in the proof of~\cite[Claim~4.6]{ALN08} as well as its predecessors~\cite[Lemma~3.4]{mendel2004euclidean}, \cite[Theorem~2]{ABN15} and~\cite[Lemma~2.9]{CMM10}, all of which are an elaboration of the important embedding technique that was introduced in~\cite{bourgain1985lipschitz}.

Given a metric space $(\MM,d_\MM)$, the distance of a point $x\in \MM$ to a nonempty subset $\sub$ of $\MM$ is denoted as usual by $d_\MM(x,\sub)=\inf\{d_\MM(x,y):\ y\in \sub\}$.  It is worthwhile to state separately the following  observation:

\begin{observation}\label{lem:bernoulli}  Let $(\MM,d_\MM)$ be a finite metric space. Given $0\le \fp\le 1$, denote by $\cZ_\fp$  the ($\fp$-Bernoulli) random subset of $\MM$ that is obtained by including independently each $x\in \MM$ in $\cZ$ with probability $\fp$, i.e., 
\begin{equation}\label{eq:bernouli formula for prob}
\forall \sub\subset \MM,\qquad \prob[\cZ_\fp=\sub]=\fp^{|\sub|}(1-\fp)^{n-|\sub|}. 
\end{equation}
Then, for every $0<r<R$ and every $x,y\in \MM$ such that $d_\MM(x,y)> \frac12 r+\frac32 R$, we have 
\begin{align}\label{eq:bernoulli goal}
\begin{split}
\prob \bigg[\cZ_\fp\neq\emptyset \quad\mathrm{and}\quad 
\big|&d_\MM(x,\cZ_\fp)-d_\MM(y,\cZ_\fp)\big|> \frac{R-r}{2}\bigg]\\
&\quad\ge \min\Big\{ (1-\fp)^{|B_\MM(x,R)|}, 1-(1-\fp)^{|B_\MM(x,r)|}\Big\}.
\end{split}
\end{align}
\end{observation}

\begin{proof} Consider the events $E$ and $F$ that are defined as follows:
\begin{align}\label{eq:def EF}
\left\{\begin{array}{ll}
E\eqdef \left\{d_\MM(x,\cZ_\fp)>R\quad \mathrm{and}\quad d_\MM(y,\cZ_\fp)\le \frac{r+R}{2}\right\},\\
F\eqdef \left\{d_\MM(x,\cZ_\fp)\le r\quad \mathrm{and}\quad d_\MM(y,\cZ_\fp)> \frac{r+R}{2} \right\}.
\end{array}\right.
\end{align}
By definition,  $E\cap F=\emptyset$ and $E\cup F\subset \{\cZ_\fp\neq\emptyset \ \wedge\ |d_\MM(x,\cZ_\fp)-d_\MM(y,\cZ_\fp)|> (R-r)/2\}$, so we have 
\begin{equation}\label{eq:at least EF}
\prob \left[\cZ_\fp\neq\emptyset \quad\mathrm{and}\quad |d_\MM(x,\cZ_\fp)-d_\MM(y,\cZ_\fp)|> \frac{R-r}{2}\right]\ge \prob[E]+\prob[F].
\end{equation}

Note that the event $\{d_\MM(x,\cZ_\fp)>R\}$ occurs if and only if $\cZ_\fp$ contains no point from  $B_\MM(x,R)$, and the event $\{d_\MM(y,\cZ_\fp)\le (r+R)/2\}$ occurs if and only if $\cZ_\fp$ contains at least one point from  $B_\MM(y,(r+R)/2)$. By the triangle inequality for $d_\MM$, the assumption $d_\MM(x,y)>r/2+3R/2=R+(r+R)/2$ of Observation~\ref{lem:bernoulli} implies that $B_\MM(x,R)\cap B_\MM(y,(r+R)/2)=\emptyset$, so  $\{d_\MM(x,\cZ_\fp)>R\}$ and $\{d_\MM(y,\cZ_\fp)\le (r+R)/2\}$ are independent events,  thanks to the definition of $\cZ_\fp$. Consequently, 
\begin{align}\label{eq:probE}
\begin{split}
\prob[E]\stackrel{\eqref{eq:def EF}}{=}\prob \big[\cZ_\fp\cap B_\MM(x,R)&=\emptyset\big] \prob\left[d_\MM(y,\cZ_\fp)\le \frac{r+R}{2}\right]\\&\stackrel{\eqref{eq:bernouli formula for prob}}{=}(1-\fp)^{|B_\MM(x,R)|}\prob\left[d_\MM(y,\cZ_\fp)\le \frac{r+R}{2}\right].
\end{split}
\end{align}
In the same vein, the events $\{ d_\MM(x,\cZ_\fp)\le r\}$ and  $\{d_\MM(y,\cZ_\fp)>(r+R)/2\}$ coincide with, respectively, the events $\{\cZ_\fp\cap B_\MM(x,r)\neq\emptyset\}$ and   $\{\cZ_\fp\cap B_\MM(y,(r+R)/2)=\emptyset\}$. The latter two events are independent because $B_\MM(x,r)\cap B_\MM(y,(r+R)/2)\subset B_\MM(x,R)\cap B_\MM(y,(r+R)/2)=\emptyset$, so we have:
\begin{align}\label{eq:probF}
\begin{split}
\prob[F]\stackrel{\eqref{eq:def EF}}{=}\prob \big[\cZ_\fp\cap B_\MM(x,r)&\neq\emptyset\big] \prob\left[d_\MM(y,\cZ_\fp)>\frac{r+R}{2} \right]\\&\stackrel{\eqref{eq:bernouli formula for prob}}{=}\left(1-(1-\fp)^{|B_\MM(x,r)|}\right)\prob\left[d_\MM(y,\cZ_\fp)>\frac{r+R}{2} \right].
\end{split}
\end{align}
The desired bound~\eqref{eq:bernoulli goal} is now justified as follows:
\begin{align*}
&\prob \bigg[\cZ_\fp\neq\emptyset \quad\mathrm{and}\quad |d_\MM(x,\cZ_\fp)-d_\MM(y,\cZ_\fp)|> \frac{R-r}{2}\bigg]\\&\stackrel{\eqref{eq:at least EF}\wedge \eqref{eq:probE}\wedge\eqref{eq:probF}}{\ge} (1-\fp)^{|B_\MM(x,R)|}\prob\left[d_\MM(y,\cZ_\fp)\le \frac{r+R}{2}\right]\\&\qquad \qquad \qquad \qquad \qquad \qquad +\left(1-(1-\fp)^{|B_\MM(x,r)|}\right)\prob\left[d_\MM(y,\cZ_\fp)>\frac{r+R}{2} \right]\\
&\qquad\quad\   \ge \min\Big\{ (1-\fp)^{|B_\MM(x,R)|}, 1-(1-\fp)^{|B_\MM(x,r)|}\Big\}\\&\qquad\qquad\qquad\qquad \qquad\qquad  \cdot  \left(\prob\left[d_\MM(y,\cZ_\fp)\le\frac{r+R}{2} \right]+\prob\left[d_\MM(y,\cZ_\fp)>\frac{r+R}{2} \right]\right)\\&\ \ \ \ \ \ \ \ \ \ \ \ \ =\min\Big\{ (1-\fp)^{|B_\MM(x,R)|}, 1-(1-\fp)^{|B_\MM(x,r)|}\Big\}. \qedhere
\end{align*}
\end{proof}

For the of proof of the first part of Theorem~\ref{thm:single scale on controlled growth}, it will  be notationally beneficial to use the convention that the distance of a point $x\in \MM$ of a metric space $(\MM,d_\MM)$ to the empty set equals $\infty$ (in our setting $\MM$ is finite, so the same purpose will be served by defining $d_\MM(x,\emptyset)$ to be any fixed number that is strictly larger than the diameter $\diam(\MM)=\sup\{d_\MM(y,z):\ y,z\in \MM\}$ of $\MM$, say, $d_\MM(x,\emptyset)=\diam(\MM)+1$).

\begin{proof}[Proof of the first part of Theorem~\ref{thm:single scale on controlled growth}] For each $i\in \{1,\ldots,\lceil \log n\rceil\}$ let $\cZ_{e^{-i}}$ be the random subset of $\MM$ from Observation~\ref{lem:bernoulli} with $\fp=e^{-i}$. Denote the probability space on which these sets are defined by $(\Omega,\prob)$; it can be the corresponding product space, though they need not be independent for the ensuing reasoning. Let $v_1,v_2,\ldots$ be an orthonormal basis of $ L_2$ and define a $1$-Lipschitz function $f:\MM\to L_2(\prob; L_2)\cong  L_2$ by
\begin{equation}\label{eq:def bourgain embedding}
\forall(x,\omega)\in \MM\times \Omega,\qquad f(x)(\omega)\eqdef \frac{1}{\sqrt{\lceil\log n\rceil}}\sum_{i=1}^{\lceil \log n\rceil} \min\big\{d_\MM\big(x,\cZ_{e^{-i}}(\omega)\big),\diam(\MM)\big\}v_i. 
\end{equation}

For $x\in \MM$ and $R>0$ let $i_R(x)\in \{1,\ldots,\lceil \log n\rceil\}$ be such that $e^{i_R(x)-1}\le |B_\MM(x,R)|<e^{i_R(x)}$. Then, for every  $K\ge 1$, every $0<r<R$, and every $x\in \cG_{\le K}(r,R)$ we have:
 \begin{equation}\label{eq:i(x) bounds}
\begin{aligned}
e^{i_R(x)}>|B_\MM(x,R)|\ge |B_\MM(x,r)|
\stackrel{\eqref{eq:local groth set def}}{\ge} \frac{|B_\MM(x,R)|}{K}\ge \frac{e^{i_R(x)-1}}{K}.
\end{aligned}
\end{equation}
Using Observation~\ref{lem:bernoulli}, the following probabilistic estimate holds for every $x\in \cG_{\le K}(r,R)$ and every $y\in \MM$ that satisfy $d_\MM(x,y)>r/2+3R/2$: 
\begin{align}\label{eq:probabilistic version bourgain growth}
\begin{split}
\prob\bigg[&\big|\min\{d_\MM(x,\cZ_{e^{-i_R(x)}}),\diam(\MM)\}
-\min\{d_\MM(y,\cZ_{e^{-i_R(x)}}),\diam(\MM)\}\big|> \frac{R-r}{2}\bigg]\\
&\stackrel{\eqref{eq:bernoulli goal}\wedge \eqref{eq:i(x) bounds}}{\ge}
\min\left\{ \left(1-e^{-i_R(x)}\right)^{e^{i_R(x)}},
1-\left(1-e^{-i_R(x)}\right)^{\frac{e^{i_R(x)-1}}{K}}\right\}\asymp \frac{1}{K}.
\end{split}
\end{align}
This directly implies the desired estimate~\eqref{eq:goak GK bourgain} as follows:
\begin{align*}
\|f(x)-f(y)\|_{\!L_2(\prob; L_2)}
&\stackrel{\eqref{eq:def bourgain embedding}}{=}
\frac{1}{\sqrt{\lceil\log n\rceil}}
\bigg(\sum_{j=1}^{\lceil \log n\rceil} \int_\Omega
\Big|
\min\big\{d_\MM\big(x,\cZ_{e^{-j}}(\omega)\big),\diam(\MM)\big\}\\
&\qquad\qquad
-\min\big\{d_\MM\big(y,\cZ_{e^{-j}}(\omega)\big),\diam(\MM)\big\}
\Big|^2 \ud\prob(\omega)\bigg)^{\frac12}\\
&\ge \frac{1}{\sqrt{\lceil\log n\rceil}}
\bigg(\int_\Omega
\Big|
\min\big\{d_\MM\big(x,\cZ_{e^{-i_R(x)}}(\omega)\big),\diam(\MM)\big\}\\
&\qquad\qquad
-\min\big\{d_\MM\big(y,\cZ_{e^{-i_R(x)}}(\omega)\big),\diam(\MM)\big\}
\Big|^2 \ud\prob(\omega)\bigg)^{\frac12}\\
&\stackrel{\eqref{eq:probabilistic version bourgain growth}}{\ge}
\frac{R-r}{\sqrt{K\log n}}.\qedhere
\end{align*}
\end{proof}

It remains to prove the second part of Theorem~\ref{thm:single scale on controlled growth}. We will start with: 
\begin{lemma}\label{lem:rR} Fix $K,D,\beta\ge 1$. Let $(\MM,d_\MM)$ be a locally finite metric space such that for every $\sub\subset \MM$ with $|\sub|\le K$ and every $\Delta>0$ there exists a $1$-Lipschitz function $f=f_{d_\MM,\sub,\Delta}:\MM\to  L_2$ that satisfies:
\begin{equation}\label{eq:repeat single scale assumption for U}
\forall x,y\in \sub,\qquad \Delta\le d_\MM(x,y)\le 3\beta\Delta\implies \|f(x)-f(y)\|_{\! L_2}\ge\frac{\Delta}{D}.
\end{equation}
Then, for every finite nonempty subset $\cU$ of $\MM$,  any $\Delta>0$ and any  $R>r>0$ that satisfy the restrictions
\begin{equation}\label{eq:Rr Delta restrictions diam U}
\Delta\ge 9Dr\qquad\mathrm{and}\qquad   R-r\ge \diam_\MM(\cU),
\end{equation}
there exists a $1$-Lipschitz function $\phi=\phi_{\cU,r,R,\Delta,d_\MM}:\MM\to  L_2$ such that 
\begin{equation}\label{eq:goal beta/2}
\forall (x,y)\in \big(\cU\cap \cG_{\le K}(r,R)\big)\times \MM,\qquad \Delta\le d_\MM(x,y)\le \beta\Delta\implies \|\phi(x)-\phi(y)\|_{\! L_2}\gtrsim \frac{\Delta}{D}.
\end{equation}
\end{lemma}

Due to the similarity of Lemma~\ref{lem:rR} and the second part of Theorem~\ref{thm:single scale on controlled growth}, it is helpful to spell out how they differ. Assumption~\eqref{eq:repeat single scale assumption for U} of Lemma~\ref{lem:rR} coincides with assumption~\eqref{eq:assume single scale in theorem} of Theorem~\ref{thm:single scale on controlled growth}. Conclusion~\eqref{eq:goal beta/2} of Lemma~\ref{lem:rR} is weaker than the desired conclusion~\eqref{eq:x in local growth y general} of  Theorem~\ref{thm:single scale on controlled growth}:  It imposes a further restriction that the point $x$ belongs not only to the set $\cG_{\le K}(r,R)$ of centers whose local growth at the input radii $R>r>0$  is at most $K$,  but also to an arbitrary but fixed bounded nonempty subset $\cU$ of $\MM$; for this,  Lemma~\ref{lem:rR} requires in~\eqref{eq:Rr Delta restrictions diam U} that  those radii  are separated by  $\diam_\MM(\cU)$, while Theorem~\ref{thm:single scale on controlled growth}  asks in~\eqref{eq:Rr restrictions} for them to be separated by a quantity that involves both the given scale $\Delta>0$ and the size of  $\cG_{\le K}(r,R)$.

\begin{proof}[Proof of Lemma~\ref{lem:rR}] Let $\NN$ be  a subset of $\cU\cap \cG_{\le K}(r,R)$ with the property that $d_\MM(a,b)>2r$ for all distinct $a,b\in \NN$, and furthermore $\NN$ is maximal with respect to inclusion among all the subsets of $\cU\cap \cG_{\le K}(r,R)$  that have this property. By the triangle inequality for $d_\MM$, the balls $\{B_\MM(a,r)\}_{a\in \NN}$ are disjoint, so we have 
\begin{equation}\label{eq:poacking}
|\NN|\min_{a\in \NN} |B_\MM(a,r)|\le \sum_{a\in \NN} |B_\MM(a,r)|= \Big|\bigcup_{a\in \NN} B_\MM(a,r)\Big|.
\end{equation}

Fix $a_{\min} \in \NN$ such that $|B_\MM(a_{\min},r)|=\min_{a\in \NN} |B_\MM(a,r)|$.  By the triangle inequality for $d_\MM$, since $\NN\subset \cU$ we have $B_\MM(a,r)\subset B_\MM(a_{\min},\diam_\MM(\cU)+r)\subset B_\MM(a_{\min},R)$ for every $a\in \NN$, as $R-r\ge \diam_\MM(\cU)$. Thanks to~\eqref{eq:poacking}, this implies  $|\NN|\le |B_\MM(a_{\min},R)|/|B_\MM(a_{\min},r)|$. As $a_{\min}\in \NN\subset \cG_{\le K}(r,R)$, we deduce that $|\NN|\le K$. 

By the assumption of Lemma~\ref{lem:rR}  applied to $\sub=\NN$, there exists a $1$-Lipschitz function $f:\MM\to  L_2$ that satisfies the following condition for every $ a,b\in \NN$:
\begin{equation}\label{eq:use single scale for Dela0}
\Delta_0\le d_\MM(a,b)\le 3\beta\Delta_0\implies \|f(a)-f(b)\|_{\! L_2}\ge\frac{\Delta_0}{D},\quad\mathrm{where}\quad  \Delta_0\eqdef \frac{19D-1}{19D}\Delta-4r, 
\end{equation}
which is valid since $\Delta\ge 9Dr> 76Dr/(19D-1)$, as $D\ge 1$, so the parameter $\Delta_0$ in~\eqref{eq:use single scale for Dela0} is indeed positive. Because $\MM$ is locally finite, $\phi(\MM)$ is a countable subset of $ L_2$, so by replacing $ L_2$ by $ L_2\oplus \R$, we may assume that there exists a vector $v\in f(\MM)^\perp\subset  L_2$ that satisfies $\|v\|_{\! L_2}=1$ and $v$ is orthogonal to $f(\MM)$. Since $f$ is $1$-Lipschitz, and by the triangle inequality for $d_\MM$ also the function $(x\in \MM)\mapsto d_\MM(x,\NN)$ is $1$-Lipschitz, the aforementioned orthogonality implies that the following  function $\phi:\MM\to L_2$ is $1$-Lipschitz:
$$
\forall x\in \MM,\qquad \phi(x)\eqdef \frac{1}{\sqrt{2}}f(x)+  \frac{d_\MM(x,\NN)}{\sqrt{2}}v.
$$
We then have the following lower bound: 
\begin{equation}\label{eq:max version of lower bound}
\forall x,y\in \MM,\qquad \|\phi(x)-\phi(y)\|_{\! L_2}\gtrsim \max \big\{\|f(x)-f(y)\|_{\! L_2}, |d_\MM(x,\NN)-d_\MM(y,\NN)|\big\}.
\end{equation}

Fix $x\in \cU\cap \cG_{\le K}(r,R)$. By the maximality of $\NN$ with respect to inclusion, there must exist $a_x\in \NN$ such that $d_\MM(x,a_x)\le 2r$, i.e., we have $d_\MM(x,\NN)\le 2r$.   Consequently, if $y\in \MM$ satisfies  
\begin{equation}\label{eq:lower bound on distance to NN}
d_\MM(y,\NN)\ge 2r+\frac{\Delta}{19D}, 
\end{equation}
then it follows from~\eqref{eq:max version of lower bound} that $\|\phi(x)-\phi(y)\|_{\! L_2}\gtrsim \Delta/D$, which is the desired conclusion of~\eqref{eq:lower bound on distance to NN}.

The above simple preparatory reasoning demonstrates that the proof of Lemma~\ref{lem:rR} will be complete if we will prove that  $\|\phi(x)-\phi(y)\|_{\! L_2}\gtrsim \Delta/D$ whenever $x,y\in \MM$ have  the following properties:
\begin{equation}\label{eq:new assumptions axay}
\left\{ \begin{array}{ll}
\Delta\le d_\MM(x,y)\le \beta\Delta,\\ \exists a_x,a_y\in \NN,\quad d_\MM(x,a_x)\le 2r\quad \mathrm{and}\quad d_\MM(y,a_y)<2r+\frac{\Delta}{19D}. 
\end{array}\right.
\end{equation}
This is indeed the case because by the triangle inequality for $d_\MM$ we have 
\begin{equation}\label{eq:bigger than Delta0}
d_\MM(a_x,a_y)\ge d_\MM(x,y)-d_\MM(x,a_x)-d_\MM(y,a_y)\stackrel{\eqref{eq:new assumptions axay}}{\ge} \Delta-4r-\frac{\Delta}{19D}\stackrel{\eqref{eq:use single scale for Dela0}}{=}\Delta_0,
\end{equation}
and
\begin{equation}\label{eq:upper beta bound distance}
d_\MM(a_x,a_y)\le d_\MM(x,y)+d_\MM(x,a_y)+d_\MM(y,a_y)\stackrel{\eqref{eq:new assumptions axay}}{\le}  \beta \Delta+4r+\frac{\Delta}{19D}\le 3\beta\Delta_0,
\end{equation}
where checking that the last step of~\eqref{eq:upper beta bound distance} is valid is straightforward using the definition~\eqref{eq:use single scale for Dela0} of $\Delta_0$ together with the assumptions $\beta\ge 1$ and $\Delta\ge 9Dr\ge 9r$. Thanks to~\eqref{eq:bigger than Delta0} and~\eqref{eq:upper beta bound distance} we may use~\eqref{eq:use single scale for Dela0} with $a=a_x$ and $b=a_y$ in combination with~\eqref{eq:max version of lower bound} to conclude that    the following lower bound holds:
\begin{align}\label{eq:to finish rR lemma}
\begin{split}
\|\phi(x)-\phi(y)\|_{ L_2}&\stackrel{\eqref{eq:max version of lower bound} }{\gtrsim} \|f(x)-f(y)\|_{ L_2}\\&\ \ \ \ge \|f(a_x)-f(a_y)\|_{ L_2}-\|f(x)-f(a_x)\|_{ L_2}-\|f(y)-f(a_y)\|_{ L_2}\\&\stackrel{\eqref{eq:use single scale for Dela0}}{\ge} \frac{\Delta_0}{D}-d_\MM(x,a_x)-d_\MM(y,a_y)\\&\!\!\!\!\!\!\!\!\!\!\!\stackrel{\eqref{eq:use single scale for Dela0}\wedge \eqref{eq:new assumptions axay}}{\ge}   \frac{\frac{19D-1}{19D}\Delta-4r}{D}-4r-\frac{\Delta}{19D}\\&\ \ \ \gtrsim \frac{\Delta}{D},
\end{split}
\end{align}
where the second step of~\eqref{eq:to finish rR lemma} uses the triangle inequality in $ L_2$, the third step of~\eqref{eq:to finish rR lemma} uses the fact that $f$ is $1$-Lipschitz, and the final step of~\eqref{eq:to finish rR lemma} is straightforward to verify using $r\le \Delta/(9D)$ and $D\ge 1$. 
\end{proof}

Lemma~\ref{lemm:get suit of functions} below provides a link between Lemma~\ref{lem:rR} and Theorem~\ref{thm:single scale on controlled growth}; it is a spatial localization principle for Euclidean single scale embeddings that  does not involve the controlled local growth centers.  The term ``spatial localization'' is used here for the following reason: The assumption of Lemma~\ref{lemm:get suit of functions} per~\eqref{eq:cX D version} is that ``one scale at a time'' embeddings exist for arbitrary subsets of small diameter, and its conclusion per~\eqref{eq:lower bound on suit in lemma} is that such an embedding guarantee is possible provided one of the points is close to a subset whose cardinality is sufficiently small without any bound on the spacial size (diameter) of that subset.  

\begin{lemma}\label{lemm:get suit of functions} There is  a universal constant $C\ge 1$ with the following property. Suppose that $\beta,D\ge 1$ and $\cd,\Delta>0$. Let $(\MM,d_\MM)$ be a  metric space and fix  $\XX,\YY\subset \MM$. Assume that for every $\emptyset\neq \cU\subset \MM$ with $\diam_\MM(\cU)\le \cd$  there is a $1$-Lipschitz function $f_\cU:\cU\to  L_2$ satisfying 
\begin{equation}\label{eq:cX D version}
\forall (x,y)\in (\cU\cap\XX)\times (\cU\cap \YY),\qquad  \Delta\le d_\MM(x,y)\le \beta \Delta\implies \|f_\cU(x)-f_\cU(y)\|_{\! L_2}\ge \frac{\Delta}{D}.
\end{equation}
Then, for any $\sub\subset \MM$ with $2\le |\sub|\le e^{\frac{\cd}{C\beta\Delta}}$ there is a $1$-Lipschitz function $\f_{\sub}:\MM\to  L_2$ satisfying  
\begin{equation}\label{eq:lower bound on suit in lemma}
\forall (x,y)\in \mathscr{X}\times \YY,\qquad \big(\Delta\le d_\MM(x,y)\le \beta\Delta\big)\wedge \big(d_\MM(x,\sub)\le (\beta+1)\Delta\big)\implies \|\f_\sub(x)-\f_\sub(y)\|_{\! L_2}\gtrsim \frac{\Delta}{D}. 
\end{equation}
\end{lemma}

Prior to proving Lemma~\ref{lemm:get suit of functions}, which  is of  independent interest and could be of value for other purposes elsewhere (in particular, another application of it will be worked out in Section~\ref{sec:vanilla localized}), we will next demonstrate how to quickly deduce the second part of Theorem~\ref{thm:single scale on controlled growth} assuming its validity: 

\begin{proof}[Proof of the second part of Theorem~\ref{thm:single scale on controlled growth} assuming Lemma~\ref{lemm:get suit of functions}] We may suppose from now  that $|\cG_{\le K}(r,R)|\ge 2$ since if $\cG_{\le K}(r,R)=\emptyset$, then  the desired conclusion~\eqref{eq:x in local growth y general} is vacuous, and if $|\cG_{\le K}(r,R)|=1$, then let $x_0\in \MM$ be such that $\cG_{\le K}(r,R)=\{x_0\}$ and the desired   conclusion~\eqref{eq:x in local growth y general} holds for $\phi(x)=d_\MM(x,x_0)\in \R$. Let $C\ge 1$ be the constant from Lemma~\ref{lemm:get suit of functions}.  Assumption~\eqref{eq:assume single scale in theorem} of Theorem~\ref{thm:single scale on controlled growth} coincides with assumption~\eqref{eq:repeat single scale assumption for U} of Lemma~\ref{lem:rR}, and  the second part of assumption~\eqref{eq:Rr restrictions} Theorem~\ref{thm:single scale on controlled growth}  implies that $R-r\ge \diam_\MM(\cU)$ for every nonempty subset $\cU$ of $\MM$ for which $\diam_\MM(\cU)\le  C\beta( \log |\cG_{\le K} (r,R)|)\Delta$. We may therefore apply Lemma~\ref{lem:rR} to deduce that for every such $\cU$ there is a $1$-Lipschitz function $f_\cU:\MM\to L_2$ for which assumption~\eqref{eq:cX D version} of Lemma~\ref{lemm:get suit of functions} holds for $\XX=\YY=\MM$ and $\cd=C\beta( \log |\cG_{\le K} (r,R)|)\Delta$, and with  $D$ replaced by a positive universal constant multiple of $D$. The desired conclusion~\eqref{eq:x in local growth y general} of Theorem~\ref{thm:single scale on controlled growth} is now a special case of the conclusion~\eqref{eq:lower bound on suit in lemma} of Lemma~\ref{lemm:get suit of functions}, applied to $\sub=\cG_{\le K} (r,R)$, which  is valid as  $2\le |\cG_{\le K} (r,R)|= e^{\cd/(C\beta\Delta)}$.  
\end{proof}

Thus, in order to complete the proof of Theorem~\ref{thm:single scale on controlled growth} it remains to prove Lemma~\ref{lemm:get suit of functions}, which we do next:

\begin{proof} [Proof of Lemma~\ref{lemm:get suit of functions}]
\begin{comment}If $|\sub|=2$, then it is simple to achieve the conclusion of Lemma~\ref{lemm:get suit of functions} for, say, $\eta=1/5$. Write $\sub=\{u,v\}$  and define $\f_\Delta^\sub(x)=(d_\MM(x,u)/\sqrt{2},d_\MM(x,v)/\sqrt{2})\in \R^2$ for  $x\in \MM$. So, $\f_\Delta^\sub$ is $1$-Lipschitz from $\MM$ to $\ell_{\!\!2 }^2$. To check~\eqref{eq:lower bound on suit in lemma}, take $x\in \MM$ with $d_\MM(x,\sub)\le \Delta/(5\log 2)\le \Delta/3$. Without loss of generality $d_\MM(x,u)\le \Delta/3$. If $y\in \MM$ satisfies $d_\MM(x,y)\ge \Delta$, then $d_\MM(y,u)-d_\MM(x,u)\ge d_\MM(x,y)-2d_\MM(x,u)\ge \Delta-2\Delta/3=\Delta/3$. Thus, the first coordinate of $\f_\Delta^\sub(y)-\f_\Delta^\sub(x)$ is at least $\Delta/(3\sqrt{2})$. We will therefore henceforth assume that $|\sub|\ge 3$. 
\end{comment}
For each  $\cS\subset \MM$ and $K>0$, let $\mathrm{Part}_\MM(\cS;K)$ denote the set of all the partitions $\Part$ of $\cS$ that  are $K$-bounded with respect to  $d_\MM$, i.e., $\diam_\MM(\Part(x))\le K$ for every $x\in \cS$. 

 By~\cite{bartal1996probabilistic}, there exist universal constants $0<\fp<1<\kappa$ such that for  any  $\sub\subset \MM$ with $|\sub|> 1$ and any $\Delta_0>0$ there exists a probability measure $\mu_{\sub,\Delta_0}$ on $\mathrm{Part}_\MM(\sub;\Delta_0)$ satisfying the following condition:

$$
\forall x\in \sub,\qquad \mu_{\sub,\Delta_0}\Big[\Part_0\in \mathrm{Part}_\MM(\sub;\Delta_0):\ B_\MM \Big(x,\frac{1}{\kappa\log |\sub|}\Delta_0\Big)\cap \sub\subset \Part_0(x)\Big] \ge \fp.
$$
We will proceed to prove that Lemma~\ref{lemm:get suit of functions} holds if we take the universal constant $C$  equal to $60\kappa$. 

Fix any subset $\sub$ of $\MM$ whose size satisfies the assumptions in the statement of Lemma~\ref{lemm:get suit of functions}, i.e., 
\begin{equation}\label{eq:sub is not too small}
2\le |\sub|\le e^{\frac{\cd}{C\beta\Delta}}=e^{\frac{\cd}{60\kappa\beta\Delta}}.
\end{equation}
It will be convenient to set the following notations: 
\begin{equation*}\label{eq:notations for mu Delta0, Omega}
\left\{\begin{array}{ll}\Delta_0\eqdef  16\kappa(\beta+1)\big(\log |\sub|\big)\Delta,\\ \Omega\eqdef \mathrm{Part}_\MM\big(\MM;8(\beta+1)(2\kappa\log |\sub|+1)\Delta\big),\\ \mu\eqdef \mu_{\sub,\Delta_0}.
\end{array}\right.
\end{equation*}
Thus, $\Omega$ is  the set of $\Delta'$-bounded partitions of the entire metric space  $\MM$, where we introduce the notation
\begin{equation}\label{eq:Delta prime}
\Delta'\eqdef 8(\beta+1)(2\kappa\log |\sub|+1)\Delta\stackrel{\eqref{eq:sub is not too small}}{\le} 8(\beta+1)\Big(\frac{1}{30\beta \Delta}\cd+1\Big)\Delta\le \cd, 
\end{equation}
and the last step of~\eqref{eq:Delta prime} is valid because~\eqref{eq:sub is not too small} implies in particular that $\cd/\Delta\ge 60\kappa\beta \log 2$, and using this lower bound on $\cd/\Delta$ the rightmost inequality in~\eqref{eq:sub is not too small} is elementary to verify, since $\beta,\kappa\ge 1$. Now, a substitution of $\mu$ into~\cite[Lemma~3.8]{LN05} shows  that  there exists a probability measure $\nu$ on $\Omega$ such that
   \begin{equation}\label{eqn:pad'}
 \forall x\in \MM,\qquad d_\MM(x,\sub)\le (\beta+1) \Delta\implies    \nu \big[\Part\in \Omega:
    \ B_\MM(x,(\beta+1)\Delta)\subset \Part(x)\big] \ge  \fp.
    \end{equation}

By~\cite[Lemma~5.2]{mendel2004euclidean} there exists a function $G=G_\Delta :  L_2 \to  L_2$ such that:
    \begin{equation}\label{eqn:G delta}
\forall x, y \in  L_2,\qquad  \left\{\begin{array}{ll}      \norm{G (x)}_{\! L_2} = \norm{G (y)}_{\! L_2} = \Delta, \\ \frac{1}{2} \min\big\{ \Delta, \norm{x-y}_{\! L_2} \big\} \le \norm{G (x) - G (y)}_{\! L_2} \le \min\big\{ \Delta, \norm{x-y}_{\! L_2} \big\}.\end{array}\right.
    \end{equation}
For every $\mathscr{U}\subsetneq \MM$ define $\mathfrak{d}_\cU=\mathfrak{d}_{\mathscr{U},\Delta,d_\MM}:\MM \to [0,1]$ by
\begin{equation}\label{eq:def rho}
\forall z\in \MM,\qquad \mathfrak{d}_{\mathscr{U}}(z)\eqdef \min\Big\{1,\frac{1}{\Delta}d_\MM(z,\MM\setminus \mathscr{U})\Big\}.
\end{equation}
As the function $(t\in \R)\mapsto \min\{1,t/\Delta\}$ is $(1/\Delta)$-Lipschitz and nondecreasing, $\|\mathfrak{d}_{\mathscr{U}}\|_{\Lip(\MM;\R)}\le 1/\Delta$.   Finally, because by~\eqref{eq:Delta prime} we have $\diam_\MM(\Part(z))\le \cd$  for every $\Part\in \Omega$ and $z\in \MM$, so we can invoke the assumption of Lemma~\ref{lemm:get suit of functions} to fix a $1$-Lipschitz function $f_{\Part(z)}:\Part(z)\to  L_2$ that satisfies~\eqref{eq:cX D version} with $\mathscr{U}=\Part(z)$, i.e., 
    \begin{equation}\label{eqn:weakly bilip in local lemma}
 \forall (x,y)\in (\Part(z)\cap \XX)\times (\Part(z)\cap \YY),\qquad   \Delta\le  d_\MM (x, y)  \le \beta\Delta \implies \norm{f_{\Part(z)}(x) - f_{\Part(z)}(y)}_{\! L_2} \ge \frac{\Delta}{D}.
    \end{equation} 
In particular, $f_{\Part(z)}(z)$ is well-defined since $z\in \Part(z)$, so we can define $$\f_\sub: \MM\to L_2(\nu; L_2)\cong  L_2$$ by setting
    \begin{equation}\label{eq:our phi}
 \forall z\in \MM,\ \forall \Part\in \Omega,\qquad        \f_{\sub} (z)(\Part) \eqdef  \frac{\mathfrak{d}_{\Part(z)}(z) }{2} G \big(f_{\Part(z)} (z)\big).
    \end{equation}

    We first claim that the following point-wise bound holds: 
 \begin{equation}\label{eq:poin wise lip}
     \forall x,y\in \MM,\ \forall \Part\in \Omega,\qquad        \|\f_{\sub}(x)(\Part) -\f_{\sub}(y)(\Part)\|_{\! L_2} \le d_\MM(x,y).
    \end{equation}
    Observe that  after~\eqref{eq:poin wise lip} will be established, by squaring both of its sides  and then integrating the resulting estimate  $\ud \nu(\Part)$, we will deduce that $\f_{\sub}$ is $1$-Lipschitz as a function $\MM$ to $L_2(\nu; L_2)\cong L_2$.

  To verify~\eqref{eq:poin wise lip},   fix $x,y\in \MM$ and $\Part\in \Omega$. Suppose first that $\Part(x)=\Part(y)\eqdef \mathscr{U}$.  Then 
 \begin{align}
 \begin{split}
|\f_{\sub} (x)(\Part) -\f_{\sub} (x)(\Part)|&=\frac1{2} \Big\|\big(\mathfrak{d}_{\mathscr{U}}(x)-\mathfrak{d}_{\mathscr{U}}(y)\big) G \big(f_{\mathscr{U}} (x)\big) +\mathfrak{d}_{\mathscr{U}}(y)\Big(G \big(f_{\mathscr{U}} (x)\big)-G \big(f_{\mathscr{U}} (y)\big)\Big)\Big\|_{\! L_2} \\  &\le 
 \frac{|\mathfrak{d}_{\mathscr{U}}(x)-\mathfrak{d}_{\mathscr{U}}(y)|}{2} \big\|G \big(f_{\mathscr{U}} (x)\big)\big\|_{\! L_2}+\frac{\mathfrak{d}_{\mathscr{U}}(y)}{2} \big\| G \big(f_{\mathscr{U}} (x)\big)-G \big(f_{\mathscr{U}} (y) \big)\big\|_{\! L_2}\\&\le \frac{\|\mathfrak{d}_{\mathscr{U}}\|_{\Lip(\MM;\R)}}{2}\Delta d_\MM(x,y)+\frac1{2} \|f_{\mathscr{U}}\|_{\Lip(\MM; L_2)}d_\MM(x,y)\le d_\MM(x,y)\label{eq:use coordinate lip},
 \end{split}
 \end{align}
  where the first step of~\eqref{eq:use coordinate lip} is  a  consequence of the definition~\eqref{eq:our phi} of $\f_{\sub}$, using the current assumption $\Part(x)=\Part(y)=\mathscr{U}$, the penultimate step of~\eqref{eq:use coordinate lip} uses~\eqref{eqn:G delta} and the fact that $0\le \mathfrak{d}_{\mathscr{U}}(\cdot)\le 1$ by~\eqref{eq:def rho}, and the  final step of~\eqref{eq:use coordinate lip} holds as $f_{\mathscr{U}}$ is $1$-Lipschitz and $\mathfrak{d}_{\mathscr{U}}$ is $(1/\Delta)$-Lipschitz. This establishes~\eqref{eq:poin wise lip} if $\Part(x)=\Part(y)$.

It remains to verify~\eqref{eq:poin wise lip}  when $\Part(x)\neq\Part(y)$, i.e., $x\in \MM \setminus\Part(y)$ and $y\in \MM \setminus\Part(x)$, which implies that
 \begin{equation}\label{eq:different clusters far}
 \max \Big\{d_\MM\big(x,\MM\setminus \Part(x)\big),d_\MM\big(y,\MM\setminus \Part(y)\big)\Big\}\le d_\MM(x,y).
 \end{equation}
Consequently, when $\Part(x)\neq\Part(y)$ we can justify~\eqref{eq:poin wise lip} as follows:
 \begin{align*}
|\f_{\sub} (x)(\Part) -\f_{\sub}(x)(\Part)|&\le |\f_{\sub} (x)(\Part)|+|\f_{\sub} (x)(\Part)|\\&\stackrel{\eqref{eq:our phi}\wedge\eqref{eqn:G delta}}{=}  \frac12\mathfrak{d}_{\Part(x)}(x) \Delta+ \frac12\mathfrak{d}_{\Part(y)}(y) \Delta\\&\stackrel{\eqref{eq:def rho}}{\le}  \frac12 d_\MM\big(x,\MM\setminus \Part(x)\big)+ \frac12d_\MM\big(y,\MM\setminus \Part(y)\big)\\& \stackrel{\eqref{eq:different clusters far}}{\le}d_\MM(x,y).
 \end{align*}
    
   Next, fix $x,y\in \MM$ that satisfy $d_\MM(x,y)\le \beta\Delta$. Using the triangle inequality for $d_\MM$, it follows that  $B_\MM(y,\Delta)\subset B_\MM(x,(\beta+1)\Delta)$. Hence, if $\Part\in \Omega$ is such that $B_\MM(x,(\beta+1)\Delta)\subset \Part(x)$, then $\Part(y)=\Part(x)$ and also $\min\{d_\MM(x,\MM\setminus\Part(x)),d_\MM(y,\MM\setminus\Part(y))\}\ge \Delta$. Recalling~\eqref{eq:def rho}, it follows that 
\begin{align}\label{eq:rho is one when padded}
\begin{split}
 \big(d_\MM(x,y)\le \beta\Delta\big)\ &\wedge\ \big(B_\MM(x,(\beta+1)\Delta)\subset \Part(x)\big)\\&\implies \big(\Part(x)=\Part(y)\big)\wedge \big(\mathfrak{d}_{\Part(x)}(x)=\mathfrak{d}_{\Part(y)}(y)=1\big),
 \end{split} 
\end{align}
for every $x,y\in \MM$  and every $\Part\in\Omega_\Delta^\MM$. Consequently, the following lower bound holds for  all $(x,y)\in \XX\times \YY$ that satisfy $d_\MM(x,\sub)\le (\beta+1)\Delta$ and  $\Delta\le d_\MM(x,y)\le \beta\Delta$: 

\begin{align}\label{eq:for lower suit}
\begin{split}
 \|\f_{\sub}(x)&-\f_{\sub}(y)\|_{\!L_2\left(\nu; L_2\right)}\\&\ge \bigg(\int_{\left\{\Part\in \Omega:\ B_\MM(x,(\beta+1)\Delta)\subset \Part(x)\right\}}\|\f_{\sub}(x)-\f_{\sub}(y)\|_{\! L_2}^2\ud \nu(\Part)\bigg)^{\frac12}\\&=
\frac{1}{2}\bigg(\int_{\left\{\Part\in \Omega:\ B_\MM(x,(\beta+1)\Delta)\subset \Part(x)\right\}}\big\|G\big(f_{\Part(x)} (x)\big)-G \big(f_{\Part(x)} (y)\big)\big\|_{\! L_2}^2\ud \nu(\Part)\bigg)^{\frac12}\\ 
&\ge \frac{1}{4}\bigg(\int_{\left\{\Part\in \Omega:\ B_\MM(x,(\beta+1)\Delta)\subset \Part(x)\right\}}\min\big\{\Delta^2,\| f_{\Part(x)} (x)-f_{\Part(x)} (y)\|_{\! L_2}^2\big\}\ud \nu(\Part)\bigg)^{\frac12}\\
&\ge \frac{\Delta}{4D}\sqrt{\nu\big(\big\{\Part\in \Omega:\ B_\MM(x,(\beta+1)\Delta)\subset \Part(x)\big\}\big)} \\&\ge\frac{\Delta\sqrt{\fp}}{4D}\asymp\frac{\Delta}{D},
\end{split}
\end{align}
where the second step of~\eqref{eq:for lower suit} holds by~\eqref{eq:our phi}  using~\eqref{eq:rho is one when padded}, which is valid as $d_\MM(x,y)\le \beta\Delta$,  the third step of~\eqref{eq:for lower suit} uses the first inequality in~\eqref{eqn:G delta}, the fourth step of~\eqref{eq:for lower suit} uses~\eqref{eqn:weakly bilip in local lemma}, which is valid as $\Delta\le d_\MM(x,y)\le \beta \Delta$ and $(x,y)\in (\Part(x)\cap \XX)\times (\Part(x)\cap \YY)$ by~\eqref{eq:rho is one when padded}, and $D\ge 1$, and  the penultimate step of~\eqref{eq:for lower suit} is where~\eqref{eqn:pad'} is used, which is valid as $d_\MM(x,\sub)\le (\beta+1)\Delta$. We have thus proved the remaining part~\eqref{eq:lower bound on suit in lemma} of Lemma~\ref{lemm:get suit of functions}.\end{proof}

\begin{remark}\label{rem:padding} {\em An inspection of the proof of Lemma~\ref{lemm:get suit of functions}  reveals that the restriction on the size of $\sub$ appears in its statement only because by~\cite{bartal1996probabilistic} the {\em padding modulus} of $\sub$ is $O(\log |\sub |)$. Specifically, using the notation of~\cite{naor2024extension}, given $0<\fp<1$ let $\PAD_{\fp}(\sub)=\PAD_{\fp}(\sub,d_\MM)$ be the smallest $K\ge 1$ such that for every $\Delta>0$ there is a distribution over $\Delta$-bounded random partitions $\Part$ of $\sub$ with the property that for all $x\in \sub$ the probability that $B_\MM(x,\Delta/K)\cap \sub$ is contained in $\Part(x)$ is at least $\fp$. A repetition of the reasoning of Lemma~\ref{lemm:get suit of functions}  gives mutatis mutandis that if one replaces its requirement $|\sub|\le e^{\cd/(C\beta\Delta)}$ by $\PAD_\fp(\sub)\le \cd/ (C\beta\Delta)$, then its conclusion~\eqref{eq:lower bound on suit in lemma} holds with $\|\f_{\sub}(x)-\f_{\sub}(y)\|_{\! L_2}\gtrsim \Delta/D$ replaced by $$\|\f_{\sub}(x)-\f_{\sub}(y)\|_{\! L_2}\gtrsim \frac{\Delta\sqrt{\fp}}{D}.$$}
\end{remark}

\section{Localized weakly bi-Lipschitz embeddings}\label{sec:vanilla localized}

The following ``localized version'' of Definition~\ref{def:npss} is a slight generalization of a  definition that appeared in~\cite[Section 7.2]{naor2014comparison}, which corresponds to the (arbitrary) choice $K=32$ below:

\begin{defn}\label{def:weak bilip} Given $K,D \ge 1$, a metric space $(\MM,d_\MM)$ is said to admit a $K$-localized weakly bi-Lipschitz embedding into a metric space $(\NN,d_\NN)$ with distortion $D$ if for every $\Delta > 0$ and every $z\in \MM$ there exists a non-constant Lipschitz function $f^z_\Delta : \MM \to \NN$ such that 
    \begin{equation}\label{eqn:local single-scale}
 \forall x,y\in B_\MM(z,K\Delta),\qquad        d_\MM (x, y) \ge \Delta \implies d_\NN \big(f^z_\Delta (x), f^z_\Delta (y)\big) \ge \frac{\|f^z_\Delta\|_{\Lip}}{D} \Delta,
    \end{equation}
\end{defn}

A key (well known) property of $ L_p$ that we will use herein is the (first part of) the following theorem: 
\begin{theorem}\label{cor:mazur lp}
If $K> 1$ and $p>2$, then $ L_p$ admits a $K$-localized weakly bi-Lipschitz embedding into $ L_2$  with distortion $D$, where $D\lesssim p2^{p/2}K^{p/2-1}$. Conversely, if $ L_p$ admits a $K$-localized weakly bi-Lipschitz embedding into $ L_2$  with distortion $D$, then necessarily $D\gtrsim 2^{p/2}K^{p/2-1}$. 
\end{theorem}

\begin{proof} The first part of Theorem~\ref{cor:mazur lp} essentially coincides with Lemma~7.6 of~\cite{naor2014comparison}, which notes the (special case $K=32$ of the) following general statement. Suppose that $(\bX, \|\cdot\|_\bX)$ and $(\bY,\|\cdot\|_\bY)$ are Banach spaces for which there exists a function $U:B_\bX\to \bY$ from the unit ball $B_\bX=\{x\in \bX:\ \|x\|_\bX\le 1\}$ to $\bY$ satisfying 
\begin{equation}\label{eq:uniform continuous condition}
\forall x,y\in B_\bX,\qquad \omega\big(\|x-y\|_\bX\big)\le \|U(x)-U(y)\|_{\bY}\le L\|x-y\|_\bX,
\end{equation}
for some $L>0$ and an increasing modulus $\omega:[0,\infty)\to [0,\infty)$. Let $\rho_\bX:\bX\to B_\bX$ be the standard retraction (e.g.~\cite[equation~(5.2)]{naor2014comparison}) from $\bX$ onto $B_\bX$, i.e.,  $\rho_\bX(x)=x/\max\{\|x\|_\bX,1\}$ for every $x\in \bX$. It is straightforward to check that  $\|\rho_\bX\|_{\Lip}\le 2$, so if we define for $z\in \bX$ and $K,\Delta>0$ a function $f_\Delta^z:\bX\to \bY$ by:
$$
\forall x\in \bX,\qquad f_\Delta^z(x)\eqdef K\Delta U \bigg(\rho_\bX \Big(\frac{1}{K\Delta}(x-z)\Big)\bigg),
$$
then $\|f_\Delta^z\|_{\Lip}\le 2L$ by the second inequality in~\eqref{eq:uniform continuous condition}. Using this together with the first inequality in~\eqref{eq:uniform continuous condition} shows that for every $x,y\in B_\bX(z,K\Delta)$, if $\|x-y\|_\bX\ge \Delta$, then  $\|f_\Delta^z(x)-f_\Delta^z(y)\|_\bY\ge K\omega(1/K)\|f_\Delta^z\|_{\Lip}\Delta/(2L)$. Thus, $\bX$ admits a   $K$-localized weakly bi-Lipschitz embedding into $\bY$  with distortion $2L/(K\omega(1/K))$. 

In particular, if in~\eqref{eq:uniform continuous condition} we have $\omega(t)= t^\alpha/\beta$ for all $t\ge 0$ and some $\alpha,\beta\ge 1$, then $\bX$ admits a   $K$-localized weakly bi-Lipschitz embedding into $\bY$  with distortion $2L\beta K^{\alpha-1}$. When $\bX= L_p$ for some $p\ge 2$ and $\bY= L_2$, one can take $U$ to be the restriction to the unit ball of $ L_p$ of the classical Mazur map~\cite{mazur1929remarque} $M_{p\to 2}: L_p\to L_2$, in which case this  holds with $\alpha=p/2$ and $\beta\asymp 2^{p/2}$, and with $L$ in~\eqref{eq:uniform continuous condition} satisfying $L\asymp p$. A derivation of these values of $\alpha,\beta,L$ for the Mazur map appears in~\cite{naor2014comparison}; see specifically equation~(5.32) there. 

We will establish the reverse direction by adjusting the proof of~\cite[Lemma~52]{naor2014comparison}.  So, suppose that $ L_p$ admits a $K$-localized weakly bi-Lipschitz embedding into $ L_2$  with distortion $D$, and our goal is to bound $D$ from below. As the complex plane $\C= \ell_{\!\!2 }^2$, embeds into $ L_p$ with distortion $1+\e$ for any $0<\e<K-1$ (e.g.~by Dvoretzky's theorem~\cite{Dvo60}, though using that theorem is overkill for this purpose), it follows that also $ L_p(\C)$ admits a $K/(1+\e)$-localized weakly bi-Lipschitz embedding into $ L_2$  with distortion $(1+\e)D$. Therefore,  it suffices to prove that if $ L_p(\C)$ admits a $K$-localized weakly bi-Lipschitz embedding into $ L_2$  with distortion $D$, then necessarily $D\gtrsim 2^{p/2}K^{p/2-1}$  (considering $\e\asymp 1/p$ above suffices).

Fix $n\in \N$.  Applying the assumption that $ L_p(\C)$ admits a $K$-localized weakly bi-Lipschitz embedding into $ L_2$  of distortion $D$ to $\Delta=n^{1/p}/K$ in Definition~\ref{def:weak bilip} , it follows that  there is $f:n^{1/p}B_{\ell_{\!\! p }^n(\C)}\to  L_2$ satisfying:
\begin{align}\label{eq:f assumption over C}
\begin{split}
\left\{\begin{array}{ll} \|f\|_{\Lip}=1,\\
\forall x,y\in n^{\frac{1}{p}} B_{\ell_{\!\! p }^n(\C)}, \qquad \|x-y\|_{\!\ell_{\!\! p }^n(\C)}\ge \frac{n^{\frac{1}{p}}}{K}\implies \|f(x)-f(y)\|_{\! L_2}\ge \frac{n^{\frac{1}{p}}}{KD}. 
\end{array}\right.
\end{split}
\end{align}

For $m\in \N$  write $\mathfrak{u}_m\eqdef e^{2\pi\i/m}\in \C$ and let $\U_m\eqdef \{1,\mathfrak{u}_m,\mathfrak{u}_m^2,\ldots, \mathfrak{u}_m^{m-1}\}\subset \C$ be the cyclic group of the  roots of unity of order $m$. By (the Hilbertian special case of) Theorem~5.2 in~\cite{NS16} (whose conclusion is a slightly simpler variant of~\cite[Theorem~4.1]{MN08}\footnote{Alternatively one could combine here Theorem~4.1 in~\cite{MN08} with the discussion in Section~4 of~\cite{EMN19}.}),  if $m\ge \sqrt{n}$ and $m$ is divisible by $8$, then every $g:\U_m^n\to  L_2$ satisfies: 
 \begin{align}\label{eq:quote cotype}
 \begin{split}
\sum_{j=1}^n \sum_{x\in \U_m^n} &\big\|g(x)-g(x_1,\ldots,x_{j-1},-x_j,x_{j+1},\ldots,x_n)\big\|_{\! L_2}^2\\ &\lesssim \frac{m^2}{2^n} \sum_{\e\in \{-1,1\}^n}\sum_{x\in \U_m^n} \big\|g(x)-g(\mathfrak{u}_m^{\e_1}x_1,\ldots,\mathfrak{u}_m^{\e_n}x_n)\big\|_{\! L_2}^2.
\end{split}
\end{align}

We will apply~\eqref{eq:quote cotype} to the restriction to $\U_m^n\subset n^{1/p}B_{\ell_{\!\! p }^n(\C)}$ of the function $f:n^{1/p}B_{\ell_{\!\! p }^n(\C)}\to  L_2$ in~\eqref{eq:f assumption over C}. For that, observe that  since $f$ is $1$-Lipschitz, every $\e=(\e_1,\ldots,\e_n)\in \{-1,1\}^n$ and  $x=(x_1,\ldots,x_n)\in \U_m^n$ satisfy
$$ \big\|f(x)-f(\mathfrak{u}_m^{\e_1}x_1,\ldots,\mathfrak{u}_m^{\e_n}x_n)\big\|_{\! L_2}\le \big\|(1-\mathfrak{u}_m^{\e_1})x_1,\ldots,(1-\mathfrak{u}_m^{\e_n})x_n\big\|_{\!\ell_{\!\!p }^n(\C)}= \sqrt{2} n^{\frac{1}{p}} \Big|1-\cos\Big(\frac{2\pi}{m}\Big)\Big|\asymp \frac{n^{\frac{1}{p}}}{m}. 
$$
Hence, the right hand side of~\eqref{eq:quote cotype} for $g=f|_{\U_m^n}$ can be bounded from above as follows: 
\begin{equation}\label{eq:upper for cotype}
\frac{m^2}{2^n} \sum_{\e\in \{-1,1\}^n}\sum_{x\in \U_m^n} \big\|f(x)-f(\mathfrak{u}_m^{\e_1}x_1,\ldots,\mathfrak{u}_m^{\e_n}x_n)\big\|_{\! L_2}^2\lesssim n^{\frac{2}{p}}m^{n}. 
\end{equation}
To bound the right hand side of~\eqref{eq:quote cotype} from below, suppose that $n\le (2K)^p$. Then, for every $x\in \U_m$ and every $j\in \{1,\ldots, n\}$ we have $\|x- (x_1,\ldots,x_{j-1},-x_j,x_{j+1},\ldots,x_n)\|_{ L_p(\C)}=2\ge n^{1/p}/K$, so using~\eqref{eq:f assumption over C} we get: 
$$
\forall x\in \U_m^m,\qquad \big\|f(x)-f(x_1,\ldots,x_{j-1},-x_j,x_{j+1},\ldots,x_n)\big\|_{\! L_2}\ge \frac{n^{\frac{1}{p}}}{KD}. 
$$
Consequently, 
\begin{equation}\label{eq:lower for cotype}
\sum_{j=1}^n \sum_{x\in \U_m^n} \big\|f(x)-f(x_1,\ldots,x_{j-1},-x_j,x_{j+1},\ldots,x_n)\big\|_{\! L_2}^2\ge \frac{n^{1+\frac{2}{p}}m^n}{(KD)^2}.  
\end{equation}

By combining~\eqref{eq:upper for cotype} and~\eqref{eq:lower for cotype} with~\eqref{eq:quote cotype} and rearranging, we deduce that   $D\gtrsim \sqrt{n}/K$. The  requirements for this to hold were that $n\le (2K)^p$, as well as that $m\ge \sqrt{n}$ and $m\equiv 0\mod 8$, so choosing  $n= \lfloor (2K)^p\rfloor$ and $m=8\left\lceil \sqrt{n}\right\rceil$ gives $D\gtrsim \sqrt{(2K)^p-1}/K\asymp 2^{p/2}K^{p/2-1}$, as required.
\end{proof} 

\begin{remark}{\em 
It is worthwhile to note (though not needed for the results herein) that by incorporating the Maurey--Pisier theorem~\cite{MP76} and the main result of~\cite{MN08} that Rademacher cotype and metric cotype coincide, the above proof of the lower bound on $D$ in Theorem~\ref{cor:mazur lp} yields mutatis mutandis the following statement. Let $(\bX,\|\cdot\|_\bX), (\bY,\|\cdot\|_\bY)$ be Banach spaces. Suppose $\bY$ has Rademacher cotype $q$ and that $p\ge 2$ equals the infimum over those $p'\ge 2$ such that $\bX$ has Rademcaher cotype $p'$. Fix $K,D>1$ and $p>q$. If $\bX$  admits a $K$-localized weakly bi-Lipschitz embedding into $\bY$  with distortion $D$, then $D\gtrsim 2^{p/q}K^{p/q-1}/C_q(\bY)$, where $C_q(\bY)$ is the Rademacher cotype $q$ constant of $\bY$.}
\end{remark}

As done in~\eqref{eq:2Delta} for the  two-sided scale-localized variant of Definition~\ref{def:npss}, we will also say in the context of Definition~\ref{def:weak bilip} that a metric space $(\MM,d_\MM)$ admits a {\em two-sided scale-localized $K$-localized weakly bi-Lipschitz embedding into a metric space $(\NN,d_\NN)$ with distortion $D>0$} if for every $\Delta > 0$ and every $z\in \MM$ there exists a non-constant Lipschitz function $g^z_\Delta : \MM \to \NN$ such that 
    \begin{equation}\label{eqn:glocal single-scale}
 \forall x,y\in B_\MM(z,K\Delta),\qquad        \Delta\le d_\MM (x, y) \le 2\Delta \implies d_\NN \big(g^z_\Delta (x),g^z_\Delta (y)\big) \ge \frac{\|g^z_\Delta\|_{\Lip}}{D} \Delta.
    \end{equation}
    
 The localization principle of Lemma~\ref{lemm:get suit of functions} allows us to relate as follows localized weakly bi-Lipschitz embeddings to the usual (global) weakly bi-Lipschitz embeddings that were discussed in Section~\ref{sec:single scale}: 

\begin{lemma}\label{lem:from localized to every scale} There exists a universal constant $\kappa\ge 2$ with the following property. Fix a nondecreasing function $D:(1,\infty)\to (1,\infty)$. Suppose that $(\MM,d_\MM)$ is a metric space admitting a two-sided scale-localized  $K$-localized weakly bi-Lipschitz embedding into $ L_2$ with distortion less than $D(K)$ for every $K>1$. Then, 
\begin{equation}\label{eq:hat bound from localized}
\forall n\in \{2,3,\ldots\},\qquad \widehat{\dd}_2^n(\MM)\le D(\kappa\log n). 
\end{equation}
\end{lemma}

\begin{proof} We will show that~\eqref{eq:hat bound from localized} holds for $\kappa=2C$, where $C\ge 1$ is the universal constant of Lemma~\ref{lemm:get suit of functions}. For every nonempty bounded subset $\cU$ of $\MM$, fix any $z_\cU\in \cU$. Then, $\cU\subset B_\MM(z_\cU,\diam_\MM(\cU))$, whence for every $\Delta>0$, by considering the restriction to $\cU$  of the function in~\eqref{eqn:glocal single-scale} with $K=\diam_\MM(\cU)/\Delta$ and $z=z_\cU$, the assumption of Lemma~\ref{lem:from localized to every scale} implies that there exists a $1$-Lipschitz function $f_{\cU,\Delta}:\MM\to L_2$ that satisfies
$$
\forall x,y\in \cU,\qquad \Delta\le d_\MM(x,y)\le 2\Delta\implies \|f_{\cU,\Delta}(x)-f_{\cU,\Delta}(y)\|_{\! L_2}\ge \frac{\Delta}{D\left(\frac{\diam_\MM(\cU)}{\Delta}\right)}. 
$$

So, for every integer $n\ge 2$ we may apply Lemma~\ref{lemm:get suit of functions} with  $\XX=\YY=\MM$ and  $\beta=2$, as well as $\cd=2C(\log n)\Delta$ and $D=D(\cd/\Delta)$, to obtain for every $\sub\subset \MM$ with $|\sub|=n$ a $1$-Lipschitz function $\f_{\sub,\Delta}:\MM\to  L_2$ that satisfies the following for every $x,y\in \MM$:
$$
\big(\Delta\le d_\MM(x,y)\le 2\Delta\big)\wedge \big(d_\MM(x,\sub)\le 3\Delta\big)\implies \|\f_{\sub,\Delta}(x)-\f_{\sub,\Delta}(y)\|_{\! L_2}\gtrsim \frac{\Delta}{D(2C\log n)}. 
$$
This is more than what is needed to deduce that $\widehat{\dd}_2^n(\MM)\le D(2C\log n)$, per~\eqref{eq:2Delta}. In particular, we obtained a function $\f_{\sub,\Delta}$ that is defined on all of $\MM$ while we only need it to be a $1$-Lipschitz function on $\sub$. For that, in the above application  Lemma~\ref{lemm:get suit of functions} it would suffice that $f_{\cU,\Delta}$ is defined on $\cU$, for which we could work here with the weakening of the notion of two-sided scale-localized  $K$-localized weakly bi-Lipschitz embedding in which the function $g^z_\Delta$ in~\eqref{eqn:glocal single-scale} is defined only on $B_\MM(z,K\Delta)$. This observation could be relevant to  future investigations, but not when $\MM$ is a Banach space since in that case Lipschitz functions can be extended from any ball to the super-space while increasing their Lipschitz constant by a factor of at most  $2$, which is seen by composing with the radial retraction as we did in the proof of Theorem~\ref{cor:mazur lp}.  
\end{proof}

With Lemma~\ref{lem:from localized to every scale} at hand, we can now summarize the best available upper bounds  on $\widehat{\dd}_2^n( L_p)$:
\begin{equation}\label{eq:current single scale bounds}
\forall n\in \{2,3\ldots,\},\qquad \widehat{\dd}_2^n( L_p)\lesssim\left\{\begin{array}{ll}\sqrt{\log n}&\mathrm{if}\  1\le p \le \sqrt{5}-1,\\
(\log n)^{\frac{2}{p}\left(\frac{1}{p}-\frac12\right)}&\mathrm{if}\ \sqrt{5}-1\le p\le 2,\\
(\log n)^{\frac{p}{2}-1}&\mathrm{if}\ 2\le p\le 3,\\
p^3\sqrt{\log n} &\mathrm{if}\ 3\le p\le \sqrt[6]{\log n},\\
\log n &\mathrm{if}\ p\ge \sqrt[6]{\log n}.\end{array}\right.
\end{equation}
Our contribution here is the range $2\le p=o(\sqrt[6]{\log n})$ of~\eqref{eq:current single scale bounds}, in which the previously best known estimate was nothing more than the $O(\log n)$ upper bound that holds by~\cite{bourgain1985lipschitz} for any $n$-point metric space. The case $2\le p\le 3$ of~\eqref{eq:current single scale bounds} follows by  substituting Theorem~\ref{cor:mazur lp} into Lemma~\ref{lem:from localized to every scale}, and if $3\le p\le \sqrt[6]{\log n}$, then~\eqref{eq:current single scale bounds} follows by  substituting Theorem~\ref{cor:lp} and the bound $\ee( L_p; L_2)\lesssim \sqrt{p}$ of~\cite{NPSS06} into Lemma~\ref{lem:sep ext product} (using Theorem~\ref{cor:mazur lp}  and Lemma~\ref{lem:from localized to every scale} here would yield a weaker result). If $1\le p \le \sqrt{5}-1$, then~\eqref{eq:current single scale bounds} is due to~\cite{ARV09}, if $\sqrt{5}-1\le p\le 2$, then~\eqref{eq:current single scale bounds} is due to~\cite{Lee05}, and if $p\ge\sqrt[6]{\log n}$, then~\eqref{eq:current single scale bounds} follows from~\cite{bourgain1985lipschitz}.

The best available lower bound on $\widehat{\dd}_2^n( L_p)$ is a universal constant multiple of $(\log n)^{1/p-1/2}$ if $1\le p\le 2$ and $\max\{\min\{p,\log n\},((\log n)/\log\log n)^{1/2-1/p}\}$ if $p\ge 2$; when $1\le p\le 2$ this follows from~\cite{Enf69} by considering as in Section~\ref{sec:lewis} the discrete $k$-dimensional hypercube with the $\ell_{\!\!p }^k$ metric,    and for $p\ge 2$ the $\Omega(\min\{p,\log n\})$ lower bound follows from~\cite{Mat97} while the $\Omega(((\log n)/\log\log n)^{1/2-1/p}$ lower bound follows from~\cite{MN08} (one cannot consider for this purpose  the planar graph example that was used in Section~\ref{sec:lewis}).  

Thus, we  know that~\eqref{eq:current single scale bounds} is optimal only when  $(p-1)\log n =O(1)$ and $p\gtrsim \log n$ (and, trivially, when $p=2$), but it plausibly not optimal for  the rest of the possible values of $p$. It would be worthwhile (and likely challenging) to obtain asymptotically sharp bounds here.

\section{Localized radially weakly bi-Lipschitz embeddings}\label{sec:all p>2}

As we will see later, Definition~\ref{def:weak bilip} (correspondingly, the first part of Theorem~\ref{cor:mazur lp}) suffices for  proving Theorem~\ref{thm:all >2}, Theorem~\ref{thm:new ext intro}, Theorem~\ref{thm:doubling ext}, Theorem~\ref{cor:lp}, and Theorem~\ref{thm:doubling sep} with a constant factor that has  a much worse (exponential) dependence on $p$.  The lower bound in the second part of Theorem~\ref{cor:mazur lp} shows that such a loss is inherent to this approach. It  is more delicate to get the stated dependence on $p$  in the above theorems. For that purpose, we will next introduce a (quite subtle, but crucial) variant of Definition~\ref{def:weak bilip} which is interesting in its own right and likely useful for other purposes beyond its applications that we derive  herein.  We will prove that the Mazur map obeys the aforementioned variant with a much better dependence on $p$, and demonstrate that this new embedding notion preserves the separation modulus.

To motivate Definition~\ref{def:radial} below, clarify its geometric meaning, and explain its nuanced difference from Definition~\ref{def:weak bilip}, we will start by examining the  following  consequence of Definition~\ref{def:weak bilip}.

Fix $K,\Delta>0$ and suppose that a metric space $(\MM,d_\MM)$  admits a $K$-localized weakly bi-Lipschitz embedding into a metric space $(\NN,d_\NN)$ with distortion $D$. Thus, for every $\Delta>0$ and every $z\in \MM$ there is a nonconstant Lipschitz function $f_\Delta^z:\MM\to \NN$ for which~\eqref{eqn:local single-scale} holds.  If $x\in \MM$ and $y,w\in B_\MM(z,K\Delta)$ satisfy $d_\MM(y,w)\ge \Delta$, then,
%\begin{align*}
% \frac{\|f^z_\Delta\|_{\Lip}}{D}\Delta \stackrel{\eqref{eqn:local single-scale} }{\le} d_\NN \big(f_\Delta^z(y),f_\Delta^z(w)\big) &\le d_\NN \big( %f^z_\Delta (y),f^z_\Delta (x)\big)+d_\NN \big(f^z_\Delta (x), f^z_\Delta (w)\big)\\&\le 2\max\Big\{ d_\NN \big(f^z_\Delta (x), f^z_\Delta (y)\big),d_\NN %\big(f^z_\Delta (x), f^z_\Delta (w)\big)\Big\}. 
%\end{align*}
\begin{align*}
 \frac{\|f^z_\Delta\|_{\Lip}}{D}\Delta \stackrel{\eqref{eqn:local single-scale} }{\le} d_\NN \big(f_\Delta^z(y),f_\Delta^z(w)\big) \le 2\max\Big\{ d_\NN \big(f^z_\Delta (x), f^z_\Delta (y)\big),d_\NN \big(f^z_\Delta (x), f^z_\Delta (w)\big)\Big\}. 
\end{align*}
Letting $B_\MM^\circ(x,r)=\{y\in \MM:\ d_\MM(x,y)<r\}$ be the open  $d_\MM$-ball centered at $x\in \MM$ of radius $r\ge 0$, we get:
$$
\forall x\in \MM, \  \forall y,w\in B_\MM(z,K\Delta), \qquad d_\MM(y,w)\ge \Delta\implies  \{y,w\}\not\subset \big( f_\Delta^z \big)^{-1} \Big(B_\NN^\circ\big(f^z_\Delta(x),\frac{\|f^z_\Delta\|_{\Lip}}{2D} \Delta\big)\Big).
$$
The contrapositive to this implication implies that
\begin{equation}\label{eq:deduce from weak bilip}
\forall x\in \MM,\qquad \diam_\MM \bigg(B_\MM(z,K\Delta)\cap\big( f_\Delta^z \big)^{-1} \Big(B_\NN^\circ\big(f^z_\Delta(x),\frac{\|f^z_\Delta\|_{\Lip}}{2D} \Delta\big)\Big)\bigg)\le\Delta.
\end{equation}

Definition~\ref{def:radial} below is important for our purposes. It modifies~\eqref{eq:deduce from weak bilip} in the following 5 ways. Firstly, it considers only  closed balls and it replaces the weak inequality in~\eqref{eq:deduce from weak bilip} by a strict inequality; both of these modifications are essentially cosmetic as they have insignificant impact on the results that are obtained herein and they are introduced for ease of remembering the definition and to streamline its subsequent implementations. Definition~\ref{def:radial} also changes the term $2D$ in~\eqref{eq:deduce from weak bilip} to $D$, which is merely a matter of notational convenience that impacts only the values of the (mostly implicit) universal constant factors in  the results that are obtained herein.  Next, it requires  the nonconstant Lipschitz function $f_\Delta^z$ to be defined on $B_\MM(z,K\Delta)$ rather than on all of $\MM$; while this could be  a genuine weakening when $\MM$ is an arbitrary  metric space, if $\MM$ is a Banach space, then (as we recalled in the proof of Theorem~\ref{cor:mazur lp})  the standard normalization mapping yields a $2$-Lipschitz retraction from $\MM$ onto $B_\MM(z,K\Delta)$, so by precomposing with this retraction we would obtain  the same property for a function that is now defined on all of $\MM$, at the cost of replacing  $D$ by $2D$; again, this would only impact the universal constant factors in  the results that are obtained herein. Finally, Definition~\ref{def:radial} replaces the $d_\MM$-diameter in~\eqref{eq:deduce from weak bilip} by the $d_\MM$-radius. While one might  think that due to~\eqref{eq:rad2diam} this is also a minor modification, it is, in fact, substantial, as it leads to an exponential improvement of the dependence on $p$, as expressed in Lemma~\ref{thm:radial} below. 

\begin{definition}\label{def:radial} Given $K,D >0$, we say that a metric space $(\MM,d_\MM)$ admits a $K$-localized {\bf radially} weakly bi-Lipschitz embedding into a metric space $(\NN,d_\NN)$ with distortion $D$ if for every $\Delta > 0$ and every $z\in \MM$ there exists a non-constant Lipschitz function $f^z_\Delta : B_\MM(z,K\Delta) \to \NN$ such that 
\begin{equation}\label{eq:radially def}
\forall x\in B_\MM(z,K\Delta),\qquad  \rad_\MM \bigg(\big( f_\Delta^z \big)^{-1} \Big(B_\NN\big(f^z_\Delta(x),\frac{\|f^z_\Delta\|_{\Lip}}{D} \Delta\big)\Big)\bigg)<\Delta.
\end{equation}
\end{definition}

By the above discussion, if a metric space $(\MM,d_\MM)$  admits a $K$-localized weakly bi-Lipschitz embedding into a metric space $(\NN,d_\NN)$ with distortion $D$, then for any $\e>0$  it also admits a $(1+\e)K$-localized {radially} weakly bi-Lipschitz embedding into $\NN$ with distortion $2(1+\e)^2D$.

Conversely, suppose that there is $L\ge 1$  such that for every ball $B_\MM(z,r)\subset \MM$ there exists an $L$-Lipschitz retraction $\rho^z_r$ from $\MM$ onto $B_\MM(z,r)$; if $\MM$ is a Banach space, then we can take $L=2$. If $f:B_\MM(z,K\Delta)\to 
\NN$ is a nonconstant Lipschitz function satisfying~\eqref{eq:radially def}, then the function $F_\Delta^z=f_\Delta^z\circ \rho_r^z:\MM\to \NN$ extends $f_\Delta^z$ and satisfies $\|F_\Delta^z\|_{\Lip}\le L\|f_\Delta^z\|_{\Lip}$. Therefore, for every $x\in B_\MM(z,K\Delta)$ we have:
\begin{align*}
 \diam_\MM \bigg(\big( & f_\Delta^z \big)^{-1} \Big(B_\NN\big(f^z_\Delta(x),\frac{\|F^z_\Delta\|_{\Lip}}{LD} \Delta\big)\Big)\bigg)\\&\le \diam_\MM \bigg(\big( f_\Delta^z \big)^{-1} \Big(B_\NN\big(f^z_\Delta(x),\frac{\|f^z_\Delta\|_{\Lip}}{D} \Delta\big)\Big)\bigg)\stackrel{\eqref{eq:radially def}\wedge \eqref{eq:rad2diam}}{<} 2\Delta.
\end{align*}
It follows that  $y\notin ( f_\Delta^z )^{-1} (B_\NN\big(f^z_\Delta(x),\|F^z_\Delta\|_{\Lip}\Delta/(LD))))$ whenever $x,y\in B_\MM(z,K\Delta)$ satisfy $d_\MM(x,y)\ge  2\Delta$. In particular, since $F_\Delta^z$ extends $f_\Delta^z$, we have 
\begin{align*}
\forall x,y\in B_\MM(z,K\Delta)&=B_\MM\Big(z,\frac{K}{2}(2\Delta)\Big),\\&\qquad d_\MM(x,y)\ge2 \Delta\implies d_\NN\big(F_\Delta^z(x),F_\Delta^z(y)\big)\ge \frac{\|F^z_\Delta\|_{\Lip}}{LD}\Delta=\frac{\|F^z_\Delta\|_{\Lip}}{2LD}2 \Delta. 
\end{align*}

Recalling Definition~\ref{def:weak bilip}, this observation shows that (under the assumption that an $L$-Lipschitz retraction from $\MM$ onto any ball in $\MM$ exists), if  $\MM$ admits a $K$-localized {radially} weakly bi-Lipschitz embedding into $\NN$ with distortion $D$, then  also $\MM$ admits a $(K/2)$-localized weakly bi-Lipschitz embedding into $\NN$ with distortion $2 LD$. In particular, if $\MM$ is a Banach space, then we deduce that it admits a $(K/2)$-localized weakly bi-Lipschitz embedding into $\NN$ with distortion $4 D$. 

The above loss of a $2 L$ factor in the distortion is of secondary importance if $L=O(1)$, as it only impacts the universal constant factors in  the results that are obtained herein. However, the above reduction of   $K$ to $K/2$  is significant, as demonstrated by Lemma~\ref{thm:radial} below, which should be contrasted with  the second part of Theorem~\ref{cor:mazur lp}, since it shows that if the word ``radially'' would be omitted from the statement of Lemma~\ref{thm:radial}, then the distortion in its conclusion would have to be $e^{\Omega(p)}$.

\begin{lemma}\label{thm:radial} If $p\ge 2$, then there are $K=K(p),D=D(p)>1$ satisfying $K-1\asymp 1/p$ and $D\asymp p$ such that  $ L_p$ admits a $K$-localized {\bf radially} weakly bi-Lipschitz embedding into $ L_2$  with distortion $D$.
\end{lemma}

Lemma~\ref{thm:radial} provides a novel property (potentially of use beyond its applications herein) of the Mazur map~\cite{mazur1929remarque}, whose classical definition was recalled in~\eqref{eq:def mazur map}; the relevant property, which implies Lemma~\ref{thm:radial}  and is, in fact, what Lemma~\ref{thm:radial}  will prove, was already stated in the Introduction as the inclusion~\eqref{eq:radial mazur normalized}. 

It is straightforward to check that $M_{p\to q}$ maps $ L_p$ bijectively onto $L_q$ with $M_{p\to q}^{-1}=M_{q\to p}$. Furthermore, $\|M_{p\to q}(x)\|_{L_q}=\|x\|_{ L_p}$  and $M_{p\to q}(sx)=s^{p/q}\sgn(s)M_{p\to q}(x)$ for every $x\in  L_p$ and $s\in \R$. The property of the Mazur map that lies at the heart of Lemma~\ref{thm:radial} is the following inequality:

\begin{lemma}\label{lem:radial mazur} If $p\ge 2$, then for every $0<\alpha\le \frac{1}{p}$ and every $0<\lambda<1$ we have\footnote{The constant $5$ in the right hand side of~\eqref{eq:new mazur inequality} is neither optimal nor does it have meaningful impact on the ensuing results. }  
\begin{equation}\label{eq:new mazur inequality}
\forall x,y\in  L_p,\qquad \|y-\alpha x\|_{\! L_p}^p\le (1-\lambda\alpha)^p\|x\|_{\! L_p}^2+\frac{5}{(1-\lambda)p\alpha} \big\|M_{p\to 2}(y)-M_{p\to 2}(x)\big\|_{\! L_2}^2.
\end{equation}
\end{lemma}
\noindent Note that since $\lambda$ appears only in the right hand side of~\eqref{eq:new mazur inequality}, the optimal way to use Lemma~\ref{lem:radial mazur} is to apply~\eqref{eq:new mazur inequality} with the value $\lambda=\lambda(x,y,p,\alpha)$  that minimizes the right hand side of~\eqref{eq:new mazur inequality} over $0\le \lambda\le 1$.

Assuming the validity of Lemma~\ref{lem:radial mazur} for the moment,  we will next show how to deduce Lemma~\ref{thm:radial}.

\begin{proof}[Proof of Lemma~\ref{thm:radial} assuming Lemma~\ref{lem:radial mazur}] We will start by demonstrating that there exists a universal constant $c>0$ such that the following inclusion holds for every $0<\alpha\le 1/p$ and $0<\sigma,\lambda<1$:\footnote{Stating~\eqref{eq:lambda sigma inclusion} for general parameters $\alpha,\lambda,\sigma$ is beneficial because its proof is simpler to read without making choices that are arbitrary at this juncture, and also this could be relevant for future investigations. Below, we will use~\eqref{eq:lambda sigma inclusion} for $\alpha=1/p$ and $\lambda=\sigma =1/2$, which suffices for the specific application herein even though these choices of $\lambda,\sigma$ are not optimal for the ensuing reasoning; its optimization leads to a different settings of $\lambda,\sigma$ (specifically, $\sigma=2/3$ and $\lambda=(5-\sqrt{13})/2=0.697...$ turn out to be  best), but this only impacts the value of the (implicit) universal constant factors in our results.}
\begin{equation}\label{eq:lambda sigma inclusion}
\forall x\in B_{\! L_p},\qquad M_{2\to p} \bigg(B_{ L_2} \Big(M_{p\to 2}(x),cp\alpha\sqrt{(1-\sigma)\lambda(1-\lambda)}\Big)\bigg)\subset B_{ L_p}\big(\alpha x,1-\sigma\lambda\alpha\big). 
\end{equation}
To verify~\eqref{eq:lambda sigma inclusion}, note first that by the mean value theorem there exists $\sigma\lambda\alpha\le \theta\le \lambda\alpha\le \frac{\lambda}{p}$ such that 
\begin{align*}
(1-\sigma \lambda\alpha)^p-(1-\lambda\alpha)^p&=p(1-\theta)^{p-1}(1-\sigma)\lambda\alpha\\&\ge (1-\sigma)\lambda p\alpha\left(1-\frac{\lambda}{p}\right)^{p-1}\ge \frac{\lambda}{e^\lambda} (1-\sigma)p\alpha\asymp \lambda(1-\sigma)p\alpha, 
\end{align*}
where we used the fact that $(p\ge 1)\mapsto (1-\lambda/p)^{p-1}$ is decreasing and tends to $e^{-\lambda}$ as $p\to \infty$. Therefore, if we define $r=r(\alpha,p,\sigma,\lambda)>0$ by
\begin{equation}\label{eq:def our r}
r\eqdef \sqrt{\frac{(1-\lambda)p\alpha}{5}\big((1-\sigma \lambda\alpha)^p-(1-\lambda\alpha)^p\big)},
\end{equation}
then $r\gtrsim p\alpha\sqrt{(1-\sigma)\lambda(1-\lambda)}$. Consequently, the inclusion~\eqref{eq:lambda sigma inclusion} would follow if we will show that
\begin{equation}\label{eq:lambda sigma inclusion r version}
\forall x\in B_{ L_p},\qquad M_{2\to p} \Big(B_{ L_2} \big(M_{p\to 2}(x),r\big)\Big)\subset B_{ L_p}\big(\alpha x,1-\sigma\lambda\alpha\big). 
\end{equation}
To justify~\eqref{eq:lambda sigma inclusion r version} we need to prove that $\|M_{2\to p}(w)-\alpha x\|_{\! L_p}\le 1-\sigma\lambda\alpha$ for every $w\in  L_2$ satisfying 
\begin{equation}\label{eq:w r close}
\|w-M_{p\to 2}(x)\|_{\! L_2}\le r.
\end{equation}
This indeed holds thanks to the following application of Lemma~\ref{lem:radial mazur} with $y= M_{2\to p}(w)$:
\begin{align*}\label{eq:distance from alpha w}
\begin{split}
\|M_{2\to p}(w)-\alpha x\|_{ L_p}&\stackrel{\eqref{eq:new mazur inequality}}{\le} \bigg((1-\lambda\alpha)^p+\frac{5\|w-M_{p\to 2}(x)\|_{ L_2}^2}{(1-\lambda)p\alpha} \bigg)^{\frac{1}{p}}\\&\stackrel{\eqref{eq:w r close}}{\le} \bigg((1-\lambda\alpha)^p+\frac{5r^2}{(1-\lambda)p\alpha} \bigg)^{\frac{1}{p}}\stackrel{\eqref{eq:def our r}}{=}1-\sigma\lambda\alpha. 
\end{split}
\end{align*}

Having established~\eqref{eq:lambda sigma inclusion}, we will use it  by noting that for $\alpha=\frac{1}{p}$ and $\lambda=\sigma=\frac12$  it implies that
\begin{equation}\label{eq:gamma universal version}
\forall x\in B_{ L_p},\qquad \rad_{ L_p }\Big(M_{2\to p} \big(B_{ L_2} (M_{p\to 2}(x),\gamma )\big)\Big)\le 1-\frac{1}{4p}, 
\end{equation}
where $\gamma>0$ is a universal constant. To deduce Lemma~\ref{thm:radial} from~\eqref{eq:gamma universal version}, fix  $z\in  L_p$ and $\Delta>0$, and define 
\begin{equation}\label{eq:def radial f}
 \forall x\in  L_p,\qquad f^z_\Delta(x)\eqdef K\Delta M_{p\to 2} \Big(\frac{1}{K\Delta}(x-z)\Big)\in  L_2, \qquad\mathrm{where}\qquad K=K(p)\eqdef 1+\frac{1}{4p}. 
\end{equation}
 Then,  $f^z_\Delta: L_p\to  L_2$  is a bijection whose inverse is given by 
\begin{equation}\label{eq:expression for inverse}
 \forall w\in  L_2,\qquad \big(f_\Delta^z\big)^{-1}(w)=z+ K\Delta M_{2\to p} \Big(\frac{1}{K\Delta}w\Big). 
\end{equation}
Since $(K\Delta)^{-1}(x-z)\in B_{\! L_p}$ for every $x\in B_{\! L_p}(z,K\Delta)$, it follows that 
\begin{equation}\label{eq:inclusion in our case}
\forall x\in B_{\! L_p}(z,K\Delta),\qquad  \rad_\MM\Big(\big(f_\Delta^z\big)^{-1}\big(B_{\! L_2} (f^z_\Delta(x),\gamma K\Delta)\big)\Big)\stackrel{\eqref{eq:gamma universal version}\wedge \eqref{eq:def radial f}\wedge \eqref{eq:expression for inverse}}{\le} \Big(1-\frac{1}{4p}\Big)K\Delta\stackrel{\eqref{eq:def radial f}}{<}\Delta.
\end{equation}
By~\cite[equation~(5.32)]{naor2014comparison} we have $\|M_{p\to 2}\|_{\Lip( L_p; L_2)}< p/\sqrt{2}$, so also  $\|f^z_\Delta\|_{\Lip( L_p; L_2)}< p$ by~\eqref{eq:def radial f}. Consequently, for every $x\in  L_p$ we have $B_{\! L_2} (f^z_\Delta(x),\gamma K\Delta))\supseteq B_{\! L_2} (f^z_\Delta(x),(\gamma K/p)\|f\|_{\Lip} \Delta))$, so, recalling Definition~\ref{def:radial}, we see that~\eqref{eq:inclusion in our case} implies that $ L_p$ admits a $K$-localized radially weakly bi-Lipschitz embedding into $ L_2$  with distortion $D$, where $K=K(p)$ is given in~\eqref{eq:def radial f}, so $K-1\asymp1/p$, and $D=D(p)=p/(K\gamma)\asymp p$.  
\end{proof}

\begin{remark} {\em The above proof of Lemma~\ref{thm:radial}  yields~\eqref{eq:radial mazur normalized} for $\beta=\sqrt{\sqrt[4]{e}-1}/\sqrt{5e}\in [0.14,0.15]$, as seen from~\eqref{eq:lambda sigma inclusion r version} while recalling~\eqref{eq:def our r} with our choices $\alpha=1/p$ and $\lambda=\sigma=1/2$, and using $\|M_{p\to 2}\|_{\Lip( L_p; L_2)}< p/\sqrt{2}$.  
}
\end{remark}

By a straightforward tensorization argument, Lemma~\ref{lem:radial mazur} can be deduced from its one-dimensional counterpart, which amounts to the following numerical fact:

\begin{lemma}\label{lem:numerical} If $p\ge 2$, then the following bound holds for every $0<\alpha\le \frac{1}{p}$ and every $0<\lambda<1$:
\begin{equation}\label{eq:alpha shift}
\forall u,v\in \R,\qquad \left||u+v|^{\frac{2}{p}}\sgn(u+v)-\alpha |u|^{\frac{2}{p}}\sgn(v)\right|^p\le (1-\lambda\alpha)^pu^2+\frac{5v^2}{(1-\lambda)p\alpha}.
\end{equation}
\end{lemma}

We will next explain how to quickly deduce Lemma~\ref{lem:radial mazur} from Lemma~\ref{lem:numerical}:

\begin{proof}[Proof of Lemma~\ref{lem:radial mazur} assuming Lemma~\ref{lem:numerical}] Fixing  $x,y\in  L_p$, define $u=u(x),v=v(x,y)\in  L_2$ by
\begin{equation}\label{eq:def uv}
u\eqdef M_{p\to 2}(x)\qquad \mathrm{and}\qquad  v\eqdef M_{p\to 2}(y)-M_{p\to 2}(x).
\end{equation}
Then,
\begin{align*}
\|y-\alpha x&\|_{\! L_p}^p=\int_0^1 |y(t)-\alpha x(t)|^p\ud t\\&\!\!\!\!\stackrel{\eqref{eq:def mazur map}\wedge \eqref{eq:def uv}}{=}\int_0^1 \big||u(t)+v(t)|^{\frac{2}{p}}\sgn\big(u(t)+v(t)\big)-\alpha |u(t)|^{\frac{2}{p}}\sgn\big(v(t)\big)\big|^p\ud t\\&\ \ \  \stackrel{\eqref{eq:alpha shift}}{\le} (1-\lambda\alpha)^p\|u\|_{\! L_2}^2+\frac{5}{(1-\lambda)p\alpha}\|v\|_{\! L_2}^2\\&\ \ \ \stackrel{\eqref{eq:def uv}}{=}(1-\lambda\alpha)^p\|M_{p\to 2}(x)\|_{\! L_2}^2+\frac{5}{(1-\lambda)p\alpha} \big\|M_{p\to 2}(y)-M_{p\to 2}(x)\big\|_{\! L_2}^2.
\end{align*}
This coincides with~\eqref{eq:new mazur inequality} because $\|M_{p\to 2}(x)\|_{\! L_2}=\|x\|_{\! L_p}$. 
\end{proof}

\begin{proof}[Proof of Lemma~\ref{lem:numerical}] If $u=0$, then~\eqref{eq:alpha shift} holds (with room to spare), so assume $u\neq 0$. By normalization, it suffices to prove~\eqref{eq:alpha shift} when $u=1$, i.e., our goal is equivalent to establishing the following estimate: 
\begin{equation}\label{eq:alpha shift x=1}
\forall v\in \R,\qquad \left||1+v|^{\frac{2}{p}}\sgn(1+v)-\alpha \right|^p\le (1-\lambda\alpha)^p+\frac{5v^2}{(1-\lambda)p\alpha}.
\end{equation}

Suppose first that $v\ge -1+\alpha^{p/2}$, in which case our goal~\eqref{eq:alpha shift x=1} becomes the following inequality:
\begin{equation}\label{eq:y large case}
(1+v)^{\frac{2}{p}}\le \alpha + \left((1-\lambda\alpha)^p+\frac{5v^2}{(1-\lambda)p\alpha}\right)^{\frac{1}{p}}.
\end{equation}
If also $|v|\le \sqrt{15(1-\lambda)p\alpha/20}$, which will see is when~\eqref{eq:y large case} is most meaningful, then we proceed as follows. The function $(t>0)\mapsto (1+t)^{1/t}$ is decreasing, so $(1+t)^{1/t}\ge (19/4)^{4/15}>e^{2/5}$ for $0<t\le 15/4$. Hence,  $(1+t)^{1/p}\ge e^{2t/(5p)}\ge 1+2t/(5p)$ for  $0\le t\le 15/4$. Using this for $t=5v^2/((1-\lambda)p\alpha)\le 15/4$, we get
\begin{equation}\label{eq:elementary1}
\left(1+\frac{5v^2}{(1-\lambda)p\alpha}\right)^{\frac{1}{p}}\ge 1+\frac{2v^2}{(1-\lambda)p^2\alpha}. 
\end{equation}
We can therefore bound the right hand side of~\eqref{eq:y large case} from below as follows:
\begin{align}\label{eq:finish interesting case}
\begin{split}
\alpha  + \bigg((1-\lambda\alpha)^p+\frac{5v^2}{(1-\lambda)p\alpha}\bigg)^{\frac{1}{p}}&\ge \alpha+(1-\lambda\alpha) \left(1+\frac{5v^2}{(1-\lambda)p\alpha}\right)^{\frac{1}{p}}\\&\!\!\!\!\!\stackrel{\eqref{eq:elementary1}}{\ge} 1+(1-\lambda)\alpha+\frac{2(1-\lambda\alpha)v^2}{(1-\lambda)p^2\alpha} \\&\ge 1+(1-\lambda)\alpha+\frac{v^2}{(1-\lambda)p^2\alpha}\\&=1+\frac{2}{p}v+\left(\frac{v}{p\sqrt{(1-\lambda)\alpha}}-\sqrt{(1-\lambda)\alpha}\right)^2\ge 
1+\frac{2}{p}v\\&\ge (1+v)^{\frac{2}{p}},
\end{split}
\end{align}
where the third step of~\eqref{eq:finish interesting case} holds because $1-\lambda\alpha\ge 1-1/p\ge 1/2$, as $0<\lambda<1$, $0\le \alpha\le 1/p$ and $p\ge 2$, and the final step of~\eqref{eq:finish interesting case} holds because $0<2/p\le 1$, so the function $(v>-1)\mapsto 1+2v/p-(1+y)^{2/p}$ attains its global minimum when $v=0$, where it vanishes.   We have thus completed the verification of~\eqref{eq:y large case} when $|v|\le \sqrt{15(1-\lambda)p\alpha/20}$. If, on the other hand, $\sqrt{15(1-\lambda)p\alpha/20}<|v|\le 1$, then~\eqref{eq:y large case} holds because
\begin{equation}\label{eq:yle 1}
(1+v)^{\frac{2}{p}}\le 2^{\frac{2}{p}}=\left(\frac14+\frac{15}{4}\right)^{\frac{1}{p}}\le \left(\frac1{4}+\frac{5v^2}{(1-\lambda)p\alpha}\right)^{\frac{1}{p}}<\alpha + \left((1-\lambda\alpha)^p+\frac{5v^2}{(1-\lambda)p\alpha}\right)^{\frac{1}{p}},
\end{equation}
where in the last step of~\eqref{eq:yle 1} we used that $(1-\lambda\alpha)^p\ge (1-1/p)^p\ge 1/4$, as $\lambda\alpha\le \alpha\le 1/p\le 1/2$. The remaining case of~\eqref{eq:y large case} is when $v\ge 1$ , which holds (with room to spare) because in this case we have 
\begin{equation*}\label{eq:elementary v}
(1+v)^{\frac{2}{p}}\le (2v)^{\frac{2}{p}}<\left(5v^2\right)^{\frac{1}{p}}< \alpha + \left((1-\lambda\alpha)^p+\frac{5v^2}{(1-\lambda)p\alpha}\right)^{\frac{1}{p}},
\end{equation*}
where the last step  is valid because $(1-\lambda)p\alpha<p\alpha\le 1$.

It remains to check~\eqref{eq:alpha shift x=1} for $v<-1+\alpha^{p/2}$, in which case set $w=-v>1-\alpha^{p/2}>0$ and~\eqref{eq:alpha shift x=1}  becomes: 
\begin{equation}\label{eq:negative case}
\alpha-|w-1|^{\frac{2}{p}}\sgn(1-w)\le \left((1-\lambda\alpha)^p+\frac{5w^2}{(1-\lambda)p\alpha}\right)^{\frac{1}{p}}.
\end{equation}
But~\eqref{eq:negative case} is very crude because $-|w-1|^{2/p}\sgn(1-w)< (5/ ((1-\lambda)p\alpha))^{1/p} w^{2/p}$, since $5/ ((1-\lambda)p\alpha)>5/(p\alpha)>5$, and also $\alpha\le (1- \alpha)<(1-\lambda\alpha)$, as $\alpha\le 1/p\le 1/2$. This completes the proof of Lemma~\ref{lem:numerical}.  
\end{proof}

\section{Localization and induction on scales for separated random partitions}\label{sec:induction on scales}

Here we will prove the general localization and induction on scales principle that was formulated as Lemma~\ref{lem:induction and localization}. All of the relevant definitions were provided in the Introduction. In particular, the  notion of radially bounded separating random partitions was introduced in Section~\ref{intro sketch}.

\begin{proof}[Proof of Lemma~\ref{lem:induction and localization}] We will prove that the following inequality  holds  for every $K>1$ and every $\Delta,\e>0$, even without the assumption~\eqref{decay assumption on sep} of Lemma~\ref{lem:induction and localization}, namely for any separable metric space $(\MM,d_\MM)$:
\begin{equation}\label{eq:for recursion}
\TSEP_{\Delta} (\sub; \MM) \le \frac{1}{K} \TSEP_{K\Delta} (\sub; \MM) + \sup_{z \in \MM} \TSEP_{\Delta} \big(\sub \cap B_\MM (z, K\Delta+\e); \MM\big).
\end{equation}  
Accepting the validity of~\eqref{eq:for recursion} for the moment (its justification appears below), we will next proceed to explain how to use it to quickly deduce Lemma~\ref{lem:induction and localization}.

Suppose first that $\TSEP_{\Delta} (\sub; \MM) <\infty$, which is when~\eqref{eqn:sum geometric series} is most meaningful.  Then,  $\TSEP_{\Delta'} (\sub; \MM) <\infty$ for every $\Delta'\ge \Delta$. We may therefore rearrange  the limit as $\e\to 0^+$ of~\eqref{eq:for recursion} with $\Delta$ replaced by $K^s\Delta$ for every integer $s\ge 0$ to obtain  the following recursive estimate: 
\begin{equation}\label{eq:for telescoping}
\frac{1}{K^s}\TSEP_{K^s\Delta} (\sub; \MM) - \frac{1}{K^{s+1}}\TSEP_{K^{s+1}\Delta} (\sub; \MM)  \le \frac{1}{K^s}\lim_{\e\to 0^+} \sup_{z \in \MM} \TSEP_{ K^s \Delta} \big(\sub \cap B_\MM (z, K^{s+1}\Delta+\e); \MM\big). 
\end{equation}
Thanks to the assumption~\eqref{decay assumption on sep}, by summing~\eqref{eq:for telescoping} over $s\in \N\cup \{0\}$ and telescoping we get~\eqref{eqn:sum geometric series}, as $K>1$.  

If $\TSEP_{\Delta} (\sub; \MM)=\infty$, then  we need to demonstrate that the right hand side of~\eqref{eqn:sum geometric series} is also infinite. This is so  because the assumption~\eqref{decay assumption on sep} implies in particular that there is $s\in \N$ for which $\SEP_{K^s\Delta} (\sub; \MM)<\infty$, whence also  $\TSEP_{K^s\Delta} (\sub; \MM)<\infty$ by~\eqref{eq:TSEP relations}. We can therefore consider the  largest nonnegative integer $s_0$ for which $\TSEP_{K^{s_0}\Delta} (\sub; \MM)=\infty$. By applying~\eqref{eq:for recursion} with $\Delta$ replaced by $K^{s_0}\Delta$, we see that since the left hand side of that inequality is infinite while the first term in the right hand side of that inequality is finite by the maximality of $s_0$, necessarily $\sup_{z \in \MM} \TSEP_{K^{s_0}\Delta} (\sub \cap B_\MM (z, K^{s_0+1}\Delta+\e); \MM)=\infty$ for every $\e>0$. Therefore the $s_0$-summand in the right hand side of~\eqref{eqn:sum geometric series} is infinite, as required.

To prove~\eqref{eq:for recursion}, fix $\e>0$. Define $\sigma,\tau>0$ by 
\begin{equation}\label{eq:def sigma tau bootstrap}
\sigma\eqdef \TSEP_{K\Delta} (\sub; \MM) \qquad\mathrm{and}\qquad \tau\eqdef \sup_{z \in \MM} \TSEP_{\Delta} \big(\sub \cap B_\MM (z, K\Delta+\e); \MM\big).
\end{equation}
We may assume from now that $\sigma,\tau>0$ because otherwise~\eqref{eq:for recursion} is vacuous.

Fix $\eta>0$.  The definition of $\sigma$ in~\eqref{eq:def sigma tau bootstrap} yields the existence of a probability space $(\Omega_0,\mu_0)$ and a sequence of strongly measurable mappings  $$\big\{\Phi^i:\Omega_0\to 2^\sub\big\}_{i=1}^\infty$$ 
 satisfying
\begin{equation}\label{eq:radial requirement Q}
\forall (i,\omega_0)\in \N\times \Omega_0\qquad  \rad_\MM\big(\Phi^i(\omega_0)\big)\le K\Delta,
\end{equation}
and furthermore if we define 
\begin{equation}\label{eq:def Q parititon}
\forall \omega_0\in \Omega_0,\qquad \mathcal{Q}^{\omega_0}\eqdef\big\{\Phi^i(\omega_0)\big\}_{i=1}^\infty,
\end{equation} 
then $\mathcal{Q}^{\omega_0}$ is a partition of $\sub$  for each $\omega_0\in \Omega_0$, and the following requirement holds:  
\begin{equation}\label{eq:separation requirement for Q eta}
\forall x,y\in \sub,\qquad \mu_0 \big[\omega_0\in \Omega_0:\ \mathcal{Q}^{\omega_0}(x)\neq \mathcal{Q}^{\omega_0}(y) \big]\le \frac{\sigma+\eta}{K\Delta}d_\MM(x,y). 
\end{equation}

Because  $\MM$ is separable, we can fix a sequence $\{u_n\}_{n=1}^\infty$ that is dense in $\MM$. Thanks to~\eqref{eq:radial requirement Q} we can then define a sequence of random  indices $\{n^i:\Omega_0\to \N\}_{i=1}^\infty\subset \N$ as follows: 
\begin{equation}\label{eq:def random index}
\forall \omega_0\in \Omega,\qquad n^i(\omega_0)\eqdef \min \big\{n\in \N:\ \Phi^i(\omega_0)\subset B_\MM(u_n,K\Delta+\e) \big\}.
\end{equation}
Note  in passing (so that we could use it freely later) that the measurability of all of these random indices is a quick consequence of the strong measurability of each of $\{\Phi^i:\Omega_0\to 2^\sub\}_{i=1}^\infty$. Indeed, for every $i,n\in \N$, 
\begin{align*}
\begin{split}
\big\{\omega_0\in \Omega_0:\ & n^i(\omega_0)=n\big\}\\
&=\big\{\omega_0\in \Omega_0:\ \Phi^i(\omega_0)\cap \big(\sub\setminus B_\MM(u_n,K\Delta+\e)\big)= \emptyset \big\}\bigcap \\&\qquad\quad \bigg( \bigcap_{k=1}^{n-1} \big\{\omega_0\in \Omega_0:\  \Phi^i(\omega_0)\cap \big(\sub\setminus B_\MM(u_k,K\Delta+\e)\big) \neq \emptyset \big\}\bigg).
\end{split}
\end{align*}

Next, by the definition of $\tau$ in~\eqref{eq:def sigma tau bootstrap}  for every $n\in \N$ there is a probability space $(\Omega_n,\mu_n)$ and a sequence of strongly measurable mappings

$$
\Big\{\Psi_{\!\!\! n}^j:\Omega_n\to 2^{\sub \cap B_\MM (u_n, K\Delta+\e)}\Big\}_{j=1}^\infty
$$ satisfying
\begin{equation}\label{eq:radial requirement Psi}
\forall (j,n)\in \N\times \N,\ \forall \omega_n\in \Omega_n\qquad  \rad_\MM\big(\Psi_{\!\!\! n}^j(\omega_n)\big)\le \Delta,
\end{equation} 
and furthermore if we define 
\begin{equation}\label{eq:def R partition}
\forall n\in \N,\ \forall \omega_n\in \Omega_n,\qquad \mathcal{R}_n^{\omega_n}\eqdef\big\{\Psi_{\!\!\! n}^j(\omega_n)\big\}_{i=1}^\infty,
\end{equation}
then $\mathcal{R}_n^{\omega_n}$ is a partition of $\sub \cap B_\MM (u_n, K\Delta+\e)$  for each $\omega_n\in \Omega_n$, and the following  requirement holds:  
\begin{equation}\label{eq:separation requirement for R eta}
\forall x,y\in \sub \cap B_\MM (u_n, K\Delta+\e),\qquad \mu_n \big[\omega_n\in \Omega_n:\ \mathcal{R}_n^{\omega_n}(x)\neq \mathcal{R}_n^{\omega_n}(y) \big]\le \frac{\tau+\eta}{\Delta}d_\MM(x,y). 
\end{equation}

We will henceforth work with the product space $(\Omega,\mu)$ that is given by: 
\begin{equation}\label{eq:def product spaces}
\Omega\eqdef\prod_{n=0}^\infty \Omega_n \qquad\mathrm{and}\qquad \mu\eqdef \bigotimes_{n=0}^\infty \mu_n.
\end{equation}
For every $i,j\in \N$ define $\Gamma^{i,j}:\Omega\to 2^\sub$ by:
\begin{equation}\label{eq:def common refinement}
\forall (i,j)\in \N\times \N,\ \forall \omega=(\omega_0,\omega_1,\ldots)  \in \Omega,\qquad \Gamma^{i,j}(\omega) \eqdef \Phi^i(\omega_0)\cap \Psi^{j}_{\!\!\!n^i(\omega_0)}\big(\omega_{n^i(\omega_0)}\big). 
\end{equation}
The strong measurability of $\Gamma^{i,j}$ follows from the assumed strong measurability of $\{\Psi^j_n\}_{j,n=1}^\infty$, together with the measurability of $\{n^i\}_{i=1}^\infty$ that we verified above, as a consequence of the assumed strong measurability of $\{\Phi^i\}_{i=1}^\infty$. Indeed, suppose that $E\subset \sub$ is closed. For $\omega=(\omega_0,\omega_1,\ldots)\in \Omega$, definition~\eqref{eq:def common refinement} shows that $\Gamma^{i,j}(\omega)\cap E\neq \emptyset$ if and only if  $\Psi^{j}_{\!\!\!n^i(\omega_0)}\big(\omega_{n^i(\omega_0)}\big)\cap E\neq \emptyset$, as $E\subset \sub$ and  $\{\Phi^i(\omega_0)\}_{i=1}^\infty$ is a partition of $\sub$.  Thus, 
\begin{multline*}
\big\{\omega\in \Omega:\ \Gamma^{i,j}(\omega)\cap E\neq \emptyset \big\}\\=\bigcup_{i=1}^\infty\bigcup_{j=1}^\infty \bigcup_{n=1}^\infty \big\{\omega_0\in \Omega_0:\ n^i(\omega_0)=n\big\}\times \Omega_1\times\ldots\times \Omega_{n-1}\times \big\{\omega_n\in \Omega_n:\ \Psi^{j}_{\!\!\!n}(\omega_{n})\cap E\neq \emptyset\big\}\times \prod_{k=n+1}^\infty \Omega_k.
\end{multline*}

Therefore, we can consider the random partition of $\sub$ that is given by
$$
\forall \omega=(\omega_0,\omega_1,\ldots)  \in \Omega,\qquad \Part^\omega\eqdef \big\{\Gamma^{i,j}(\omega)\big\}_{i,j=1}^\infty, 
$$
which satisfies $\rad_\MM(\Part^\omega(x))\le \Delta$ for every $x\in \sub$ thanks to~\eqref{eq:radial requirement Psi} and~\eqref{eq:def common refinement}. Finally, every $x,y\in \sub$ satisfy:   
\begin{align*}
&1-\mu\big[\omega\in \Omega:\ \Part^\omega(x)\neq \Part^\omega(y)\big]\\&\ \ \ \ \  \stackrel{\eqref{eq:def common refinement}}{=} \sum_{i=1}^\infty \sum_{j=1}^\infty \mu\big[\omega\in \Omega:\ \Phi^i(\omega_0)\supseteq \{x,y\}\ \wedge\ \Psi^{j}_{\!\!\!n^i(\omega_0)}\big(\omega_{n^i(\omega_0)}\big)\supseteq\{x,y\}\big]\\&\ \ \ \ \ \ \  = 
\sum_{i=1}^\infty \sum_{j=1}^\infty\sum_{n=1}^\infty  \mu\big[\omega\in \Omega:\ \Phi^i(\omega_0)\supseteq \{x,y\}\ \wedge\ n^i(\omega_0)=n\ \wedge \   \Psi^{j}_{\!\!\!n}(\omega_n)\supseteq\{x,y\}\big]\\ \ \ \ \ \ \ \ \ 
&\stackrel{\eqref{eq:def product spaces}}{=} \sum_{i=1}^\infty \sum_{j=1}^\infty\sum_{n=1}^\infty  \mu_0\big[\omega\in \Omega_0:\ \Phi^i(\omega_0)\supseteq \{x,y\}\ \wedge\  n^i(\omega_0)=n\big]\mu_n\big[\omega_n\in \Omega_n:\ \Psi^{j}_{\!\!\!n}(\omega_n)\supseteq\{x,y\}\big]\\
&\ \stackrel{\eqref{eq:def R partition}}{=} \sum_{n=1}^\infty \Big(1-\mu_n\big[\omega_n\in \Omega_n:\ \mathcal{R}_n^{\omega_n}(v)\neq \mathcal{R}_n^{\omega_n}(y)\big]\Big)\\&\qquad\qquad \qquad\qquad \cdot\sum_{i=1}^\infty \mu_0\big[\omega\in \Omega_0:\ \Phi^i(\omega_0)\supseteq \{x,y\}\ \wedge\ n^i(\omega_0)=n\big]\\
&\!\!\!\!\!\!\!\!\!\stackrel{\eqref{eq:separation requirement for R eta}\wedge \eqref{eq:def Q parititon}}{\ge}  \Big(1-\min\big\{1,\frac{\tau+\eta}{\Delta}d_\MM(x,y)\big\}\Big)\mu_0 \big[\omega_0\in \Omega_0:\ \mathcal{Q}^{\omega_0}(x)\neq \mathcal{Q}^{\omega_0}(y) \big]\\
&\   \stackrel{\eqref{eq:separation requirement for Q eta}}{\ge} \Big(1-\min\big\{1,\frac{\tau+\eta}{\Delta}d_\MM(x,y)\big\}\Big)\Big(1-\min\big\{1,\frac{\sigma+\eta}{K\Delta}d_\MM(x,y)\big\}\Big)\\&\ \ \ \ge 1-\frac{\frac{\sigma+\eta}{K}+\tau+\eta}{\Delta}d_\MM(x,y).
\end{align*}
By taking the limit of this estimate as $\eta\to 0^+$ we conclude that~\eqref{eq:for recursion} indeed holds.  
\end{proof}

\section{Analytic separation is preserved under localized radially weakly bi-Lipschitz embeddings}\label{sec:radial preserves sep}

We will need to impose a stronger measurability requirement from random partitions to be able to easily  use localized radially weakly bi-Lipschitz embeddings to transfer separating random partitions from one metric space to another; this is the mechanism by which  the proof of Theorem~\ref{thm:neighborhoods in Lp} will be completed.

Given a Polish metric space $(\MM,d_\MM)$ and a probability space $(\Omega,\prob)$, call a sequence of set-valued mappings 
$\{\Gamma^i:\Omega\to 2^\MM\}_{i=1}^\infty$ an analytic random partition of $\MM$ if for every $i\in \N$ and every analytic subset $A$ of $\MM$ the set $\{\omega\in \Omega:\ A\cap \Gamma^i(\omega)\}$ is $\prob$-measurable. For this, we recall the standard terminology  that a metric space is called Polish if it is separable and complete and a subset of a Polish metric space is analytic if it is a continuous image of a Polish metric space; see~\cite{Kec95} for a thorough treatment and~\cite{Lor01} for the history.

Because closed subsets of a Polish metric space are analytic, any analytic random partition is in particular a random partition per the definition that we recalled in Section~\ref{sec:sep intro}. The following theorem from~\cite{naor2024extension} provides examples of analytic random partitions with good (optimal) separation properties:

\begin{theorem}\label{thm:sep for lp spaces} For each $k\in \N$ and   $1\le p\le \infty$ there is $1\le \sigma=\sigma(p,k)\lesssim k^{\max\left\{\frac{1}{p},\frac12\right\}}$ such that for any $\Delta>0$ there exists a $\Delta$-bounded  $\sigma$-separating analytic random partition of $\ell_{\!\!p }^k$. 
\end{theorem}

The Introduction of~\cite{naor2024extension} states  Theorem~\ref{thm:sep for lp spaces} without mentioning that the corresponding random partitions are analytic, but their analyticity is stated in~\cite[Lemma~119]{naor2024extension}, which is what is applied in the proof of~\cite[Lemma~125]{naor2024extension} to derive the measurability of the random partition that Theorem~\ref{thm:sep for lp spaces} uses.\footnote{To notice that the statement of~\cite[Lemma~119]{naor2024extension} provides the measurability that we need, recall the important classical theorem of Luzin~\cite{Luz17}  (see also e.g.~\cite[Theorem~21.10]{Kec95})  that analytic sets are universally measurable, i.e., they are measurable with respect to every complete $\sigma$-finite Borel measure on the given Polish metric space.} We will use below only the case $p=2$ of Theorem~\ref{thm:sep for lp spaces}, for which the underlying construction for finite subsets of $\ell_{\!\!2}^k$ is due to~\cite{ccg98}, and its extension to random partitions of all of $\ell_{\!\!2}^k$  is due to~\cite{LN05,naor2024extension}. The short proof of following basic and useful lemma  clarifies why it is beneficial to consider analytic random partitions:

\begin{lemma}\label{lem:pullback} Fix $\sigma,\Delta,R,L>0$. Let $(\NN,d_\NN)$ be a Polish metric space admitting an analytic random partition that is  $(LR)$-bounded and $\sigma$-separating.   Suppose that $(\MM,d_\MM)$ is a Polish metric space, $\cS\subset \MM$ is a Borel subset of $\MM$, and that $\f:\cS\to \NN$  is an $L$-Lipschitz function that satisfies   the following property: 
\begin{equation}\label{eq:condition in analytic lemma}
\forall x\in \cS,\qquad \rad_\MM\Big(\f^{-1} \big(B_{\NN}(\f(x),LR)\big)\Big) \le \Delta. 
\end{equation}
Then, $(\cS,d_\MM)$ admits an analytic random partition that is  $\Delta$-radially bounded and $\sigma\frac{\Delta}{R}$-separating. Hence,
$$
\TSEP_\Delta(\cS;\MM)\le \sigma \frac{\Delta}{R}. 
$$
\end{lemma}

\begin{proof} Fix a probability space $(\Omega,\mu)$ and an analytic random partition 
\begin{equation}\label{eq:take analytic partition}
\Part=\big\{\Gamma^i:\Omega\to 2^{\NN}\big\}_{i=1}^\infty
\end{equation}
that is $(LR)$-bounded and $\sigma$-separating. If $E\subset \cS$ is analytic, then  $\f(E)\subset \NN$ is analytic as $\f$ is  continuous, whence by the assumed analyticity of the random partition~\eqref{eq:take analytic partition} for every $i\in \N$ the set 
$$
\big\{\omega\in \Omega:\ \f^{-1} \big(\Gamma^i(\omega)\big)\cap E\neq \emptyset\big\}=\big\{\omega\in \Omega:\ \Gamma^i(\omega)\cap\f(E)\neq \emptyset\big\}
$$
is $\mu$-measurable.  So, the following sequence of set-valued mappings  is an analytic random partition of $\cS$:
\begin{equation}\label{eq:the pullback of analytic}
\mathcal{Q}=\big\{(\omega\in \Omega)\mapsto \mathcal{Q}^\omega\big\}\eqdef \big\{(\omega\in \Omega)\mapsto \f^{-1}\big(\Gamma^i(\omega)\big)\subset \cS\big\}_{i=1}^\infty, 
\end{equation}
A different way to write~\eqref{eq:the pullback of analytic} is:
\begin{equation}\label{eq:easier notation Q}
\forall \omega\in \Omega,\ \forall x\in \cS,\qquad \mathcal{Q}^\omega(x)\stackrel{\eqref{eq:take analytic partition}\wedge\eqref{eq:the pullback of analytic}}{=}\f^{-1}\Big(\Part^\omega\big(\f(x)\big)\Big).
\end{equation}

Because $\Part$ is $(RL)$-bounded by assumption, $\diam_\NN(\Part^\omega(\f(x)))\le RL$ for every $\omega\in \Omega$ and $x\in \MM$, whence  $\Part^\omega(\f(x))\subset B_\NN(\f(x),RL)$. Thanks to~\eqref{eq:easier notation Q}, this implies that $\mathcal{Q}^\omega(x)$ is contained in $\f^{-1}(B_\NN(\f(x),RL))$. By invoking the assumption~\eqref{eq:condition in analytic lemma} we conclude that the random partition $\mathcal{Q}$ is $\Delta$-radially bounded with respect to the super-space $(\MM,d_\MM)$.  Finally, for every $x,y\in \cS$ we have 
\begin{align}\label{eq:sigma/R}
\begin{split}
\mu\big[\omega\in \Omega:\  \mathcal{Q}^\omega(x)\neq \mathcal{Q}^\omega(y)\big]&\stackrel{\eqref{eq:easier notation Q}}{=}\mu\Big[\omega\in \Omega:\  \Part^\omega\big(\f(x)\big)\neq \Part^\omega\big(\f(y)\big)\Big]\\& \qquad\qquad \le \sigma \frac{d_\NN\big(\f(x),\f(y)\big)}{LR} \le \Big(\frac{\sigma\Delta}{R}\Big) \frac{d_\MM(x,y)}{\Delta}, 
\end{split}
\end{align}
where the second step of~\eqref{eq:sigma/R} uses the assumption that the random partition $\Part$ is $\sigma$-separating and $(LR)$-bounded,   and the third step of~\eqref{eq:sigma/R}  uses the assumption that  $\f$ is $L$-Lipschitz. We already checked that $\mathcal{Q}$ is $\Delta$-radially bounded with respect to $(\MM,d_\MM)$, so the required estimate $\TSEP_\Delta(\cS;\MM)\le \sigma\Delta/R$ follows from~\eqref{eq:sigma/R}. \end{proof}

The following lemma sets the stage for our subsequent application of Lemma~\ref{lem:pullback}; its short proof proceed by combining   two important Euclidean results, namely the Kirszbraun Lipschitz extension theorem~\cite{Kir34} and  the Johnson--Lindenstrauss dimension lemma~\cite{JL82}.

\begin{lemma}\label{lem:JL+K} For any $\sub\subset  L_2$ with $2\le |\sub|<\infty$ there is an integer $1\le k\lesssim \log |\sub|$ and $H: L_2\to \ell_{\!\! 2 }^k$ satisfying:
\begin{equation}\label{eq:additive gpal}
\forall x,y\in  L_2,\qquad \frac12 \|x-y\|_{\! L_2}-\frac32 d_{ L_2}(x,\sub)-\frac32 d_{ L_2}(y,\sub)\le \|H(x)-H(y)\|_{\!\ell_{\!\!2 }^k}\le \|x-y\|_{\! L_2}.
\end{equation}
Therefore, the following inclusion holds for any $r>0$ and any point $x$ in the $r$-neighborhood $B_{\! L_2}(\sub,r)$ of $\sub$:
\begin{equation}\label{eq:inclusion for us dim reduction}
B_{ L_2}(\sub,r)\cap H^{-1} \big(B_{\!\ell_{\!\!2 }^k}(H(x),r)\big)\subset B_{ L_2}(x,8r).
\end{equation}
\end{lemma}

For~\eqref{eq:inclusion for us dim reduction}, recall our notation~\eqref{eq:def of r neighborhood notation} for neighborhoods of subsets in a metric space.

\begin{proof}[Proof of Lemma~\ref{lem:JL+K}] \eqref{eq:inclusion for us dim reduction} follows from the first inequality in~\eqref{eq:additive gpal}. Indeed, fix $r>0$ and $x\in B_{ L_2}(\sub,r)$, i.e., $d_{ L_2}(x,\sub)\le r$. Consider any point $y$ that belongs to the set that appears in the left hand side of\eqref{eq:inclusion for us dim reduction}. Thus,  $y\in B_{ L_2}(\sub,r)$, i.e., $d_{ L_2}(y,\sub)\le r$, and $H(y)\in B_{\!\ell_{\!\!2 }^k}(H(x),r)$, i.e. $\|H(x)-H(y)\|_{\! L_2}\le r$. By substituting  these $3$ bounds into the first inequality in~\eqref{eq:additive gpal} and rearranging we arrive at $\|x-y\|_{\! L_2}\le 8r$, i.e., $y\in B_{ L_2}(x,8r)$, as required. Lemma~\ref{lem:JL+K} will therefore be proven once we establish~\eqref{lem:JL+K}, which we will proceed to do next.

The Johnson--Lindenstrauss lemma~\cite{JL82} yields a positive integer $k\lesssim \log |\sub|$ and $h: L_2\to \ell_{\!\!2 }^k$ such that 
\begin{equation}\label{eq:use JLh}
\forall x,y\in \sub,\qquad \frac12\|x-y\|_{\! L_2}\le \|h(x)-h(y)\|_{\!\ell_{\!\!2 }^k}\le \|x-y\|_{\! L_2}.
\end{equation}
Kirszbraun's Lipschitz extension theorem~\cite{Kir34} provides a function $H: L_2\to \ell_{\!\!2 }^k$ satisfying:
\begin{equation}\label{eq:H extends}
\forall a\in \sub,\quad H(a)=h(a)\qquad\mathrm{and}\qquad \forall x,y\in  L_2,\quad \|H(x)-H(y)\|_{\!\ell_{\!\!2 }^k}\le \|x-y\|_{\! L_2}.
\end{equation}

Given $x,y\in  L_2$, the second inequality in~\eqref{eq:additive gpal}, i.e., the fact that $H$ is $1$-Lipschitz,  coincides with the second condition in~\eqref{eq:H extends}. For the first inequality in~\eqref{eq:additive gpal} fix $a,b\in  L_2$ with
\begin{equation}\label{eq:nearest points state}
a,b\in \sub\qquad\mathrm{and}\qquad \|x-a\|_{\! L_2}=d_{ L_2}(x,\sub)\qquad\mathrm{and}\qquad \|y-b\|_{\! L_2}=d_{ L_2}(y,\sub).
\end{equation}
Then,
\begin{eqnarray*}\label{eq:to rearrange JL}
\|x-y\|_{\! L_2}&\le & \|x-a\|_{\! L_2}+ \|a-b\|_{\! L_2}+\|b-y\|_{\! L_2}\\&\stackrel{\eqref{eq:use JLh}\wedge \eqref{eq:nearest points state}}{\le}& d_{ L_2}(x,\sub)+2\|h(a)-h(b)\|_{\!\ell_{\!\!2 }^k} +d_{ L_2}(y,\sub)\\&\stackrel{\eqref{eq:H extends}\wedge \eqref{eq:nearest points state}}{=}& d_{ L_2}(x,\sub)+2\|H(a)-H(b)\|_{\!\ell_{\!\!2 }^k} +d_{ L_2}(y,\sub)\\
&\le& d_{ L_2}(x,\sub)+d_{ L_2}(y,\sub) \\&& \qquad\quad +2\big(\|H(a)-H(x)\|_{\!\ell_{\!\!2 }^k}+\|H(x)-H(y)\|_{\!\ell_{\!\!2 }^k}+\|H(y)-H(b)\|_{\!\ell_{\!\!2 }^k}\big)\\
&\stackrel{\eqref{eq:H extends}}{\le}& d_{ L_2}(x,\sub)+2\big(\|a-x\|_{\! L_2}+\|H(x)-H(y)\|_{\!\ell_{\!\!2 }^k}+\|y-b\|_{\! L_2}\big) +d_{ L_2}(y,\sub)\\
&\stackrel{\eqref{eq:nearest points state}}{=}& d_{ L_2}(x,\sub)+2\big(d_{ L_2}(x,\sub)+\|H(x)-H(y)\|_{\!\ell_{\!\!2 }^k}+d_{ L_2}(x,\sub)\big) +d_{ L_2}(y,\sub).
\end{eqnarray*}
which rearranges to give the desired lower bound on $\|H(x)-H(y)\|_{\!\ell_{\!\!2 }^k}$ in the first inequality of~\eqref{eq:additive gpal}.
\end{proof}

\begin{lemma}\label{lem:weakly radiale plus dim reduction} Fix $K,D\ge  1$. Suppose that  $(\MM,d_\MM)$ is a metric space that admits a $K$-localized radially weakly bi-Lipschitz embedding into $ L_2$ with distortion $D$.  If $\sub\subset \MM$ satisfies $2\le |\sub|<\infty$ and also
\begin{equation}\label{eq:sub is in ball}
\sub\subset B_\MM\Big(z,\Big(K-\frac{1}{8D}\Big)\Delta\Big)
\end{equation}
for some $\Delta>0$ and $z\in \MM$, then there is an integer $k=O(\log |\sub|)$  a $1$-Lipschitz function 
$$\f:B_\MM\Big(\sub,\frac{1}{8D}\Delta\Big)\to \ell_{\!\!2 }^k$$ 
such that
\begin{equation}\label{eq:reduced dim inclusion}
\forall x\in B_\MM\Big(\sub,\frac{1}{8D}\Delta\Big),\qquad \rad_\MM\Bigg(\f^{-1} \bigg(B_{\ell_{\!\!2 }^k}\Big(\f(x),\frac{1}{8D}\Delta\Big)\bigg)\Bigg) <\Delta. 
\end{equation}
\end{lemma}

\begin{proof} Writing $r\eqdef \Delta/(8D)$, we have $B_\MM(\sub,r)\subset B_\MM(z,K\Delta)$  thanks to~\eqref{eq:sub is in ball} and the triangle inequality for $d_\MM$. Hence,
recalling Definition~\ref{def:radial}, the assumption of Lemma~\ref{lem:weakly radiale plus dim reduction} implies that there is a $1$-Lipschitz function $f:B_\MM(\sub,r)\to  L_2$ and $0<R<\Delta$ such that for any $x\in B_\MM(\sub,r)$ there is $y_x\in \MM$ satisfying:  
\begin{equation}\label{eq:use radially weak f}
f^{-1}\big(B_{ L_2}(f(x),8r)\big)=f^{-1} \bigg(B_{ L_2}\Big(f(x),\frac{1}{D}\Delta\Big)\bigg)\subset B_\MM(y_x,R). 
\end{equation}

By applying  Lemma~\ref{lem:JL+K} to $f(\sub)$ we get an integer $1\le k\lesssim \log |\sub|$ and  a $1$-Lipschitz   function $H: L_2\to \ell_{\!\! 2 }^k$ that satisfies the following inclusion for any $v\in B_{ L_2}(f(\sub),r)$: 
$$
 B_{ L_2}(f(\sub),r)\cap H^{-1} \big(B_{\!\ell_{\!\!2 }^k}(H(v),r)\big)\subset B_{ L_2}(v,8r). 
$$
Since $f$ is $1$-Lipschitz, $f(x)\in B_{ L_2}(f(\sub),r)$ for every $v\in B_\MM(\sub,r)$, so the following holds a special case: 
\begin{equation}\label{eq:use the JL+K inclusion}
\forall x\in B_\MM(\sub,r),\qquad  B_{ L_2}(f(\sub),r)\cap H^{-1} \Big(B_{\!\ell_{\!\!2 }^k}\big(\f(x),r\big)\Big)\subset B_{ L_2}(f(x),8r),
\end{equation}
where we define $\f\eqdef H\circ f:B_\MM(\sub,r)\to \ell_{\!\!2 }^k$. Then, $\f$ is $1$-Lipschitz as both $f$ and $H$ are $1$-Lipschitz. Observe furthermore that as $f$ is $1$-Lipschitz and its domain is $B_\MM(\sub,r)$,  all of its values belong to $B_{ L_2}(f(\sub),r)$. Therefore~\eqref{eq:use the JL+K inclusion} can we rewritten as follows: 
\begin{equation}\label{eq:use the JL+K inclusion 2}
\forall x\in B_\MM(\sub,r),\qquad  f\big(B_{\MM}(f(\sub),r)\big)\cap H^{-1} \Big(B_{\ell_{\!\!2 }^k}\big(\f(x),r\big)\Big)\subset B_{ L_2}(f(x),8r),
\end{equation}
By applying $f^{-1}$ to both sides of~\eqref{eq:use the JL+K inclusion 2} and then using~\eqref{eq:use radially weak f}, we conclude that  
\begin{align*}
\forall x\in  B_\MM\Big(\sub,\frac{1}{8D}&\Delta\Big)=B_\MM(\sub,r),\\& \f^{-1} \bigg(B_{\ell_{\!\!2 }^k}\Big(\f(x),\frac{1}{8D}\Delta\Big)\bigg)=f^{-1}
\bigg(H^{-1}\Big(B_{\ell_{\!\!2 }^k}\big(\f(x),r\big)\Big)\bigg)\subset B_\MM(y_x,R). 
\end{align*}
The $d_\MM$-radius of $\f^{-1} (B_{\ell_{\!\!2 }^k}(\f(x),\delta/(8D)))$ is therefore at most $R<\Delta$ if $x\in B_\MM(\sub,\Delta/(8D))$. In other words, the desired conclusion~\eqref{eq:reduced dim inclusion} of Lemma~\ref{lem:weakly radiale plus dim reduction}  indeed holds.
\end{proof}

The following theorem implies Theorem~\ref{thm:neighborhoods in Lp} because if $\MM=L_p$ for some $p>2$, then by Lemma~\ref{thm:radial} its assumptions hold with $D\asymp p$ and $K-1\asymp 1/p$.

\begin{theorem}\label{thm:get sep of neioghborhood doubling} Fix $K,D,\ul>1$. Let $(\MM,d_\MM)$ be a Polish metric space admitting a $K$-localized radially weakly bi-Lipschitz embedding into $ L_2$ with distortion $D$. Then, every $\ul$-doubling Borel subset $\cD$ of $\MM$  satisfies:
\begin{equation}\label{eq:sep hat of neighborhood doubling-cases}
\forall \Delta>0,\qquad \TSEP_\Delta \Big(B_\MM\big(\cD,\frac{1}{9D}\Delta\big);\MM\Big)\lesssim \frac{{\textstyle \sqrt{\log \ul}}}{K-1}\cdot \left\{\begin{array}{ll} \frac{\sqrt{\log \frac{1}{K-1}}}{K-1}&\mathrm{if}\ 1<K\le 1+\frac{1}{D},\\K D\sqrt{\log (KD)}&\mathrm{if}\ K>1+\frac{1}{D}.\\\end{array}\right.
\end{equation}
If furthermore $2\le |\cD|<\infty$, then 
\begin{equation}\label{eq:finite cD case-cases}
\forall \Delta>0,\qquad  \TSEP_\Delta \Big(B_\MM\big(\cD,\frac{1}{9D}\Delta\big);\MM\Big)\lesssim \frac{{\textstyle \sqrt{\log |\cD|}}}{K-1}\cdot \left\{\begin{array}{ll} \frac{1}{K-1}&\mathrm{if}\ 1<K\le 1+\frac{1}{D},\\K D&\mathrm{if}\ K>1+\frac{1}{D}.\\\end{array}\right.
\end{equation}
\end{theorem}

Note that  if $|\cD|=n<\infty$, then conclusion~\eqref{eq:sep hat of neighborhood doubling} of Theorem~\ref{thm:get sep of neioghborhood doubling} is stronger than its conclusion~\eqref{eq:finite cD case} unless the doubling constant $\ul$ is very large; specifically, this occurs if and only if
$$
\ul =n^{o\left(\frac{1}{\log (eKD)}\right)}.
$$    

\begin{proof}[Proof of Theorem~\ref{thm:get sep of neioghborhood doubling} ] It suffices to prove  Theorem~\ref{thm:get sep of neioghborhood doubling}  while imposing the further  assumption $K\ge 1+1/D$ because if $1<K<1+1/D<2$, then we can replace $D$ by the larger quantity $1/(K-1)$. Thus, we will assume from now that $K\ge 1+1/D$ and our goal becomes to prove that 
\begin{equation}\label{eq:sep hat of neighborhood doubling}
\TSEP_\Delta \Big(B_\MM\big(\cD,\frac{1}{9D}\Delta\big);\MM\Big)\lesssim \frac{KD\sqrt{\log(KD)}}{K-1}{\textstyle \sqrt{\log \ul}},
\end{equation}
and correspondingly if $2\le |\cD|<\infty$, then 
\begin{equation}\label{eq:finite cD case}
\TSEP_\Delta \Big(B_\MM\big(\cD,\frac{1}{9D}\Delta\big);\MM\Big)\lesssim \frac{KD}{K-1}{\textstyle \sqrt{\log |\cD|}}.
\end{equation}

We will start by proving~\eqref{eq:sep hat of neighborhood doubling}. Fix $\Delta>0$. Our goal is to eventually  apply Lemma~\ref{lem:induction and localization} with $\sub$ replaced by $B_\MM(\cD,\Delta/(9D))$, and with $K$ replaced by $K_*$, where for convenience we set the following notation:
\begin{equation}\label{eq:def Kstar}
K_*\eqdef K-\frac{1}{2D}.
\end{equation}
Then, our assumption $K\ge 1+1/D$  ensures that 
\begin{equation}\label{eq:K star minus 1}
K_*-1\asymp K-1.
\end{equation}
In particular, $K_*>1$, so this  will be a valid instantiation  of Lemma~\ref{lem:induction and localization}.

Fix $s\ge 0$  and define three auxiliary parameters $\alpha=\alpha(s,\Delta,K,D), \beta=\beta(s,\Delta,K,D),\e=\e(s,\Delta,K,D)$ by:
\begin{equation}\label{eq:def Delta prime}
\alpha\eqdef \Big(K-\frac{1}{8D}\Big)K_*^s\Delta  \qquad\mathrm{and}\qquad  \beta\eqdef \frac{1}{8D}K_*^s\Delta -\frac{1}{9D}\Delta\qquad\mathrm{and}\qquad \e_0\eqdef \frac{\frac{3}{8}K_*^s-\frac{1}{9}}{D}\Delta. 
\end{equation}
Observe for later use that because $s\ge 0$,  $K>K_*>1$ and $D>1$, we have:  
\begin{equation}\label{eq:alpha over beta}
\alpha\asymp KK_*^{s}\Delta\qquad\mathrm{and}\qquad \alpha>\beta\asymp \frac{K_*^s}{D}\Delta\asymp\frac{\alpha}{KD}\qquad\mathrm{and}\qquad \e_0\asymp\frac{K_*^s}{D}\Delta.
\end{equation} 
Fix $z\in \MM$. By a standard iteration of the assumed $\ul$-doubling property of $\cD$, there exists a finite subset  $\sub$ of $\cD\cap B_\MM(z,\alpha)$ which is $\beta$-dense in $\cD\cap B_\MM(z,\alpha)$, i.e., $B_\MM(\sub,\beta)\supseteq \cD\cap B_\MM(z,\alpha)$, and  whose size satisfies:

\begin{equation}\label{eq:net size doubling}
2\le |\sub|\le \ul^{\left\lceil \log_2 \left(\frac{2\alpha}{\beta}\right)\right\rceil}. 
\end{equation}
(Briefly, for completeness (see also~\cite[Lemma~4.1.12]{HK15}): If $\ell\in \N$ satisfies $2^\ell\ge 2\alpha/\beta$, then  iterate $\ell$ times the doubling condition for $\cD$ to get the existence of  a cover $\cD\cap B_\MM(z,\alpha)$ by $\ul^\ell$ balls of radius $\alpha/2^\ell\le \beta/2$, and then let $\sub$ consist of  one point from the intersection of each of these balls with $\cD\cap B_\MM(z,\alpha)$.)
Observe that by the definition~\eqref{eq:def Delta prime} of $\alpha$, the aforementioned inclusion  $\sub\subset \cD\cap B_\MM(z,\alpha)$ can be rewritten as 
\begin{equation}\label{eq:display two inclusions for using lemma}
\sub \subset \cD\cap B_\MM\Big(z,\Big(K-\frac{1}{8D}\Big)K_*^s\Delta\Big).
\end{equation}
Also, the  parameters in~\eqref{eq:def Kstar} and~\eqref{eq:def Delta prime} were chosen judiciously so that  following inclusion holds:
\begin{equation}\label{neighborhood in ball inclusion}
B_\MM\big(\cD,\frac{1}{9D}\Delta\big)\cap  B_\MM (z, K_*^{s+1}\Delta+\e_0)\subset B_\MM\big(\sub,\frac{1}{8D}K_*^s\Delta\big).
\end{equation}
Indeed, if  $x\in \MM$ satisfies   $d_\MM(x,y)\le \Delta/(9D)$ for some $y\in \cD$, and also $d_\MM(x,z)\le K_*^{s+1}\Delta+\e_0$, then 
\begin{equation*}
d_\MM(y,z)\le d_\MM(x,y)+d_\MM(x,z)\le \frac{1}{9D}\Delta+ K_*^{s+1}\Delta+\e_0\stackrel{\eqref{eq:def Kstar}\wedge \eqref{eq:def Delta prime}}{=} \alpha.
\end{equation*}
Therefore, $y\in \cD\cap B_\MM(z,\alpha)$. As $\sub$ is $\beta$-dense in $\cD\cap B_\MM(z,\alpha)$, there exists $c\in \sub$ with $d_\MM(c,y)\le \beta$, whence: 
$$d_\MM(x,c)\le d_\MM(x,y)+d_\MM(y,c)\le \frac{1}{9D}\Delta+\beta\stackrel{\eqref{eq:def Delta prime}}{=}\frac{1}{8D}K_*^s\Delta.$$ 
This implies that $x$ belongs to the right hand side of~\eqref{neighborhood in ball inclusion}, as required.

Thanks to~\eqref{eq:display two inclusions for using lemma}, we may invoke  Lemma~\ref{lem:weakly radiale plus dim reduction} with $\Delta$ replaced by $K_*^s\Delta$ to get an integer $k$ satisfying
\begin{equation}\label{eq:after dim reduction bound doubling}
k\asymp \log |\sub|\stackrel{\eqref{eq:net size doubling}}{\lesssim} (\log \ul)\log\left(\frac{\alpha}{\beta}\right)\stackrel{\eqref{eq:alpha over beta}}{\lesssim} (\log \ul)\log(KD),
\end{equation}
and a $1$-Lipschitz function $\f=\f_{s,\e}:B_\MM\big(\cD,\frac{1}{9D}\Delta\big)\cap  B_\MM (z, K_*^{s+1}\Delta+\e_0)\to \ell_{\!\!2 }^k$, such that
\begin{equation}\label{eq:reduced dim inclusion s eps}
\forall x\in B_\MM\big(\cD,\frac{1}{9D}\Delta\big)\cap  B_\MM (z, K_*^{s+1}\Delta+\e_0),\qquad \rad_\MM\Bigg(\f^{-1} \bigg(B_{\ell_{\!\!2 }^k}\Big(\f(x),\frac{1}{8D}K_*^s\Delta\Big)\bigg)\Bigg) <K_*^s\Delta. 
\end{equation}

By~\eqref{eq:reduced dim inclusion s eps} we may apply Lemma~\ref{lem:pullback} to $\cS=B_\MM(\cD,\Delta/(9D))\cap  B_\MM (z, K_*^{s+1}\Delta+\e_0)\subset \MM$ and the target space $\NN=\ell_{\!\!2 }^k$, for which the assumption of Lemma~\ref{lem:pullback} holds for $\sigma\lesssim \sqrt{K}$ by Theorem~\ref{thm:sep for lp spaces}, with the parameters  $L=1$,  $R=K_*^s\Delta/(8D)$ and  $\Delta$ replaced by $K_*^s\Delta$, to get that for every $z\in \MM$ we have: 
\begin{equation*}
 \TSEP_{K_*^s\Delta} \Big(B_\MM\big(\cD,\frac{1}{9D}\Delta\big)\cap  B_\MM (z, K_*^{s+1}\Delta+\e_0);\MM\Big) \lesssim D{\textstyle \sqrt{k}}\stackrel{\eqref{eq:after dim reduction bound doubling}}{\lesssim} D{\textstyle \sqrt{(\log \ul)\log(KD)}}.
\end{equation*}
Consequently,
\begin{align}\label{eq:foir substitution into induction on scales}
\begin{split}
\lim_{\e\to 0^+} \sup_{z \in \MM} &\TSEP_{ K_*^s \Delta} \Big(B_\MM\big(\cD,\frac{1}{9D}\Delta\big) \cap B_\MM (z, K_*^{s+1}\Delta+\e); \MM\Big)\\&\le \sup_{z\in \MM} \TSEP_{K_*^s\Delta} \Big(B_\MM\big(\cD,\frac{1}{9D}\Delta\big)\cap  B_\MM (z, K_*^{s+1}\Delta+\e_0);\MM\Big)\\&\lesssim D{\textstyle \sqrt{(\log \ul)\log(KD)}}.
\end{split}
\end{align}

In order to be able to use Lemma~\ref{lem:induction and localization} in conjunction with~\eqref{eq:foir substitution into induction on scales}, we must first check that its assumption~\eqref{decay assumption on sep} holds with $\sub$ replaced by $B_\MM(\cD,\Delta/(9D))$, i.e., we need to verify that:
\begin{equation}\label{eq:sep of neighborhood}
  \lim_{T \to \infty} \frac{1}{T}\SEP_{T} \Big(B_\MM\big(\cD,\frac{1}{9D}\Delta\big)\Big) =0.
\end{equation}
The ensuing justification of~\eqref{eq:sep of neighborhood} is suboptimal from the quantitative perspective. We chose this route as it is quick and for our purposes a qualitative statement (namely, without providing a rate of convergence) such as~\eqref{eq:sep of neighborhood} suffices. See Remark~\ref{Rem:sep of neighborhood optimal} below for an improved (optimal) statement.

A well-known classical result (contained in~\cite{Ass83}; see also~\cite[Chapter~12]{Hei01} and the discussion in~\cite{DLP13}, or e.g.~the proof of~\cite[Theorem~5.2]{NS11}) asserts that there are $k=k(\ul)\in \N$ and $0<\eta=\eta(\ul)\le 1$ such that for any $T>0$ there exists a function  $f_T:\MM\to \ell_{\!\!2 }^k$ which satisfies: 
\begin{equation}\label{eq:weakly defined on M}
\left\{\begin{array}{ll} \|f_T\|_{\Lip(\MM;\ell_{\!\!2 }^k)}\le 1,\\  \forall a,b\in \cD,\qquad d_\MM(a,b)\ge \frac14 T \implies d_\MM\big(f_T(a),f_T(b)\big)\ge \eta T.\end{array}\right. 
\end{equation}

We will next proceed to demonstrate that the following inclusion holds: 
\begin{align}\label{eq:inclusion for f_tau}
\begin{split}
\forall T\ge \frac{\Delta}{2\eta D},\forall x\in  &B_\MM\big(\cD,\frac{1}{9D}\Delta\big),\\  &B_\MM\big(\cD,\frac{1}{9D}\Delta\big)\cap f_T^{-1}\Big( B_{\!\ell_{\!\!2 }^k}\big(f_T(x),\frac{\eta}{2}T\big)\Big)\subset B_\MM\big(x,\frac12 T\big).
\end{split}
\end{align}
After~\eqref{eq:inclusion for f_tau} will be established, we will proceed by fixing $T\ge \Delta/(2\eta D)$ and using Lemma~\ref{lem:pullback} with both $\MM$ and $\cS$ replaced by $B_\MM(\cD,\Delta/(9D))$, the function $\f$ replaced by the restriction of $f_T$ to $B_\MM(\cD,\Delta/(9D))$, the target space $\NN=\ell_{\!\!2 }^k$, for which the assumption of Lemma~\ref{lem:pullback} holds for $\sigma\lesssim \sqrt{K}$ by Theorem~\ref{thm:sep for lp spaces}, and with the parameters  $L=1\ge \|f_T\|_{\Lip(\MM;\ell_{\!\!2 }^k)}$,  $R=\eta T/2$ and  $\Delta$ replaced by $T/2$, to get that: 
\begin{equation}\label{eq:sep of neighborhood of doubling bad bound}
\forall T\ge \frac{\Delta}{\eta D},\qquad  \SEP_T\Big(B_\MM\big(\cD,\frac{1}{9D}\Delta\big)\Big) \stackrel{\eqref{eq:TSEP relations}}{\le} 2\TSEP_{\frac12T}\Big(B_\MM\big(\cD,\frac{1}{9D}\Delta\big)\Big) \lesssim \frac{\sqrt{k}}{\eta}\lesssim_\ul 1,
\end{equation}
which handily implies~\eqref{eq:sep of neighborhood}. 

Thus, it remains to prove~\eqref{eq:inclusion for f_tau}, which we will next do by contradiction. The contrapositive of~\eqref{eq:inclusion for f_tau} is that there exist  $T>0$, as well as   $x,y\in \MM$ and $a,b\in \cD$ that satisfy:
\begin{equation}\label{eq:contrapositive ab}
\left\{\begin{array}{ll} T\ge \frac{\Delta }{\eta D},\\ \|f_T(x)-f_T(y)\|_{\!\ell_{\!\!2 }^k}\le \frac{\eta T}{2},\\ d_\MM(x,a),d_\MM(y,b)\le \frac{\Delta}{9D},\\  d_\MM(x,y)>\frac{T}{2}. \end{array}\right.
\end{equation}
It follows in particular from~\eqref{eq:contrapositive ab} that $d_\MM(a,b)$ is sufficiently large so that the second part of~\eqref{eq:weakly defined on M} applies:  
\begin{equation}\label{eq:a,b are indeed far}
d_\MM(a,b)\ge d_\MM(x,y)-d_\MM(x,a)-d_\MM(y,b)\stackrel{\eqref{eq:contrapositive ab}}{>} \frac{T}2-\frac{2\Delta }{9D}\stackrel{\eqref{eq:contrapositive ab}}{\ge} \frac{T}{2} -\frac{2\eta T}{9}>\frac{T}4.
\end{equation}
We may therefore use~\eqref{eq:weakly defined on M} to deduce the following contradictory chain of inequalities: 
\begin{align*}
 \frac{\eta T}{2}  &\stackrel{\eqref{eq:contrapositive ab}}{\ge }  \|f_T(x)-f_T(y)\|_{\!\ell_{\!\!2 }^k} \\&\ge \|f_T(a)-f_T(b)\|_{\!\ell_{\!\!2 }^k}-\|f_T(x)-f_T(a)\|_{\!\ell_{\!\!2 }^k}-\|f_T(y)-f_T(b)\|_{\!\ell_{\!\!2 }^k} \\&\ge \|f_T(a)-f_T(b)\|_{\!\ell_{\!\!2 }^k}-\|f_T\|_{\Lip(\MM;\ell_{\!\!2 }^k)}\big(d_\MM(x,a)+d_\MM(x,b)\big)\\&
\stackrel{\eqref{eq:a,b are indeed far}\wedge \eqref{eq:weakly defined on M}}{\ge} \eta T-\frac{2\Delta}{9D} \stackrel{\eqref{eq:contrapositive ab}}{\ge }\eta T-\frac{4\eta T}{9}.  
\end{align*}

Having checked the assumption~\eqref{decay assumption on sep} of Lemma~\ref{lem:induction and localization} with $\sub$ replaced by $B_\MM(\cD,\Delta/(9D))$ is indeed satisfied, substituting~\eqref{eq:foir substitution into induction on scales} into the conclusion~\eqref{eqn:sum geometric series} of  Lemma~\ref{lem:induction and localization}  and  $K$ replaced by $K_*$ gives: 
$$
\TSEP_{\Delta} \Big(B_\MM\big(\cD,\frac{1}{9D}\Delta\big);\MM\Big)\lesssim D{\textstyle \sqrt{(\log \ul)\log(KD)}}\sum_{s=0}^\infty \frac{1}{K_*^s}=\frac{K_*D\sqrt{\log (KD)}}{K_*-1}{\textstyle{\log \ul}},
$$
which implies the desired estimate~\eqref{eq:sep hat of neighborhood doubling} by the definition~\eqref{eq:def Kstar} of $K_*$ and~\eqref{eq:K star minus 1}.

The justification of the finitary variant~\eqref{eq:finite cD case} is identical to the above reasoning, except that in this case one can simply work with $\sub=\cD$, and then replace the bound on $k$ in~\eqref{eq:after dim reduction bound doubling} by $k\lesssim \log |\cD|$. 
\end{proof}

\begin{remark}\label{Rem:sep of neighborhood optimal}{\em  An  optimal version of~\eqref{eq:sep of neighborhood of doubling bad bound} is the following statement. There is a universal constant $C>0$ such that if  $(\MM,d_\MM)$ is a metric space and  $\cD\subset \MM$ is complete and $\ul$-doubling for some $\ul\ge 2$, then:
\begin{equation}\label{eq:sharp sep of neighborhood}
\forall r>0,\ \forall T\ge Cr,\qquad \SEP_T\big(B_{\MM}(\cD,r)\big)\lesssim \log \ul. 
\end{equation}
If $\cD$ is compact, then~\eqref{eq:sharp sep of neighborhood} follows from~\cite[Theorem~3.17]{LN05} applied to any $\ul^{O(1)}$-doubling nondegenerate measure $\sigma$ on $\cD$; the existence of such a measure is due to~\cite{VK84}, though the aforementioned dependence of its doubling constant on $\ul$ is not stated there but instead follows from an inspection of its proof (see also~\cite[Chapter~13]{Hei01}). If  $\cD$ is complete but not necessarily compact, then~\eqref{eq:sharp sep of neighborhood}  is not proved in the literature but we will next indicate two ways to justify it. One such route is to use~\cite[Corollary~3.12]{LN05} to obtain a $O(\log \ul)$-padded partition of $\cD$, then use~\cite[Lemma~3.8]{LN05} to extend it to a $O(\log \ul)$-padded random partition of  $B_{\MM}(\cD,r)$, and finally use~\cite[Theorem~2.2]{LN03} to transform the extended random partition to a random partition of $B_{\MM}(\cD,r)$ that is $O(\log \ul)$-separating.  However, \cite{LN03} is an unpublished manuscript and one needs to justify the measurability that we need of each of the above steps (we checked that this is indeed the case, but it is somewhat tedious and does not appear in the literature). An alternative way to prove~\eqref{eq:sharp sep of neighborhood} is to follow the approach in the proof of~\cite[Theorem~3.17]{LN05} when $\sigma$ is now the    $\ul^{O(1)}$-doubling nondegenerate measure $\sigma$ on $\cD$ whose existence is proved in~\cite{LS98} (since $\cD$ is complete). As this $\sigma$ can have infinite mass (unlike the compact case), one needs to incorporate a (somewhat involved) reasoning that is a suitable adaptation of what was done in~\cite[Chapter~4]{naor2024extension}; an upshot of this route is that the measurability that we need is already worked out in~\cite{naor2024extension}.   Both of the above ways to demonstrate~\eqref{eq:sharp sep of neighborhood} result in digressions that are lengthier than how we proceeded above to establish~\eqref{eq:sep of neighborhood of doubling bad bound}, albeit suboptimally.}
\end{remark}

\section{Deduction  of Theorem~\ref{thm:ext from sep of neighborhood} from~\texorpdfstring{\cite{LN05}}{LN05}}\label{sec:exte neighborhood version proof}

Here will explain how to  derive Theorem~\ref{thm:ext from sep of neighborhood}  from one of the main results of~\cite{LN05}.

\begin{proof}[Proof of Theorem~\ref{thm:ext from sep of neighborhood}] Fix $L>0$ and denote
$$
\sigma\eqdef \sup_{\Delta>0} \SEP_\Delta\Big(B_\MM\big(\sub,\frac{1}{L}\Delta\big)\Big). 
$$
We may assume that $\sigma<\infty$ as otherwise~\eqref{eq:lip extension from sep of neighborhood} is vacuous. Fix $\Delta>0$. By the definition of $\sigma$, for any $\e>0$ there exists a probability space $(\Omega,\prob)=(\Omega_{\e,\Delta,L},\prob_{\e,\Delta,L})$ and a random partition 
$$
\Part\eqdef \Part_{\e,\Delta,L}=\left\{\Gamma^i:\Omega\to 2^{B_\MM\left(\sub,\frac{1}{L}\Delta\right)}\right\}_{i=1}^\infty
$$
of $B_\MM(\sub,\Delta/L)$ that is $\Delta$-bounded and $(\sigma+\e)$-separating. Because $\sub$ is nonempty and locally compact, by~\cite[Lemma~115]{naor2024extension} for each $i\in \N$ there exists a $\prob$-to-Borel measurable function $\gamma^i:\Omega\to \sub$ that satisfies 
\begin{equation}\label{eq:measurable nearest point}
\forall \omega\in \Omega,\qquad \Gamma^i(\omega)\neq\emptyset \implies d_\MM\big(\gamma^i(\omega),\Gamma^i(\omega)\big)=d_\MM\big(\sub,\Gamma^i(\omega)\big).
\end{equation}

For each $x\in \MM\setminus B_\MM(\sub,\Delta/L)$ fix any  $c_x\in \sub$ with $d_\MM(x,c_x)<2d_\MM(x,\sub)$. Denote by  $\Gamma^x:\Omega\to 2^\MM$ and $\gamma^x:\Omega\to \sub$ the constant (set-valued and point-valued, respectively) functions that are given by setting $\Gamma^x(\omega)=\{x\}$ and $\gamma^x(\omega)=c_x$ for every $\omega\in \Omega$. Then, the following forms a stochastic decomposition\footnote{Formally, Definition~3.1 of~\cite{LN05} also requires that $\Gamma^i(\omega)$ is a Borel subset of $\MM$ for every $i\in \N$ and every $\omega\in \Omega$. However, this assumption is never used in~\cite{LN05}. To see this, note that such measurability is not required in the definition of gentle partition of unity in~\cite[Section~2]{LN05}. The gentle partition of unity that is constructed in~\cite{LN05} from the  stochastic decomposition is given in~\cite[equation~(6)]{LN05}. For it to be obey the measurability requirements of~\cite[Section~2]{LN05} one needs that for each fixed $x\in \MM\setminus \sub$ the function $(\omega \in \Omega)\mapsto \1_{\Gamma^i(\omega)}(x)\in\{0,1\}$ is $\prob$-measurable, which is equivalent to the set $\{\omega\in \Omega:\ x\in \Gamma^i(\omega)\}$ being $\prob$-measurable, and this is a special case of the strong measurability which is imposed by the definition of random partition that is used in the present work (following~\cite{naor2024extension}). That said, all of the random partitions that we use herein are into Borel subsets of the given metric space, as is evident by their constructions (including those that are cited, specifically from~\cite{naor2024extension}, for Theorem~\ref{thm:sep for lp spaces}).} of $\MM$ with respect to $\sub$ in the sense of~\cite[Definition~3.1]{LN05}:
\begin{equation}\label{eq:stochastic}
\Big(\Omega, \prob, \{\Gamma^i(\cdot),\gamma^i(\cdot)\}_{i=1}^\infty \cup \{\Gamma^x(\cdot),\gamma^x(\cdot)\}_{x\in \MM\setminus B_\MM\left(\sub,\frac{1}{L}\Delta\right)}\Big). 
\end{equation}

By combining~\cite[Lemma~2.1]{LN05} and~\cite[Theorem~4.1, part~{\em 3.}]{LN05} it suffices to check that the stochastic decomposition~\eqref{eq:stochastic} is $\Delta$-bounded and $(1/(2L),1/\max\{2L,\sigma+\e\})$-separating with respect to $\sub$, in the sense of~\cite[Definition~3.2]{LN05} and~\cite[Definition~3.7]{LN05}, respectively. The former requirement is immediate as~\eqref{eq:stochastic}  adds singleton clusters to the random partition $\Part$. The latter requirements mean that   if we write
\begin{equation}\label{eq:P with singletons far away}
\forall \omega\in \Omega,\qquad \mathcal{Q}^\omega\eqdef \Part^\omega\cup \big\{\{x\}\big\}_{x\in \MM\setminus B_\MM\left(\sub,\frac{1}{L}\Delta\right)},
\end{equation}
then
\begin{equation}\label{eq:adding singletons keeps separating}
\forall x,y\in \MM,\qquad d_\MM(\{x,y\},\sub) \le \frac{1}{2L}\Delta\implies \prob  \big[\mathcal{Q}^\omega(x)\neq \mathcal{Q}^\omega(y)\big] \le \frac{\max\{2L,\sigma+\e\}}{\Delta}d_\MM(x,y). 
\end{equation}
To verify~\eqref{eq:adding singletons keeps separating}, fix $x,y\in \MM$ with $d_\MM(\{x,y\},\sub)\le \Delta/(2L)$. Assume that also $d_\MM(x,y)<\Delta/(2L)$, as otherwise the right hand side of~\eqref{eq:adding singletons keeps separating} is at least $1$. So,
\[
\max\{d_\MM(x,\sub),d_\MM(y,\sub)\}\le  d_\MM(\{x,y\},\sub)+d_\MM(x,y)\le \frac{\Delta}{L},
\]
i.e., $x,y\in B_\MM(\sub,\Delta/L)$, whence~\eqref{eq:adding singletons keeps separating} follows from the definition~\eqref{eq:P with singletons far away} of $\mathcal{Q}$ and the fact that $\Part$ is assumed to be a $(\sigma+\e)$-separating random partition of $B_\MM(\sub,\Delta/L)$. 
\end{proof}

\begin{remark} {\em We assumed in Theorem~\ref{thm:ext from sep of neighborhood}  that $\sub$ is locally compact only for invoking~\cite[Lemma~115]{naor2024extension} to get measurable functions $\{\gamma^i:\Omega\to \sub\}_{i=1}^\infty$ that satisfy~\eqref{eq:measurable nearest point} (an inspection of the proof of~\cite[Lemma~115]{naor2024extension} reveals that all that is required for this is that $\sub$ is $\sigma$-compact). One can alternatively stipulate the existence of such functions (even with the weaker requirement $d_\MM(\gamma^i(\omega),\Gamma^i(\omega))<2d_\MM(\sub,\Gamma^i(\omega))$ in~\eqref{eq:measurable nearest point}) as part of the assumption of Theorem~\ref{thm:ext from sep of neighborhood}. This would be analogous to the route that was pursued in~\cite{LN05}, which axiomatizes the minimal assumptions that are needed for its proofs. However, it is hard to maintain such an assumption under various natural operations, including those that are applied herein, such as preimages under Lipschitz functions and intersections. The foundational reworking of~\cite{naor2024extension} circumvents such issues by imposing stronger measurability requirements that are readily seen to be  perserved under natural geometric operations; adopting this approach facilitates the constructions herein. }
\end{remark}

%\subsection*{Acknowledgements} A.~N.~was supported  by NSF grant DMS-2453936, BSF grant 2018223, and a Simons Investigator award. K.~R.~was supported by %an NSF GRFP fellowship.
%%      --------------------------- BIBLIOGRAPHY ----------------------------
%%      ---------------------------------------------------------------------
%% PUT HERE THE BIBLIOGRAPHY IN YOUR FAVOURITE FORMAT
%% Please check that the format of the bibliography is uniform and coherent

\bibliographystyle{abbrv}
\bibliography{distortion-separation-growth-ars}

\def\cprime{$'$}
\begin{thebibliography}{10}

\bibitem{ABN15}
I.~Abraham, Y.~Bartal, and O.~Neiman.
\newblock Local embeddings of metric spaces.
\newblock {\em Algorithmica}, 72(2):539--606, 2015.

\bibitem{AP10}
A.~B. Aleksandrov and V.~V. Peller.
\newblock Functions of operators under perturbations of class {${\bf S}_p$}.
\newblock {\em J. Funct. Anal.}, 258(11):3675--3724, 2010.

\bibitem{ANN25}
A.~Andoni, A.~Naor, and O.~Neiman.
\newblock On isomorphic dimension reduction in $\ell_1$.
\newblock Preprint, 2025.

\bibitem{ALN07}
S.~Arora, J.~R. Lee, and A.~Naor.
\newblock Fr\'echet embeddings of negative type metrics.
\newblock {\em Discrete Comput. Geom.}, 38(4):726--739, 2007.

\bibitem{ALN08}
S.~Arora, J.~R. Lee, and A.~Naor.
\newblock Euclidean distortion and the sparsest cut.
\newblock {\em J. Amer. Math. Soc.}, 21(1):1--21, 2008.

\bibitem{ARV09}
S.~Arora, S.~Rao, and U.~Vazirani.
\newblock Expander flows, geometric embeddings and graph partitioning.
\newblock {\em J. ACM}, 56(2):Art. 5, 37, 2009.

\bibitem{Ass83}
P.~Assouad.
\newblock Plongements lipschitziens dans {${\bf R}\sp{n}$}.
\newblock {\em Bull. Soc. Math. France}, 111(4):429--448, 1983.

\bibitem{AR98}
Y.~Aumann and Y.~Rabani.
\newblock An {$O(\log k)$} approximate min-cut max-flow theorem and
  approximation algorithm.
\newblock {\em SIAM J. Comput.}, 27(1):291--301, 1998.

\bibitem{Bal92}
K.~Ball.
\newblock Markov chains, {R}iesz transforms and {L}ipschitz maps.
\newblock {\em Geom. Funct. Anal.}, 2(2):137--172, 1992.

\bibitem{Ball-Ribe}
K.~Ball.
\newblock The {R}ibe programme.
\newblock Number 352, pages Exp. No. 1047, viii, 147--159. 2013.
\newblock S\'{e}minaire Bourbaki. Vol. 2011/2012. Expos\'{e}s 1043--1058.

\bibitem{BCL94}
K.~Ball, E.~A. Carlen, and E.~H. Lieb.
\newblock Sharp uniform convexity and smoothness inequalities for trace norms.
\newblock {\em Invent. Math.}, 115(3):463--482, 1994.

\bibitem{bartal1996probabilistic}
Y.~Bartal.
\newblock Probabilistic approximation of metric spaces and its algorithmic
  applications.
\newblock In {\em Proceedings of 37th Conference on Foundations of Computer
  Science}, pages 184--193. IEEE, 1996.

\bibitem{BG14}
Y.~Bartal and L.-A. Gottlieb.
\newblock Dimension reduction techniques for $\ell_{\!\!p}$, $1\le p<\infty$,
  with applications.
\newblock Preprint available at \url{https://arxiv.org/pdf/1408.1789v2}
  (observe that this is version v2), 2014.

\bibitem{Beg99}
B.~Begun.
\newblock A remark on almost extensions of {L}ipschitz functions.
\newblock {\em Israel J. Math.}, 109:151--155, 1999.

\bibitem{BL00}
Y.~Benyamini and J.~Lindenstrauss.
\newblock {\em Geometric nonlinear functional analysis. {V}ol. 1}, volume~48 of
  {\em American Mathematical Society Colloquium Publications}.
\newblock American Mathematical Society, Providence, RI, 2000.

\bibitem{bourgain1985lipschitz}
J.~Bourgain.
\newblock On {L}ipschitz embedding of finite metric spaces in {H}ilbert space.
\newblock {\em Israel Journal of Mathematics}, 52:46--52, 1985.

\bibitem{Bourgain-superreflexivity}
J.~Bourgain.
\newblock The metrical interpretation of superreflexivity in {B}anach spaces.
\newblock {\em Israel J. Math.}, 56(2):222--230, 1986.

\bibitem{Bou87}
J.~Bourgain.
\newblock Remarks on the extension of {L}ipschitz maps defined on discrete sets
  and uniform homeomorphisms.
\newblock In {\em Geometrical aspects of functional analysis (1985/86)}, volume
  1267 of {\em Lecture Notes in Math.}, pages 157--167. Springer, Berlin, 1987.

\bibitem{BN25}
M.~Braverman and A.~Naor.
\newblock Quantitative {W}asserstein rounding.
\newblock Preprint, 2025.

\bibitem{CL93}
E.~A. Carlen and E.~H. Lieb.
\newblock Optimal hypercontractivity for {F}ermi fields and related
  noncommutative integration inequalities.
\newblock {\em Comm. Math. Phys.}, 155(1):27--46, 1993.

\bibitem{Cha95}
F.~Chaatit.
\newblock On uniform homeomorphisms of the unit spheres of certain {B}anach
  lattices.
\newblock {\em Pacific J. Math.}, 168(1):11--31, 1995.

\bibitem{ANR24}
A.~Chang, A.~Naor, and K.~Ren.
\newblock Random zero sets with local growth guarantees.
\newblock Preprint available at \url{https://arxiv.org/pdf/2410.21931}, 2024.

\bibitem{ccg98}
M.~Charikar, C.~Chekuri, A.~Goel, S.~Guha, and S.~Plotkin.
\newblock Approximating a finite metric by a small number of tree metrics.
\newblock In {\em Proceedings 39th Annual Symposium on Foundations of Computer
  Science (Cat. No.98CB36280)}, pages 379--388, 1998.

\bibitem{CLLM10}
M.~Charikar, T.~Leighton, S.~Li, and A.~Moitra.
\newblock Vertex sparsifiers and abstract rounding algorithms.
\newblock In {\em 2010 {IEEE} 51st {A}nnual {S}ymposium on {F}oundations of
  {C}omputer {S}cience---{FOCS} 2010}, pages 265--274. IEEE Computer Soc., Los
  Alamitos, CA, 2010.

\bibitem{CMM10}
M.~Charikar, K.~Makarychev, and Y.~Makarychev.
\newblock Local global tradeoffs in metric embeddings.
\newblock {\em SIAM J. Comput.}, 39(6):2487--2512, 2010.

\bibitem{Cor24}
W.~Corr\^ea.
\newblock Uniform homeomorphisms between spheres induced by interpolation
  methods.
\newblock {\em Proc. Amer. Math. Soc.}, 152(5):2157--2167, 2024.

\bibitem{CFGT25}
W.~Corr\^ea, V.~Ferenczi, R.~Gesing, and P.~Tradacete.
\newblock Extremes of interpolation scales of {B}anach spaces.
\newblock {\em J. Funct. Anal.}, 289(2):Paper No. 110924, 41, 2025.

\bibitem{Dah95}
M.~Daher.
\newblock Hom\'eomorphismes uniformes entre les sph\`eres unit\'e{} des espaces
  d'interpolation.
\newblock {\em Canad. Math. Bull.}, 38(3):286--294, 1995.

\bibitem{DLP13}
J.~Ding, J.~R. Lee, and Y.~Peres.
\newblock Markov type and threshold embeddings.
\newblock {\em Geom. Funct. Anal.}, 23(4):1207--1229, 2013.

\bibitem{Dvo60}
A.~Dvoretzky.
\newblock Some results on convex bodies and {B}anach spaces.
\newblock In {\em Proc. {I}nternat. {S}ympos. {L}inear {S}paces ({J}erusalem,
  1960)}, pages 123--160. Jerusalem Academic Press, Jerusalem, 1961.

\bibitem{Enf69}
P.~Enflo.
\newblock On the nonexistence of uniform homeomorphisms between
  {$L\sb{p}$}-spaces.
\newblock {\em Ark. Mat.}, 8:103--105, 1969.

\bibitem{EGKRTT14}
M.~Englert, A.~Gupta, R.~Krauthgamer, H.~R\"acke, I.~Talgam-Cohen, and
  K.~Talwar.
\newblock Vertex sparsifiers: new results from old techniques.
\newblock {\em SIAM J. Comput.}, 43(4):1239--1262, 2014.

\bibitem{EMN19}
A.~Eskenazis, M.~Mendel, and A.~Naor.
\newblock Nonpositive curvature is not coarsely universal.
\newblock {\em Invent. Math.}, 217(3):833--886, 2019.

\bibitem{God17}
G.~Godefroy.
\newblock De {G}rothendieck \`a {N}aor: une promenade dans l'analyse
  m\'{e}trique des espaces de {B}anach.
\newblock {\em Gaz. Math.}, (151):13--24, 2017.

\bibitem{Gru60}
B.~Gr\"unbaum.
\newblock Some applications of expansion constants.
\newblock {\em Pacific J. Math.}, 10:193--201, 1960.

\bibitem{GKL03}
A.~Gupta, R.~Krauthgamer, and J.~R. Lee.
\newblock Bounded geometries, fractals, and low-distortion embeddings.
\newblock In {\em 44th Symposium on Foundations of Computer Science {(FOCS}
  2003), 11-14 October 2003, Cambridge, MA, USA, Proceedings}, pages 534--543.
  {IEEE} Computer Society, 2003.

\bibitem{Hei01}
J.~Heinonen.
\newblock {\em Lectures on analysis on metric spaces}.
\newblock Universitext. Springer-Verlag, New York, 2001.

\bibitem{HK15}
J.~Heinonen, P.~Koskela, N.~Shanmugalingam, and J.~T. Tyson.
\newblock {\em Sobolev spaces on metric measure spaces}, volume~27 of {\em New
  Mathematical Monographs}.
\newblock Cambridge University Press, Cambridge, 2015.
\newblock An approach based on upper gradients.

\bibitem{Him75}
C.~J. Himmelberg.
\newblock Measurable relations.
\newblock {\em Fund. Math.}, 87:53--72, 1975.

\bibitem{HSZ22}
J.~Huang, F.~Sukochev, and D.~Zanin.
\newblock Operator {$\theta$}-{H}\"older functions with respect to
  {$\Vert\cdot\Vert_p$}, {$0 < p \leqslant \infty$}.
\newblock {\em J. Lond. Math. Soc. (2)}, 105(4):2436--2477, 2022.

\bibitem{JL82}
W.~B. Johnson and J.~Lindenstrauss.
\newblock Extensions of {L}ipschitz mappings into a {H}ilbert space.
\newblock In {\em Conference in modern analysis and probability ({N}ew {H}aven,
  {C}onn., 1982)}, volume~26 of {\em Contemp. Math.}, pages 189--206. Amer.
  Math. Soc., Providence, RI, 1984.

\bibitem{JLS86}
W.~B. Johnson, J.~Lindenstrauss, and G.~Schechtman.
\newblock Extensions of {L}ipschitz maps into {B}anach spaces.
\newblock {\em Israel J. Math.}, 54(2):129--138, 1986.

\bibitem{JS01}
W.~B. Johnson and G.~Schechtman.
\newblock Finite dimensional subspaces of {$L_p$}.
\newblock In {\em Handbook of the geometry of {B}anach spaces, {V}ol. {I}},
  pages 837--870. North-Holland, Amsterdam, 2001.

\bibitem{Kal08}
N.~J. Kalton.
\newblock The nonlinear geometry of {B}anach spaces.
\newblock {\em Rev. Mat. Complut.}, 21(1):7--60, 2008.

\bibitem{Kec95}
A.~S. Kechris.
\newblock {\em Classical descriptive set theory}, volume 156 of {\em Graduate
  Texts in Mathematics}.
\newblock Springer-Verlag, New York, 1995.

\bibitem{Kir34}
M.~Kirszbraun.
\newblock {\"U}ber die zusammenziehende und lipschitzsche transformationen.
\newblock {\em Fundamenta Mathematicae}, 22(1):77--108, 1934.

\bibitem{KLMN05}
R.~Krauthgamer, J.~R. Lee, M.~Mendel, and A.~Naor.
\newblock Measured descent: a new embedding method for finite metrics.
\newblock {\em Geom. Funct. Anal.}, 15(4):839--858, 2005.

\bibitem{krauthgamer2025lipschitz}
R.~Krauthgamer and N.~Petruschka.
\newblock Lipschitz decompositions of finite $\ell_p$ metrics.
\newblock {\em Preprint available at \url{https://arxiv.org/pdf/2502.01120}},
  2025.

\bibitem{KPS25}
R.~Krauthgamer, N.~Petruschka, and S.~Sapir.
\newblock The {P}ower of {R}ecursive {E}mbeddings for $\ell_{\!\! p}$
  {M}etrics.
\newblock Preprint available at \url{https://arxiv.org/pdf/2503.18508}, 2025.

\bibitem{Lee05}
J.~R. Lee.
\newblock On distance scales, embeddings, and efficient relaxations of the cut
  cone.
\newblock In {\em Proceedings of the {S}ixteenth {A}nnual {ACM}-{SIAM}
  {S}ymposium on {D}iscrete {A}lgorithms}, pages 92--101. ACM, New York, 2005.

\bibitem{LN03}
J.~R. Lee and A.~Naor.
\newblock Metric decomposition, smooth measures, and clustering.
\newblock Unpublished manuscript, available on request, 2003.

\bibitem{LN04-comptes}
J.~R. Lee and A.~Naor.
\newblock Absolute {L}ipschitz extendability.
\newblock {\em C. R. Math. Acad. Sci. Paris}, 338(11):859--862, 2004.

\bibitem{LN05}
J.~R. Lee and A.~Naor.
\newblock Extending {L}ipschitz functions via random metric partitions.
\newblock {\em Invent. Math.}, 160(1):59--95, 2005.

\bibitem{Lew78}
D.~R. Lewis.
\newblock Finite dimensional subspaces of {$L\sb{p}$}.
\newblock {\em Studia Math.}, 63(2):207--212, 1978.

\bibitem{Lin64}
J.~Lindenstrauss.
\newblock On nonlinear projections in {B}anach spaces.
\newblock {\em Michigan Math. J.}, 11:263--287, 1964.

\bibitem{LT77}
J.~Lindenstrauss and L.~Tzafriri.
\newblock {\em Classical {B}anach spaces. {I}}, volume Band 92 of {\em
  Ergebnisse der Mathematik und ihrer Grenzgebiete [Results in Mathematics and
  Related Areas]}.
\newblock Springer-Verlag, Berlin-New York, 1977.
\newblock Sequence spaces.

\bibitem{LT79}
J.~Lindenstrauss and L.~Tzafriri.
\newblock {\em Classical {B}anach spaces. {II}}, volume~97 of {\em Ergebnisse
  der Mathematik und ihrer Grenzgebiete [Results in Mathematics and Related
  Areas]}.
\newblock Springer-Verlag, Berlin-New York, 1979.
\newblock Function spaces.

\bibitem{linial1995geometry}
N.~Linial, E.~London, and Y.~Rabinovich.
\newblock The geometry of graphs and some of its algorithmic applications.
\newblock {\em Combinatorica}, 15:215--245, 1995.

\bibitem{Lor01}
G.~G. Lorentz.
\newblock Who discovered analytic sets?
\newblock {\em Math. Intelligencer}, 23(4):28--32, 2001.

\bibitem{LS98}
J.~Luukkainen and E.~Saksman.
\newblock Every complete doubling metric space carries a doubling measure.
\newblock {\em Proc. Amer. Math. Soc.}, 126(2):531--534, 1998.

\bibitem{Luz17}
N.~N. Luzin.
\newblock Sur la classification de {M}. {B}aire.
\newblock {\em C. R. Math. Acad. Sci. Paris}, 164:91--94, 1917.

\bibitem{MM16}
K.~Makarychev and Y.~Makarychev.
\newblock Metric extension operators, vertex sparsifiers and {L}ipschitz
  extendability.
\newblock {\em Israel J. Math.}, 212(2):913--959, 2016.

\bibitem{MP84}
M.~B. Marcus and G.~Pisier.
\newblock Characterizations of almost surely continuous {$p$}-stable random
  {F}ourier series and strongly stationary processes.
\newblock {\em Acta Math.}, 152(3-4):245--301, 1984.

\bibitem{MAT90}
J.~Matou\v{s}ek.
\newblock Extension of {L}ipschitz mappings on metric trees.
\newblock {\em Comment. Math. Univ. Carolin.}, 31(1):99--104, 1990.

\bibitem{Mat97}
J.~Matou\v{s}ek.
\newblock On embedding expanders into {$l_p$} spaces.
\newblock {\em Israel J. Math.}, 102:189--197, 1997.

\bibitem{MP76}
B.~Maurey and G.~Pisier.
\newblock S\'eries de variables al\'eatoires vectorielles ind\'ependantes et
  propri\'et\'es g\'eom\'etriques des espaces de {B}anach.
\newblock {\em Studia Math.}, 58(1):45--90, 1976.

\bibitem{mazur1929remarque}
S.~Mazur.
\newblock Une remarque sur l'hom{\'e}omorphie des champs fonctionnels.
\newblock {\em Studia Mathematica}, 1(1):83--85, 1929.

\bibitem{mendel2004euclidean}
M.~Mendel and A.~Naor.
\newblock Euclidean quotients of finite metric spaces.
\newblock {\em Advances in Mathematics}, 189(2):451--494, 2004.

\bibitem{MN06}
M.~Mendel and A.~Naor.
\newblock Some applications of {B}all's extension theorem.
\newblock {\em Proc. Amer. Math. Soc.}, 134(9):2577--2584, 2006.

\bibitem{MN08}
M.~Mendel and A.~Naor.
\newblock Metric cotype.
\newblock {\em Ann. of Math. (2)}, 168(1):247--298, 2008.

\bibitem{MN13-ext}
M.~Mendel and A.~Naor.
\newblock Spectral calculus and {L}ipschitz extension for barycentric metric
  spaces.
\newblock {\em Anal. Geom. Metr. Spaces}, 1:163--199, 2013.

\bibitem{Moi09}
A.~Moitra.
\newblock Approximation algorithms for multicommodity-type problems with
  guarantees independent of the graph size.
\newblock In {\em 2009 50th {A}nnual {IEEE} {S}ymposium on {F}oundations of
  {C}omputer {S}cience---{FOCS} 2009}, pages 3--12. IEEE Computer Soc., Los
  Alamitos, CA, 2009.

\bibitem{Nao01}
A.~Naor.
\newblock A phase transition phenomenon between the isometric and isomorphic
  extension problems for {H}\"older functions between {$L_p$} spaces.
\newblock {\em Mathematika}, 48(1-2):253--271, 2001.

\bibitem{Naor-Ribe}
A.~Naor.
\newblock An introduction to the {R}ibe program.
\newblock {\em Japanese Journal of Mathematics}, 7(2):167--233, 2012.

\bibitem{naor2014comparison}
A.~Naor.
\newblock Comparison of metric spectral gaps.
\newblock {\em Analysis and Geometry in Metric Spaces}, 2(1):1--52, 2014.

\bibitem{Nao15-nonextend}
A.~Naor.
\newblock Uniform nonextendability from nets.
\newblock {\em C. R. Math. Acad. Sci. Paris}, 353(11):991--994, 2015.

\bibitem{naor2017probabilistic}
A.~Naor.
\newblock Probabilistic clustering of high dimensional norms.
\newblock In {\em Proceedings of the Twenty-Eighth Annual ACM-SIAM Symposium on
  Discrete Algorithms}, pages 690--709. SIAM, 2017.

\bibitem{Nao18}
A.~Naor.
\newblock Metric dimension reduction: a snapshot of the {R}ibe program.
\newblock In {\em Proceedings of the {I}nternational {C}ongress of
  {M}athematicians---{R}io de {J}aneiro 2018. {V}ol. {I}. {P}lenary lectures},
  pages 759--837. World Sci. Publ., Hackensack, NJ, 2018.

\bibitem{Nao21-almost}
A.~Naor.
\newblock Impossibility of almost extension.
\newblock {\em Adv. Math.}, 384:Paper No. 107761, 34, 2021.

\bibitem{naor2024extension}
A.~Naor.
\newblock {\em Extension, separation and isomorphic reverse isoperimetry},
  volume~11 of {\em Memoirs of the European Mathematical Society}.
\newblock European Mathematical Society (EMS), Berlin, 2024.

\bibitem{NPSS06}
A.~Naor, Y.~Peres, O.~Schramm, and S.~Sheffield.
\newblock Markov chains in smooth {B}anach spaces and {G}romov-hyperbolic
  metric spaces.
\newblock {\em Duke Math. J.}, 134(1):165--197, 2006.

\bibitem{NR17}
A.~Naor and Y.~Rabani.
\newblock On {L}ipschitz extension from finite subsets.
\newblock {\em Israel J. Math.}, 219(1):115--161, 2017.

\bibitem{NR25-v1}
A.~Naor and K.~Ren.
\newblock $\ell_{\!\! p}$ has nontrivial {E}uclidean distortion growth when
  $2<p<4$.
\newblock Preprint available at \url{https://arxiv.org/pdf/2410.21931v1}
  (observe that this is version v1), 2025.

\bibitem{NS16}
A.~Naor and G.~Schechtman.
\newblock Metric {$X_p$} inequalities.
\newblock {\em Forum Math. Pi}, 4:e3, 81, 2016.

\bibitem{NS25}
A.~Naor and G.~Schechtman.
\newblock Lipschitz almost-extension and nonexistence of uniform embeddings of
  balls in {S}chatten classes.
\newblock Preprint, 2025.

\bibitem{NS11}
A.~Naor and L.~Silberman.
\newblock Poincar\'e{} inequalities, embeddings, and wild groups.
\newblock {\em Compos. Math.}, 147(5):1546--1572, 2011.

\bibitem{NR03}
I.~Newman and Y.~Rabinovich.
\newblock A lower bound on the distortion of embedding planar metrics into
  {E}uclidean space.
\newblock {\em Discrete Comput. Geom.}, 29(1):77--81, 2003.

\bibitem{OS94}
E.~Odell and T.~Schlumprecht.
\newblock The distortion problem.
\newblock {\em Acta Math.}, 173(2):259--281, 1994.

\bibitem{Ost16}
M.~Ostrovskii.
\newblock Metric characterizations of some classes of {B}anach spaces.
\newblock In {\em Harmonic analysis, partial differential equations, complex
  analysis, {B}anach spaces, and operator theory. {V}ol. 1}, volume~4 of {\em
  Assoc. Women Math. Ser.}, pages 307--347. Springer, [Cham], 2016.

\bibitem{Ostrovskii-book}
M.~I. Ostrovskii.
\newblock {\em Metric embeddings}, volume~49 of {\em De Gruyter Studies in
  Mathematics}.
\newblock De Gruyter, Berlin, 2013.
\newblock Bilipschitz and coarse embeddings into Banach spaces.

\bibitem{Rao99}
S.~Rao.
\newblock Small distortion and volume preserving embeddings for planar and
  {E}uclidean metrics.
\newblock In {\em Proceedings of the {F}ifteenth {A}nnual {S}ymposium on
  {C}omputational {G}eometry ({M}iami {B}each, {FL}, 1999)}, pages 300--306.
  ACM, New York, 1999.

\bibitem{Ray02}
Y.~Raynaud.
\newblock On ultrapowers of non commutative {$L_p$} spaces.
\newblock {\em J. Operator Theory}, 48(1):41--68, 2002.

\bibitem{Ric15}
E.~Ricard.
\newblock H\"older estimates for the noncommutative {M}azur maps.
\newblock {\em Arch. Math. (Basel)}, 104(1):37--45, 2015.

\bibitem{Sob41}
A.~Sobczyk.
\newblock Projections in {M}inkowski and {B}anach spaces.
\newblock {\em Duke Math. J.}, 8:78--106, 1941.

\bibitem{Tom89}
N.~Tomczak-Jaegermann.
\newblock {\em Banach-{M}azur distances and finite-dimensional operator
  ideals}, volume~38 of {\em Pitman Monographs and Surveys in Pure and Applied
  Mathematics}.
\newblock Longman Scientific \& Technical, Harlow; copublished in the United
  States with John Wiley \& Sons, Inc., New York, 1989.

\bibitem{VK84}
A.~L. Vol'berg and S.~V. Konyagin.
\newblock A homogeneous measure exists on any compactum in {${\bf R}^n$}.
\newblock {\em Dokl. Akad. Nauk SSSR}, 278(4):783--786, 1984.

\bibitem{Woj91}
P.~Wojtaszczyk.
\newblock {\em Banach spaces for analysts}, volume~25 of {\em Cambridge Studies
  in Advanced Mathematics}.
\newblock Cambridge University Press, Cambridge, 1991.

\end{thebibliography}

\end{document}